\documentclass[fleqn]{elsarticle}
\usepackage{amsmath,amsfonts,graphicx,epstopdf,float,amssymb,bm,amsthm,multirow,appendix,bbm}

\usepackage{geometry}
\numberwithin{equation}{section}
\newtheorem{thm}{Theorem}[section]

\newtheorem{rmk}{Remark}[section]
\newtheorem{lem}{Lemma}[section]

\newtheorem*{prf}{Proof}
\graphicspath{{FIGS/}}

\begin{document}

\begin{frontmatter}

\title{Structure-preserving algorithms for the two-dimensional sine-Gordon equation with Neumann boundary conditions}

\author{Wenjun Cai$^{a}$, Chaolong Jiang$^{b}$, Yushun Wang$^{a,*}$}
\address[1]{Jiangsu Key Laboratory for NSLSCS, School of Mathematical Sciences, Nanjing Normal University}
\address[2]{School of Statistics and Mathematics, Yunnan University of Finance and Economics}

\begin{abstract}
	
	This paper presents two kinds of strategies to construct structure-preserving algorithms with homogeneous Neumann boundary conditions for the sine-Gordon equation, while most existing structure-preserving algorithms are only valid for zero or periodic boundary conditions. The first strategy is based on the conventional second-order central difference quotient but with a cell-centered grid, while the other is established on the regular grid but incorporated with summation by parts (SBP) operators. Both the methodologies can provide conservative semi-discretizations with different forms of Hamiltonian structures and the discrete energy. However, utilizing the existing SBP formulas, schemes obtained by the second strategy can directly achieve higher-order accuracy while it is not obvious for schemes  based on the cell-centered grid to make accuracy improved easily. Further combining the symplectic Runge-Kutta method and the scalar auxiliary variable (SAV) approach, we construct symplectic integrators and linearly implicit energy-preserving schemes for the two dimensional sine-Gordon equation, respectively. Extensive numerical experiments demonstrate their effectiveness with the homogeneous Neumann boundary conditions.

\end{abstract}


\end{frontmatter}

\begin{figure}[b]
	\small \baselineskip=10pt
	\rule[2mm]{1.8cm}{0.2mm} \par
	$^{*}$Corresponding author.\\
	E-mail address: wangyushun@njnu.edu.cn (Y. Wang).
\end{figure}

\section{Introduction}

In this paper, we consider the undamped two-dimensional sine-Gordon equation
\begin{equation}\label{eq-1-1}
u_{tt}=\Delta u-\phi(x,y)\sin u, \quad (x,y)\in \Omega\subset \mathbb{R}^2
\end{equation}
with initial conditions 
\begin{equation}\label{eq-1-2}
u(x,y,0)=f_1(x,y),\quad u_t(x,y,0)=f_2(x,y),
\end{equation}
and homogeneous Neumann boundary conditions along the boundary $\Gamma$ of $\Omega$,
\begin{equation}\label{eq-1-3}
\nabla u\cdot\bm n=0,\quad (x,y)\in\Gamma.
\end{equation}
The function $\phi$ can be interpreted as a Josephson current density, $f_1$ and $f_2$ are wave modes or kinks and their velocity, respectively. The sine-Gordon equation \eqref{eq-1-1} can be the governing equation of wide variety of physical processes, in particular the propagation of magnetic flux in a Josephson junction consisting of two layer of superconducting material separated by an isolating barrier \cite{j65}.

Conservation of the total energy is the main property of the undamped sine-Gordon system \eqref{eq-1-1}-\eqref{eq-1-3}. Taking the inner product of \eqref{eq-1-1} with $u_t$ and utilizing the homogeneous Neumann boundary condition, we have
\begin{equation}\label{eq-1-4}
\frac{d}{dt}\mathcal{H}(t)=0, \quad \mbox{with}\quad \mathcal{H}(t)=\frac{1}{2}\iint_\Omega \big(u_t^2+|\nabla u|^2+2\phi(1-\cos(u))\big) dxdy.
\end{equation}
Besides the conservation of energy, the sine-Gordon system also possesses a canonical Hamiltonian structure. By introducing the velocity $v=u_t$ and further denoting $z=(u,v)^T$, we can rewrite the system as
\begin{equation}\label{eq-1-5}
z_t=J\frac{\delta\mathcal{H}}{\delta z},\quad J=\left(\begin{array}{cc}
0 & -1 \\ 
1 & 0
\end{array} \right),
\end{equation}
where $\frac{\delta H}{\delta z}$ is the vector of variational derivatives with respect to $u$ and $v$. An important character associated with the Hamiltonian structure \eqref{eq-1-5} is the conservation of sympleticity, which for the sine-Gordon equation reads
\begin{equation}
\frac{d}{dt}\omega(t)=0,\quad\mbox{with}\quad \omega(t)=\iint_\Omega (du\wedge dv) dxdy,
\end{equation}
where $\wedge$ is the wedge product in exterior algebra, and $du, dv$ are solutions of variation equations of the system \eqref{eq-1-5}. For details, we refer readers to \cite{m94,br01}.

When designing numerical algorithms for Hamiltonian partial differential equations (PDEs), it is naturally to preserve the conservative structures of the given continuous system in discrete analogs.  As two fundamental structure-preserving algorithms (SPAs), the symplectic and energy-preserving integrators have attracted continuous attention in the recent decades, emerging lots of applications on Hamiltonian PDEs, such as the nonlinear wave equation \cite{m94,ztstw10,lws18}, the  Schr\"odinger-type equation \cite{iks01,khwf09,lmw18}, the KdV equation \cite{am05,cm07}, the Maxwell's equations \cite{cll10,st11,chwg15}, the Camassa-Holm equation \cite{mm12,gw16,csw16} and so on. However, it is noticed that most of the existing SPAs are based on zero Dirichlet or periodic boundary conditions. Although the use of such kinds of boundary conditions can enable the numerical solution to propagate straightly in a finite window, and contribute a long-time verification of numerical conservation of invariants, it cannot simulate all real physical processes where other boundary conditions are requisite. For example, when solving the heat equation, the ideal heat insulator with the heat diffusion has to be described by the homogeneous Neumann boundary conditions. As a consequence, it is desirable to extend the SPAs to more general boundary conditions. In the current work, we consider the two-dimensional sine-Gordon system \eqref{eq-1-1}-\eqref{eq-1-3} and start the investigation from homogeneous Neumann boundary conditions.

There have many efficient numerical schemes concerning the two-dimensional sine-Gordon equation \eqref{eq-1-1} with Neumann boundary conditions on both damped and undamped circumstances, which are based on several different  methodologies, such as the finite difference schemes \cite{cl81,skv05,b07}, the finite element schemes \cite{ahh91}, the meshless methods \cite{ds08,kdt17} and so on. Although these proposed methods can well resolve the Neumann boundary condition,  none of them are strictly SPAs and thereby cannot guarantee a long-time stability as well as the preservation of conservative quantities, for example, the system energy. The only existing results of energy preservation  are obtained for the one-dimensional sine-Gordon equation \cite{slwl17,lws18}, which are based on the finite element method and the finite difference method, respectively.  As we all know, the semi-discretization by the finite element method can naturally preserve discrete analogs of conservation laws, regarding different boundary conditions. However, only when utilizing temporal geometric integrations can complete the construction of fully discrete SPAs. In the study of \cite{slwl17}, authors straightly applied this fact on the one-dimensional sine-Gordon equation and derived an energy-preserving integrator. A similar idea can be found in \cite{czw17} for the simulation of two-dimensional seismic wave equations. While for the finite difference based schemes in \cite{lws18}, authors used Taylor expansions to achieve a consistent accuracy in all spatial grid points for semi-discretizations and further combined generalized AVF methods to obtain energy-preserving schemes with second and fourth orders, respectively. Although one can construct SPAs by extending the finite element semi-discretization \cite{slwl17} to the two-dimensional sine-Gordon equation directly, for cases considered in rectangular domains, it is preferable to derive SPAs in the framework of finite difference method due to its simplicity of implementation. Besides, the Taylor expansions will become more tedious when generalizing the idea in \cite{lws18} to get schemes of higher-order accuracy for two-dimensional problems. Therefore, in the current work, we propose another two strategies to construct SPAs for the sine-Gordon equation \eqref{eq-1-1} with Neumann boundary conditions.

The first strategy is rather simple which is based on a cell-centered grid and a second-order central difference method for the approximation of both spatial derivatives and boundary conditions. The resulting differential matrix for the Laplace operator is automatically symmetric, making the semi-discrete ODE system Hamiltonian. Subsequently, symplectic and energy-preserving integrators can be derived. The second strategy is based on the summation by parts (SBP) formula, as a counterpart of integration by parts, which is a crucial important feature for differential operators in the derivation of energy conservation laws. SBP operators have predominantly been developed in the context of high-order finite difference methods \cite{ks74,s94} where the nodal distribution in computational space is uniform, see the excellent review papers \cite{fhz14,sn14} and the references therein. In this paper, we use the SBP finite difference operators in \cite{mn04} to approximate second derivatives, where SBP operators up to eighth-order accuracy are given explicitly. Moreover, the corresponding semi-discretizations still possess non-canonical Hamiltonian structures, in contrast to the canonical ones with respect to that obtained by the first strategy. Nevertheless, it does not affect the construction of symplectic and energy-preserving integrators, except the definitions of the discrete Hamiltonians and equivalently the discrete energy functions.

When the original PDE system is discretized to a Hamiltonian ODEs, various time integration methods can thereby be applied to get an energy-preserving integrator, such as the discrete gradient method \cite{mqr99}, the averaged vector field (AVF) method \cite{qm08}, the Hamitonian boundary value method (HBVM) \cite{bit12}, and so on. However, it can be noticed that most of the resulting energy-preserving integrators are fully implicit and therefore time-consuming. In this paper, we utilize the idea of the scalar auxiliary variable (SAV) approach to design linearly implicit energy-preserving schemes for the two-dimensional sine-Gordon equation \eqref{eq-1-1}. The SAV approach is recently proposed by Shen et al. \cite{sxy17,sxy18} to deal with nonlinear terms in gradient flows, and leads to efficient and robust energy stable schemes which only require to solve decoupled linear equations with constant coefficients. To the best of our knowledge, there is no result concerning the SAV approach for the conservative Hamiltonian system. Taking the sine-Gordon equation \eqref{eq-1-1} for example, we first explore the feasibility of the SAV approach and then devise several linearly implicit energy-preserving schemes based on the aforementioned two strategies of semi-discretizations, respectively. While for the corresponding symplectic integrators, we simply choose the implicit midpoint rule for the time integration, not only to verify the effectiveness of the conservative semi-discretization but also to compare the efficiency of the SAV based schemes with a typical implicit one, although explicit symplectic integrators can be easily constructed by symplectic partitioned Runge-Kutta methods \cite{hlw06}.

The rest of this paper is organized as follows. In Section 2, we give a conservative semi-discretization based on the cell-centered grid and derive the corresponding discrete energy conservation law as well as a canonical Hamiltonian structure. In Section 3, we utilize the SBP operator to construct an alternative semi-discretization with a modified energy conservation law and a non-canonical Hamiltonian structure. The detailed application of the SAV approach on semi-discrete systems obtained from both the two strategies are discussed in Section 4. Extensive numerical experiments are given in Section 5 to demonstrate the effectiveness of the proposed methods in energy preservation and long-time performance with homogeneous Neumann boundary conditions. In Appendix A, we further discuss the implementation of SBP operator for non-homogeneous Neumann boundary conditions, with which the convergence tests in Section 5 are carried on, and a fourth-order SBP formula is also presented.

\section{Conservative semi-discretization on the cell-centered grid}

In this section, we will give a conservative finite difference semi-discretization on the cell-centered grid, which both satisfies a discrete energy conservation law and a canonical Hamiltonian structure.

 Let $\Omega=[a,b]\times[c,d]$ and $N_x$, $N_y$ be positive integers, the domain $\Omega$ is uniformly partitioned with $h_x=(b-a)/N_x$, $h_y=(d-c)/N_y$ and variables are stored at each cell center as follows:
\[
\Omega_h=\Big\{(x_j,y_k)\big|x_j=a+\big(j-\frac{1}{2}\big)h_x, y_k=c+\big(k-\frac{1}{2}\big)h_y, 1\leq j\leq N_x, 1\leq k\leq N_y\Big\}.
\]
The time interval $(0,T]$ is discretized uniformly by $t^n=n\tau$, $n=0,1,\cdots, N$, where $\tau=T/N$, and $N$ is a positive integer. Let $V_h=\{u=(u_{j,k})|(x_j,y_k)\in\Omega_h\}$ be a space of grid function. For $u, v\in V_h$, we define the following difference operators:
\[
\begin{aligned}
&\delta_x^+u_{j,k}=\frac{u_{j+1,k}-u_{j,k}}{h_x},\quad \delta_x^-u_{j,k}=\frac{u_{j,k}-u_{j-1,k}}{h_x},\\
&\delta_y^+u_{j,k}=\frac{u_{j,k+1}-u_{j,k}}{h_y},\quad \delta_y^-u_{j,k}=\frac{u_{j,k}-u_{j,k-1}}{h_y},\\
&\nabla_h^+u_{j,k}=(\delta_x^+u_{j,k},\delta_y^+u_{j,k})^T,\quad \nabla_h^-u_{j,k}=(\delta_x^-u_{j,k},\delta_y^-u_{j,k})^T,\\
&\Delta_h u_{j,k}=\frac{u_{j+1,k}-2u_{j,k}+u_{j-1,k}}{h_x^2}+\frac{u_{j,k+1}-2u_{j,k}+u_{j,k-1}}{h_y^2}.
\end{aligned}
\]
The corresponding discrete $L^2$ inner product and discrete $L^2$ norm are defined as
\begin{equation}\label{eq-2-0}
\begin{aligned}
&(u,v)_h=\sum_{j=1}^{N_x}\sum_{k=1}^{N_y}u_{j,k}v_{j,k}h_xh_y,\quad \|u\|_h^2=(u,u)_h,\\
&(\nabla_h^+u,\nabla_h^+u)_h=(\delta_x^+u,\delta_x^+u)_h+(\delta_y^+u,\delta_y^+u)_h
\end{aligned}
\end{equation}
Notice that for the homogeneous Neumann boundary conditions \eqref{eq-1-3}, the values on the center of the ghost cells outside the boundary are given by the central difference quotient such that
\begin{equation}\label{eq-2-1}
u_{0,k}=u_{1,k}, \quad u_{N_x+1,k}=u_{N_x,k}, \quad u_{j,0}=u_{j,1}, \quad u_{j,N_y+1}=u_{j,N_y},
\end{equation}
where the ghost points are denoted by the subscripts $0, N_x+1, N_y+1$. For the derivation of structure-preserving algorithms, the symmetry as well as the summation by parts properties of above difference operators are essential. It is well-known that with the periodic boundary condition, those requirements can be easily fulfilled on a regular grid. However, the situation is quite different for the homogeneous Neumann boundary condition \eqref{eq-2-1} on the cell-centered grid. For clear illustration, we first express the difference operators in the form of matrix-vector multiplication such that
\begin{equation}\label{eq-2-2}
\begin{aligned}
&\delta_x^+u_{j,k}=(D_x^+u)_{j,k},\quad \delta_x^-u_{j,k}=(D_x^-u)_{j,k},\quad  \delta_y^+u_{j,k}=(u(D_y^+)^T)_{j,k},\quad  \delta_y^-u_{j,k}=(u(D_y^-)^T)_{j,k},
\end{aligned}
\end{equation}
where $D_x^+$, $D_x^-$, $D_y^+$, $D_y^-$ are related first-order differential matrices, respectively.  The major difference between the mostly used periodic boundary condition and the homogeneous Neumann boundary condition is now reflected in the corresponding differential matrix. Utilizing the boundary condition \eqref{eq-2-1}, the entries in $D_\alpha^+$ and $D_\alpha^-$, $\alpha=x$ or $y$ are given by
\begin{equation}\label{eq-2-3}
D_\alpha^+=\frac{1}{h_\alpha}\left(\begin{array}{ccccc}
-1 & 1 &  &  &  \\ 
& -1 & 1 &  &  \\ 
&  & \ddots & \ddots &  \\ 
&  &  & -1 & 1 \\ 
&  &  &  & 0
\end{array}\right),\quad D_\alpha^-=\frac{1}{h_\alpha}\left(\begin{array}{ccccc}
0 &  &  &  &  \\ 
-1 & 1 &  &  &  \\ 
& \ddots & \ddots &  &  \\ 
&  & -1 & 1 &  \\ 
&  &  & -1 & 1
\end{array}\right).
\end{equation}
 It is clearly found that the last row of $D_\alpha^+$ and the first row of $D_\alpha^-$ are all zeros, while for the periodic case these are all circulant matrices. Similar observation can be noticed on the associated matrix of the discrete Laplace
\begin{equation}\label{eq-2-4}
\Delta_hu_{j,k}=(D_x^2u)_{j,k}+(uD_y^2)_{j,k},
\end{equation}
where
\begin{equation}\label{eq-2-5}
D_\alpha^2=\frac{1}{h_\alpha^2}\left(\begin{array}{ccccc}
-1 & 1 &  &  &  \\ 
1 & -2 & 1 &  &  \\ 
& \ddots & \ddots &  &  \\ 
&  & 1 & -2 & 1 \\ 
&  &  & 1 & -1
\end{array}\right).
\end{equation}
Nevertheless, $D_\alpha^2$ is still symmetric. Based on the matrix forms of above difference operators, we have the following lemmas.

\begin{lem}
	
With the homogeneous Neumann boundary, the first-order difference operators $\delta_\alpha^+$ and $\delta_\alpha^-$, $\alpha=x$ or $y$ cannot be commutative for all $1\leq j\leq N_x,~ 1\leq k\leq N_y$, 
which in the matrix form is
\begin{equation}\label{eq-2-6}
 D_\alpha^+D_\alpha^-\neq D_\alpha^-D_\alpha^+. 
\end{equation}
Also we have $D_\alpha^+\neq -(D_\alpha^-)^T$. Furthermore, the discrete Laplace operator do not satisfies
\begin{equation}\label{eq-2-7}
\Delta_hu_{j,k}=\delta_x^+\delta_x^-u_{j,k}+\delta_y^+\delta_y^-u_{j,k},
\end{equation}
for all $1\leq j\leq N_x,~ 1\leq k\leq N_y$, which in the matrix form means
\begin{equation}\label{eq-2-77}
D_\alpha^2\neq D_\alpha^+D_\alpha^-.
\end{equation}
Instead, we have
\begin{equation}\label{eq-2-8}
D_\alpha^2=-(D_\alpha^+)^TD_\alpha^+=-(D_\alpha^-)^TD_\alpha^-.
\end{equation}

\end{lem} 

\begin{prf}
	These properties can be verified directly from the matrix forms of the corresponding difference quotients \eqref{eq-2-2}-\eqref{eq-2-5}.
\end{prf}

\begin{rmk}
	Notice that with the periodic boundary condition or zeros boundary condition, these invalid properties \eqref{eq-2-6} and \eqref{eq-2-77} are satisfied on a regular grid and play important roles in the construction of conservative schemes.
\end{rmk}

\begin{lem}
	For any grid function $u\in V_h$, we have
	\begin{equation}\label{eq-2-9}
	(\Delta_hu,u)_h=-(\nabla_h^+u,\nabla_h^+u)_h=-(\nabla_h^-u,\nabla_h^-u)_h.
	\end{equation}
\end{lem}

\begin{prf}
	Utilizing the results of \eqref{eq-2-4}, \eqref{eq-2-5} and \eqref{eq-2-8}, we obtain
	\[
	\begin{aligned}
	(\Delta_hu,u)_h&=(D_x^2u+uD_y^2,u)_h\\
	&=(D_x^2u,u)_h+(uD_y^2,u)_h\\
	&=-(D_x^+u,D_x^+u)_h-(u(D_y^+)^T,u(D_y^+)^T)_h\\
	&=-(\nabla_h^+u,\nabla_h^+u)_h.
	\end{aligned}	
	\]
	Similar proof can be done for $(\Delta_hu,u)_h=-(\nabla_h^-u,\nabla_h^-u)_h.$
\end{prf}

We shall present the semi-discretization for the two dimensional sine-Gordon equation \eqref{eq-1-1} with homogeneous Neumann boundary condition \eqref{eq-1-3} as follows:
\begin{equation}\label{eq-2-10}
\big(U_{j,k}\big)_{tt}=\Delta_h U_{j,k}-\Phi_{j,k}\sin U_{j,k},\quad \forall~1\leq j\leq N_x,~ 1\leq k\leq N_y,
\end{equation}
which can be also arranged in the matrix form as
\begin{equation}\label{eq-2-11}
U_{tt}=\Delta_hU-\Phi\cdot\sin U,
\end{equation}
where $U=(U_{j,k})$, and ``$\cdot$" represents the element multiplication of matrices.

\begin{thm}
	The scheme \eqref{eq-2-11} admits a semi-discrete energy conservation law
	\begin{equation}\label{eq-2-12}
	\frac{d}{dt}H(t)=0,\quad \mbox{with}\quad H(t)=\frac{1}{2}\Big(\|U_t\|_h^2+\|\nabla_h^+U\|_h^2+2(\Phi\cdot(1-\cos U),1)_h\Big).
	\end{equation}
\end{thm}

\begin{prf}
	Taking the discrete inner product of \eqref{eq-2-11} with $U_t$ yields
	\[
	(U_{tt},U_t)_h=(\Delta_hU,U_t)_h-(\Phi\cdot\sin U,U_t)_h.
	\]
	Using the identity \eqref{eq-2-9}, we have
		\[
	(U_{tt},U_t)_h=-(\nabla_h^+U,\nabla_h^+U_t)_h-(\Phi\cdot\sin U,U_t)_h, 
	\]
	which is equivalent to
	\[
	\frac{1}{2}\frac{d}{dt}\Big((U_t,U_t)_h+(\nabla_h^+U,\nabla_h^+U)_h+2(\Phi\cdot(1-\cos U),1)_h\Big)=0.
	\]
	This completes the proof.
\end{prf}

\begin{thm}
	The scheme \eqref{eq-2-11} admits a semi-discrete canonical Hamiltonian structure.
\end{thm}

\begin{prf}
	Rearranging the matrix form of $U$ into a vector form $\bm u$ such that
	\[
	\bm u=\big(U_{1,1}, U_{2,1},\cdots,U_{N_x,1}, U_{1,2},\cdots, U_{N_x,N_y}\big)^T.
	\] 
	Further expanding the discrete Laplace operator, we can transform \eqref{eq-2-11} as
	\begin{equation}\label{eq-2-13}
	\bm u_{tt}=D\bm u-\bm \phi\cdot\sin\bm u,
	\end{equation}
	where $D=I_y\otimes D_x^2+D_y^2\otimes I_x$ is symmetric, and $\otimes$ is the Kronecker product, $I_x, I_y$ are the identical matrices of dimensional $N_x\times N_x$ and $N_y\times N_y$, respectively. Consequently, after introducing an auxiliary variable $\bm v=\bm u_t$, the system \eqref{eq-2-13} can be rewritten in the form of Hamiltonian structure
	\begin{equation}
	\left(\begin{array}{c}
	\bm u\\ 
	\bm v
	\end{array} \right)_t=\left(\begin{array}{cc}
	0 & I \\ 
	-I & 0 
	\end{array} \right)\left(\begin{array}{c}
	\nabla_u H\\ 
	\nabla_v H
	\end{array} \right),
	\end{equation} 
	where $H=\frac{1}{2}\big(\bm v^T\bm v-\bm u^TD\bm u+2\bm \phi^T(1-\cos\bm u)\big)$ and $I$ is the identical matrix of dimension $(N_x\times N_y)^2$. 
\end{prf}

\begin{rmk}
	As demonstrated above, the semi-discretization based on a cell-centered grid and central difference quotients is extremely simple and efficient for dealing with the homogeneous Neumann boundary condition. However, it only has second-order accuracy in space and cannot straightly generalize to higher-order schemes. As a consequence, we introduce an alternative discrete strategy as follows.
\end{rmk}

\section{Conservative semi-discretization with SBP operators}

The second strategy to preserve the semi-discrete properties is inspired by the SBP operators on the regular grid for approximating the second derivative \cite{mn04}. In this section, we first present a simple SBP operator of second-order accuracy to illustrate the conservative property of the corresponding semi-discretization for the sine-Gordon equation, and leave the detailed definition of SBP operators for general Neumann boundary conditions as well as a fourth-order approximation in Appendix A.

 Distinguished from the cell-centered grid, variables are normally stored at each cell boundary as follows:
\[
\Omega_h^\prime=\Big\{(x_j,y_k)\big|x_j=a+jh_x, y_k=c+kh_y, 0\leq j\leq N_x, 0\leq k\leq N_y\Big\},
\]
and the space of grid function becomes $V_h^\prime=\{u=(u_{j,k})|(x_j,y_k)\in\Omega_h^\prime\}$.

Applying the central difference quotients to approximate the homogeneous Neumann boundary condition yields
\begin{equation}\label{eq-3-1}
u_{-1,k}=u_{1,k}, \quad u_{j,-1}=u_{j,1}, \quad u_{N_x+1,k}=u_{N_x-1,k}, \quad u_{j,N_y+1}=u_{j,N_y-1},
\end{equation}
where the ghost points are denoted by the subscripts $-1, N_x+1, N_y+1$. Subsequently, the matrix forms of the first-order finite difference operators $\delta_\alpha^+$, $\delta_\alpha^-$ as well as the Laplace operator $\bar{\Delta}_h$ now become
\begin{align}\label{eq-3-2}
&\bar{D}_\alpha^+=\frac{1}{h_\alpha}\left(\begin{array}{ccccc}
-1 & 1 &  &  &  \\ 
& -1 & 1 &  &  \\ 
&  & \ddots & \ddots &  \\ 
&  &  & -1 & 1 \\ 
&  &  & 1 & -1
\end{array}\right),\quad \bar{D}_\alpha^-=\frac{1}{h_\alpha}\left(\begin{array}{ccccc}
1 & -1 &  &  &  \\ 
-1 & 1 &  &  &  \\ 
& \ddots & \ddots &  &  \\ 
&  & -1 & 1 &  \\ 
&  &  & -1 & 1
\end{array}\right),
\end{align}
and
\begin{equation}\label{eq-3-3}
\bar{\Delta}_hu_{j,k}=(\bar{D}_x^2u)_{j,k}+(u\bar{D}_y^2)_{j,k},
\end{equation}
where
\begin{align}\label{eq-3-33}
&\bar{D}_\alpha^2=\frac{1}{h_\alpha^2}\left(\begin{array}{ccccc}
-2 & 2 &  &  &  \\ 
1 & -2 & 1 &  &  \\ 
& \ddots & \ddots &  &  \\ 
&  & 1 & -2 & 1 \\ 
&  &  & 2 & -2
\end{array}\right).
\end{align}
It can be easily verified that $\bar{D}_\alpha^+\bar{D}_\alpha^-\neq \bar{D}_\alpha^-\bar{D}_\alpha^+$ and $\bar{D}_\alpha^+\neq -(\bar{D}_\alpha^-)^T$, as the results obtained in the case of above cell-center grid. Besides, the symmetry of the second-order difference operator $\bar{D}_\alpha^2$ is also violated, and thereby no matrix $\bar{D}_\alpha$ exists such that $\bar{D}_\alpha^2=-(\bar{D}_\alpha)^T\bar{D}_\alpha$, that is the reason why the direct extension of the central difference scheme to the homogeneous Neumann boundary condition cannot lead to a conservative semi-discretization. However, we can split $\bar{D}_\alpha^2$ as $\bar{D}_\alpha^2=\Lambda_\alpha^{-1} D_\alpha^2$, where $D_\alpha^2$ is defined the same as \eqref{eq-2-5} and
\begin{equation}\label{lam}
\Lambda_\alpha=\left(\begin{array}{ccccc}
	\frac{1}{2} &  &  &  &  \\ 
	& 1 &  &  &  \\ 
	& & \ddots &  &  \\ 
	&  & & 1 &  \\ 
	&  &  & & \frac{1}{2}
\end{array}\right),
\end{equation}
is called the weighed matrix. Such splitting is just inspired by the SBP operators and can be generalized to derive high-order schemes (see Appendix A). As a consequence, we define another weighted inner product for grid functions $u,v\in V_h^\prime$ as follows:
\begin{equation}\label{eq-3-0}
\begin{aligned}
&(u,v)_{\Lambda}=(\Lambda_xu\Lambda_y,v)_h,\quad \|u\|_\Lambda^2=(u,u)_{\Lambda},\\
&(\nabla_h^+u,\nabla_h^+v)_{\Lambda}=(\Lambda_x^{-1/2}\delta_x^+u,\Lambda_x^{-1/2}\delta_x^+v)_{\Lambda}+(\delta_y^+u\Lambda_y^{-1/2},\delta_y^+v\Lambda_y^{-1/2})_{\Lambda},
\end{aligned}
\end{equation}
where $\Lambda_\alpha=\Lambda_\alpha^T>0$, $\alpha=x$ or $y$.

\begin{rmk}
	Notice that in one dimensional case, the above two kinds of definitions \eqref{eq-2-0} and \eqref{eq-3-0} of discrete inner products correspond to different numerical quadrature rules. For the definition \eqref{eq-2-0}, it belongs to the rectangle formula with the cell-centered nodes being its quadrature points. While the second one \eqref{eq-3-0} is actually the composite trapezoidal rule, and the diagonal elements in $\Lambda$ are the relevant quadrature weights.
\end{rmk}

\begin{lem}
	For any grid function $u\in V_h^\prime$, we have
	\begin{equation}\label{eq-3-4}
	(\bar{\Delta}_h u,u)_\Lambda=-(\nabla_h^+u,\nabla_h^+u)_\Lambda=-(\nabla_h^-u,\nabla_h^-u)_\Lambda.
	\end{equation}
\end{lem}

\begin{prf}
	Using the definition of the weighted inner product \eqref{eq-3-0}, the relations between $\bar{D}_\alpha^2$ and $D_\alpha^2$ and the result of \eqref{eq-2-8}, we can derive
	\[
	\begin{aligned}
	(\bar{\Delta}_h u,u)_\Lambda&=(\Lambda_x^{-1} D_x^2u+uD_y^2\Lambda_y^{-1} ,u)_\Lambda\\
	&=(D_x^2u\Lambda_y+\Lambda_xuD_y^2 ,u)_h\\
	&=-(D_x^+u\Lambda_y,D_x^+u)_h-(\Lambda_xu(D_y^+)^T ,u(D_y^+)^T)_h\\
	&=-(\Lambda_x^{-1}D_x^+u,D_x^+u)_\Lambda-(uD_y^+\Lambda_y^{-1} ,uD_y^+)_\Lambda\\
	&=-(\Lambda_x^{-1/2}D_x^+u,\Lambda_x^{-1/2}D_x^+u)_\Lambda-(uD_y^+\Lambda_y^{-1/2} ,uD_y^+\Lambda_y^{-1/2})_\Lambda\\
	&=-(\nabla_h^+u,\nabla_h^+u)_\Lambda
	\end{aligned}
	\]
Similar proof can be obtained for the second equality in \eqref{eq-3-4}.
\end{prf}

The semi-discrete scheme based on the regular grid becomes
\begin{equation}\label{eq-3-5}
\big(U_{j,k}\big)_{tt}=\bar{\Delta}_h U_{j,k}-\Phi_{j,k}\sin U_{j,k},\quad \forall~0\leq j\leq N_x,~ 0\leq k\leq N_y,
\end{equation}
which also admits a discrete energy conservation law and a non-canonical Hamiltonian structure as follows.

\begin{thm}
	The scheme \eqref{eq-3-5} preserves a semi-discrete energy in the sense of a weighted norm, that is
	\begin{equation}\label{eq-3-6}
	\frac{d}{dt}H(t)=0,\quad\mbox{with}\quad H(t)=\frac{1}{2}\Big(\|U_t\|_\Lambda^2+\|\nabla_h^+U\|_\Lambda^2+2(\Phi\cdot(1-\cos U),1)_\Lambda\Big)
	\end{equation}
\end{thm}

\begin{prf}
	Taking the weighted inner product of \eqref{eq-3-5} with $U_t$, we have
	\[
	(U_{tt},U_t)_\Lambda=(\bar{\Delta}_hU,U_t)_\Lambda-(\Phi\cdot\sin U,U_t)_\Lambda.
	\]
	Using the identity \eqref{eq-3-4}, we have
	\[
	(U_{tt},U_t)_\Lambda=-(\nabla_h^+U,\nabla_h^+U_t)_\Lambda-(\Phi\cdot\sin U, U_t)_\Lambda, 
	\]
	which is equivalent to
	\[
	\frac{1}{2}\frac{d}{dt}\Big((U_t,U_t)_\Lambda+(\nabla_h^+U,\nabla_h^+U)_\Lambda+2(\Phi\cdot(1-\cos U),1)_\Lambda\Big)=0.
	\]
	This completes the proof.
\end{prf}

\begin{thm}
	The scheme \eqref{eq-3-5} has a non-canonical Hamiltonian structure as follows:
	\begin{equation}\label{eq-3-7}
	\left(\begin{array}{c}
	\bm u\\ 
	\bm v
	\end{array} \right)_t=\left(\begin{array}{cc}
	0 & \bar{\Lambda}^{-1} \\ 
	-\bar{\Lambda}^{-1} & 0 
	\end{array} \right)\left(\begin{array}{c}
	\nabla_u H\\ 
	\nabla_v H
	\end{array} \right),
	\end{equation}
	where 
	\[
	H=\frac{1}{2}\big(\bm v^T\bar{\Lambda}\bm v-\bm u^T\bar{D}\bm u\big)+\mathbbm{1}^T\bar{\Lambda}\big(\bm \phi\cdot(1-\cos\bm u)\big),
	\]
	 is equivalent to that in \eqref{eq-3-6} with an alternative vector representation, and $\bar{\Lambda}=\Lambda_y\otimes\Lambda_x$, $\bar{D}=\Lambda_y\otimes D_x^2+D_y^2\otimes\Lambda_x$, $\mathbbm{1}=(1,1,\cdots,1)^T$ is an identical vector of dimension $N_x\times N_y$.
\end{thm}

\begin{prf}
	Multiplying $\Lambda_x$, $\Lambda_y$ on both sides of \eqref{eq-3-5}, we have
	\[
	\Lambda_xU_{tt}\Lambda_y=D_x^2U\Lambda_y+\Lambda_xUD_y^2-\Lambda_x(\Phi\cdot\sin U)\Lambda_y.
	\]
	Using the identity of Kronecker product for arbitrary matrices $A, B, X$
	\[
	\mbox{vec}(AXB)=(B^T\otimes A)\mbox{vec}(X),
	\]
	where \mbox{vec}($X$) denotes the vectorization of the matrix $X$ formed by stacking the columns of $X$ into a single column vector, we can derive
	\[
	(\Lambda_y\otimes\Lambda_x)\bm u_{tt}=(\Lambda_y\otimes D_x^2+D_y^2\otimes\Lambda_x)\bm u-(\Lambda_y\otimes\Lambda_x)(\bm\phi\cdot\sin\bm u).
	\]
	Let $\bm v=\bm u_t$, we can obtain the non-canonical Hamiltonian structure \eqref{eq-3-7}. 
\end{prf}

\section{Symplectic integrators and energy-preserving schemes}

It is natural to generalize the semi-discrete properties to the fully-discrete analogs and construct symplectic integrators and energy-preserving schemes, respectively. Since the semi-discrete Hamiltonian structures have been proved for both methods based on the cell-center grid and the second-order SBP operator, it is straightforward to derive symplectic integrators by applying any symplectic Runge-Kutta method in time integration. Without loss of generality, we simply adopt the symplectic implicit midpoint method and result two symplectic integrators named {\bf SIM/CC} and {\bf SIM/SBP2}, respectively. For the construction of energy-preserving schemes, we utilize the SAV approach and derive linearly implicit schemes with constant coefficient matrix, which is given in details as follows.

\subsection{Linearly implicit energy-preserving scheme based on the SAV approach}

Let $\mathcal{H}_1=\iint_\Omega\phi(1-\cos(u))dxdy$ where $\phi(x,y)$ is a function independent of time. We assume $\phi(x,y)>0$ and introduce a scalar auxiliary variable $r=\sqrt{\mathcal{H}_1}$. Then the system \eqref{eq-1-5} can be rewritten as
\begin{equation}\label{eq-5-1}
\begin{aligned}
& u_t=v,\\
& v_t=\Delta u-b(u)r,\\
& r_t=\frac{1}{2}\big(b(u),u_t\big),
\end{aligned}
\end{equation}
where $b(u)=\frac{\phi\sin u}{\sqrt{\mathcal{H}_1(u)}}$ and $(\cdot,\cdot)$ represents the continuous $L^2$-inner product. Taking the inner product of the above system with $v_t, u_t$ and $2r$, respectively, and summing them together, we obtain the modified energy conservation law
\begin{equation}\label{eq-5-2}
\frac{d}{dt}\tilde{\mathcal{H}}=0,\quad \mbox{with} \quad\tilde{\mathcal{H}}=\frac{1}{2}\iint_\Omega \big(v^2+|\nabla u|^2\big)dxdy+r^2,
\end{equation}
which is equivalent to the original one in the continuous circumstance. 

Next, utilizing the two kinds of strategies to semi-discretize the system \eqref{eq-5-1}, we have
\begin{equation}\label{eq-5-1-1}
\begin{aligned}
& U_t=V,\\
& V_t=\Delta_h U-b(U)R,\\
& R_t=\frac{1}{2}\big(b(U),U_t\big)_h,
\end{aligned}
\end{equation}
for the cell-centered case and
\begin{equation}\label{eq-5-1-2}
\begin{aligned}
& U_t=V,\\
& V_t=\bar{\Delta}_h U-b(U)R,\\
& R_t=\frac{1}{2}\big(b(U),U_t\big)_\Lambda,
\end{aligned}
\end{equation}
for the regular grid with SBP operators, respectively. In analogous to the continuous case, taking the corresponding discrete inner products of \eqref{eq-5-1-1} and \eqref{eq-5-1-2} with $V_t, U_t$ and $2R$, and summing them together, we can obtain two modified semi-discrete energy conservation law
\begin{equation}
\frac{d}{dt}\tilde{H}=0,\quad\mbox{with}\quad \tilde{H}=\frac{1}{2}\big(\|V\|_h^2+\|\nabla_h^+U\|_h^2\big)+R^2,
\end{equation} 
and 
\begin{equation}
\frac{d}{dt}\tilde{H}=0,\quad\mbox{with}\quad \tilde{H}=\frac{1}{2}\big(\|V\|_\Lambda^2+\|\nabla_h^+U\|_\Lambda^2\big)+R^2.
\end{equation} 

\begin{rmk}
	Notice that the only differences between \eqref{eq-5-1-1} and \eqref{eq-5-1-2} are the definitions of discrete Laplace and the inner products. For the conventional application of the SAV approach, one only need the standard discrete $L^2$-inner product to guarantee the energy conservation or dissipation. However, for the second discrete strategy in Section 3, one must take the weighted inner product to match with the SBP operator. Consequently, the definitions of discrete energy are distinguished as well.
\end{rmk}

Based on the semi-discretizations \eqref{eq-5-1-1} and \eqref{eq-5-1-2}, we can construct linearly implicit energy-preserving schemes. For simplicity, we only take \eqref{eq-5-1-1} for illustration and the construction for \eqref{eq-5-1-2} is very similar.

A second-order linearly implicit scheme is obtained by the application of the Crank-Nicolson method and an extrapolation technique on \eqref{eq-5-1-1}, which yields
\begin{equation}\label{eq-5-3}
\begin{aligned}
& \frac{U^{n+1}-U^n}{\tau}=V^{n+1/2},\\
& \frac{V^{n+1}-V^n}{\tau}=\Delta_h U^{n+1/2}-\bar{B}^nR^{n+1/2},\\
& R^{n+1}-R^n=\frac{1}{2}\big(\bar{B}^n,  U^{n+1}-U^n \big)_h,
\end{aligned}
\end{equation}
where $\bar{B}^n=b(\bar{U}^{n+1/2})$, $\bar{U}^{n+1/2}=\frac{1}{2}(3U^n-U^{n-1})$ is an explicit extrapolation with a second-order approximation of $U(t^{n+1/2})$. Taking the inner product of the first two equations in \eqref{eq-5-3} with $(V^{n+1}-V^n)/\tau$, $(U^{n+1}-U^n)/\tau$, multiplying the last equation with $(R^{n+1}+R^n)/\tau$, and adding them together, we obtain the energy conservation law:
\begin{equation}\label{eq-5-4-1}
\tilde{H}^{n+1}=\tilde{H}^{n},\quad \mbox{with} \quad\tilde{H}^n=\frac{1}{2}\big(\|V^n\|_h^2+\|\nabla_h^+U^n\|_h^2\big)+(R^n)^2.
\end{equation}
Similarly, we can obtain the fully discrete scheme based on \eqref{eq-5-1-2}. We denote the two resulting schemes by {\bf SAV/CC} and {\bf SAV/SBP2}, respectively, which are both energy-preserving with second-order accuracy.

 To demonstrate the linearly implicit advantage with constant coefficient matrix, we first eliminate $V^{n+1}$, $R^{n+1}$ from \eqref{eq-5-3} and obtain 
\begin{equation}\label{eq-5-4}
\frac{U^{n+1}-U^n}{\tau}=V^n+\frac{\tau}{2}\Big(\Delta_h U^{n+1/2}-\bar{B}^nR^{n}- \frac{\bar{B}^n}{4}(\bar{B}^n,U^{n+1}-U^n)_h\Big),
\end{equation}
which can be further arranged as
\begin{equation}\label{eq-5-5}
(I-\frac{\tau^2}{4}\Delta_h)U^{n+1}+\frac{\tau^2}{8}\bar{B}^n(\bar{B}^n,U^{n+1})_h=(I+\frac{\tau^2}{4}\Delta_h)U^n+\frac{\tau^2}{8}\bar{B}^n (\bar{B}^n,U^n)_h-\frac{\tau^2}{2}\bar{B}^nR^{n}+\tau V^n.
\end{equation}
Denote the righthand side of \eqref{eq-5-5} by $C^n$. Multiplying \eqref{eq-5-5} with $(I-\frac{\tau^2}{4}\Delta_h)^{-1}$, then taking the discrete inner product with $\bar{B}^n$, we can obtain
\begin{equation}\label{eq-5-6}
(\bar{B}^n,U^{n+1})_h+\frac{\tau^2}{8}\gamma^n(\bar{B}^n,U^{n+1})_h=\big(\bar{B}^n,(I-\frac{\tau^2}{4}\Delta)^{-1}C^n\big)_h
\end{equation}
where $\gamma^n=(\bar{B}^n,(I-\frac{\tau^2}{4}\Delta_h)^{-1}\bar{B}^n)$. Hence
\begin{equation}\label{eq-5-7}
(\bar{B}^n,U^{n+1})_h=\dfrac{\big(\bar{B}^n,(I-\frac{\tau^2}{4}\Delta_h)^{-1}C^n\big)_h}{1+\frac{\tau^2}{8}\gamma^n}.
\end{equation}
Substituting the result into \eqref{eq-5-5} will give the solution of $U^{n+1}$, and subsequently $V^{n+1}$, $R^{n+1}$ from \eqref{eq-5-3}. In summary, we only have to solve a linear system with constant coefficient matrix $(I-\frac{\tau^2}{4}\Delta_h)$ twice, for which fast solvers are available.

\section{Numerical experiments}

In this section, we present extensive numerical examples to demonstrate the effectiveness of the proposed schemes in the conservation of discrete energy for the two-dimensional sine-Gordon equation \eqref{eq-1-1}-\eqref{eq-1-3}. All the examples are taken with the homogeneous Neumann boundary condition \eqref{eq-1-3} except the first example regarding the convergence test. Spatial steps are uniformly choose as $h_x=h_y=h$ for simplicity. 

\subsection{Accuracy tests}

In this experiment, we consider the sine-Gordon equation \eqref{eq-1-1} with $\phi(x,y)=1$ and initial conditions
\[
\begin{aligned}
&f_1(x,y)=4\tan^{-1}\exp(x+y),\quad -7\leq x, y\leq 7,\\
&f_2(x,y)=-\frac{4\exp(x+y)}{1+\exp(2x+2y)},\quad -7\leq x, y\leq 7,
\end{aligned}
\] 
and non-homogeneous Neumann boundary conditions
\[
\begin{aligned}
&\frac{\partial u}{\partial x}=g_1(x,y,t)=\frac{4\exp(x+y+t)}{\exp(2t)+\exp(2x+2y)},~~\mbox{for}~x=-7 ~\mbox{and}~ x=7, ~-7\leq y\leq 7,~ t>0,\\
&\frac{\partial u}{\partial y}=g_2(x,y,t)=\frac{4\exp(x+y+t)}{\exp(2t)+\exp(2x+2y)},~~\mbox{for}~y=-7 ~\mbox{and}~ y=7, ~-7\leq x\leq 7,~ t>0.\\
\end{aligned}
\]
The analytical solution of this problem is given by
\[
u(x,y,t)=4\tan^{-1}\exp(x+y-t).
\]
The implementation of above non-homogeneous boundary conditions for the SBP operator based schemes {\bf SIM/SBP2} and {\bf SAV/SBP2} is presented in Appendix A, where the derivation of energy conservation law for general SBP operators are given in details. Furthermore, we show an explicit form of a fourth-order accuracy SBP operator and the corresponding symplectic integrator is denoted by {\bf SIM/SBP4} while the energy-preserving scheme is denoted by {\bf SAV/SBP4}. 

We run the simulations till $t=1$ by all the presented schemes with different mesh sizes and compute the $L^\infty$, $L^2$-errors as well as the related convergence orders, which are listed in Table~\ref{tab-5-1}-\ref{tab-5-3}, respectively. It can be seen from Table~\ref{tab-5-1} and Table~\ref{tab-5-2} that the energy-preserving schemes {\bf SAV/CC}, {\bf SAV/SBP2} and the symplectic integrators {\bf SIM/CC}, {\bf SIM/SBP2} are of second-order accuracy both in spatial and temporal directions, because we have simultaneously decrease the mesh sizes $h$ and $\tau$ by half. While in the test of Table~\ref{tab-5-3}, the temporal step is refined gradually by a quarter, such that the results of convergence orders imply a fourth-order accuracy in space and a second-order accuracy in time for both schemes {\bf SAV/SBP4} and {\bf SIM/SBP4}.
For simplicity,  we only present the corresponding results of the four second-order schemes in the following experiments,  which however are still valid for schemes {\bf SAV/SBP4} and {\bf SIM/SBP4}.

\begin{table}[H]
	\centering
	\caption{Spatial and temporal accuracy tests of schemes {\bf SAV/CC} and {\bf SIM/CC}.}\label{tab-5-1}
	\begin{tabular*}{\textwidth}[h]{@{\extracolsep{\fill}}c l l l l l} \hline
		& $(h,\tau)$ & $L^\infty$-error & order & $L^2$-error & order  \\ \hline
		\multirow{4}{*}{{\bf SAV/CC}}  & $(\frac{1}{4},\frac{1}{100})$  & 1.0728e-02 & - & 2.2309e-03 & - \\[1ex]
		& $(\frac{1}{8},\frac{1}{200})$ & 2.8264e-03 & 1.92 & 5.5632e-04 & 2.00  \\[1ex]
		& $(\frac{1}{16},\frac{1}{400})$ & 7.2008e-04 & 1.97 & 1.3899e-04 & 2.00  \\[1ex]
		& $(\frac{1}{32},\frac{1}{800})$ & 1.8135e-04 & 1.99 & 3.4742e-05 & 2.00 \\ [1ex] \hline		
		
		\multirow{4}{*}{{\bf SIM/CC}}  & $(\frac{1}{4},\frac{1}{100})$  & 1.0720e-02 & - & 2.2300e-03 & - \\[1ex]
		& $(\frac{1}{8},\frac{1}{200})$ & 2.8243e-03 & 1.92 & 5.5608e-04 & 2.00  \\[1ex]
		& $(\frac{1}{16},\frac{1}{400})$ & 7.1955e-04 & 1.97 & 1.3893e-04 & 2.00  \\[1ex]
		& $(\frac{1}{32},\frac{1}{800})$ & 1.8121e-04 & 1.99 & 3.4727e-05 & 2.00 \\ [1ex] \hline		
	\end{tabular*}
\end{table}

\begin{table}[H]
	\centering
	\caption{Spatial and temporal accuracy tests of schemes {\bf SAV/SBP2} and {\bf SIM/SBP2}.}\label{tab-5-2}
	\begin{tabular*}{\textwidth}[h]{@{\extracolsep{\fill}}c l l l l l} \hline
		& $(h,\tau)$ & $L^\infty$-error & order & $L^2$-error & order  \\ \hline
		\multirow{4}{*}{{\bf SAV/SBP2}}  & $(\frac{1}{4},\frac{1}{100})$  & 1.6972e-02 & - & 2.2915e-03 & - \\[1ex]
		& $(\frac{1}{8},\frac{1}{200})$ & 4.3501e-03 & 1.96 & 5.6636e-04 & 2.02  \\[1ex]
		& $(\frac{1}{16},\frac{1}{400})$ & 1.0952e-03 & 1.99 & 1.4089e-04 & 2.01  \\[1ex]
		& $(\frac{1}{32},\frac{1}{800})$ & 2.7426e-04 & 2.00 & 3.5144e-05 & 2.00 \\ [1ex] \hline		
		
		\multirow{4}{*}{{\bf SIM/SBP2}}  & $(\frac{1}{4},\frac{1}{100})$  & 1.6967e-02 & - & 2.2905e-03 & - \\[1ex]
		& $(\frac{1}{8},\frac{1}{200})$ & 4.3484e-03 & 1.96 & 5.6612e-04 & 2.02  \\[1ex]
		& $(\frac{1}{16},\frac{1}{400})$ & 1.0948e-03 & 1.99 & 1.4083e-04 & 2.01  \\[1ex]
		& $(\frac{1}{32},\frac{1}{800})$ & 2.7416e-04 & 2.00 & 3.5128e-05 & 2.00 \\ [1ex] \hline		
	\end{tabular*}
\end{table}

\begin{table}[H]
	\centering
	\caption{Spatial and temporal accuracy tests of schemes {\bf SAV/SBP4} and {\bf SIM/SBP4}.}\label{tab-5-3}
	\begin{tabular*}{\textwidth}[h]{@{\extracolsep{\fill}}c l l l l l} \hline
		& $(h,\tau)$ & $L^\infty$-error & order & $L^2$-error & order  \\ \hline
		\multirow{4}{*}{\bf SAV/SBP4}  & $(\frac{1}{4},\frac{1}{100})$  & 7.5504e-04 & - & 8.7254e-05 & - \\[1ex]
		& $(\frac{1}{8},\frac{1}{400})$ & 5.7723e-05 & 3.71 & 5.3786e-06 & 4.02  \\[1ex]
		& $(\frac{1}{16},\frac{1}{1600})$ & 3.5637e-06 & 4.02 & 3.1915e-07 & 4.07  \\[1ex]
		& $(\frac{1}{32},\frac{1}{6400})$ & 2.2924e-07 & 3.96 & 1.9562e-08 & 4.03 \\ [1ex] \hline		
		
		\multirow{4}{*}{\bf SIM/SBP4}  & $(\frac{1}{4},\frac{1}{100})$  & 7.4948e-04 & - & 8.5909e-05 & - \\[1ex]
		& $(\frac{1}{8},\frac{1}{400})$ & 5.7329e-05 & 3.71 & 5.2929e-06 & 4.02  \\[1ex]
		& $(\frac{1}{16},\frac{1}{1600})$ & 3.5312e-06 & 4.02 & 3.1357e-07 & 4.08  \\[1ex]
		& $(\frac{1}{32},\frac{1}{6400})$ & 2.2719e-07 & 3.96 & 1.9209e-08 & 4.03 \\ [1ex] \hline		
	\end{tabular*}
\end{table}

\subsection{Superposition of two line solitons}

\begin{figure}[H]
	\centering
	\includegraphics[width=0.24\linewidth]{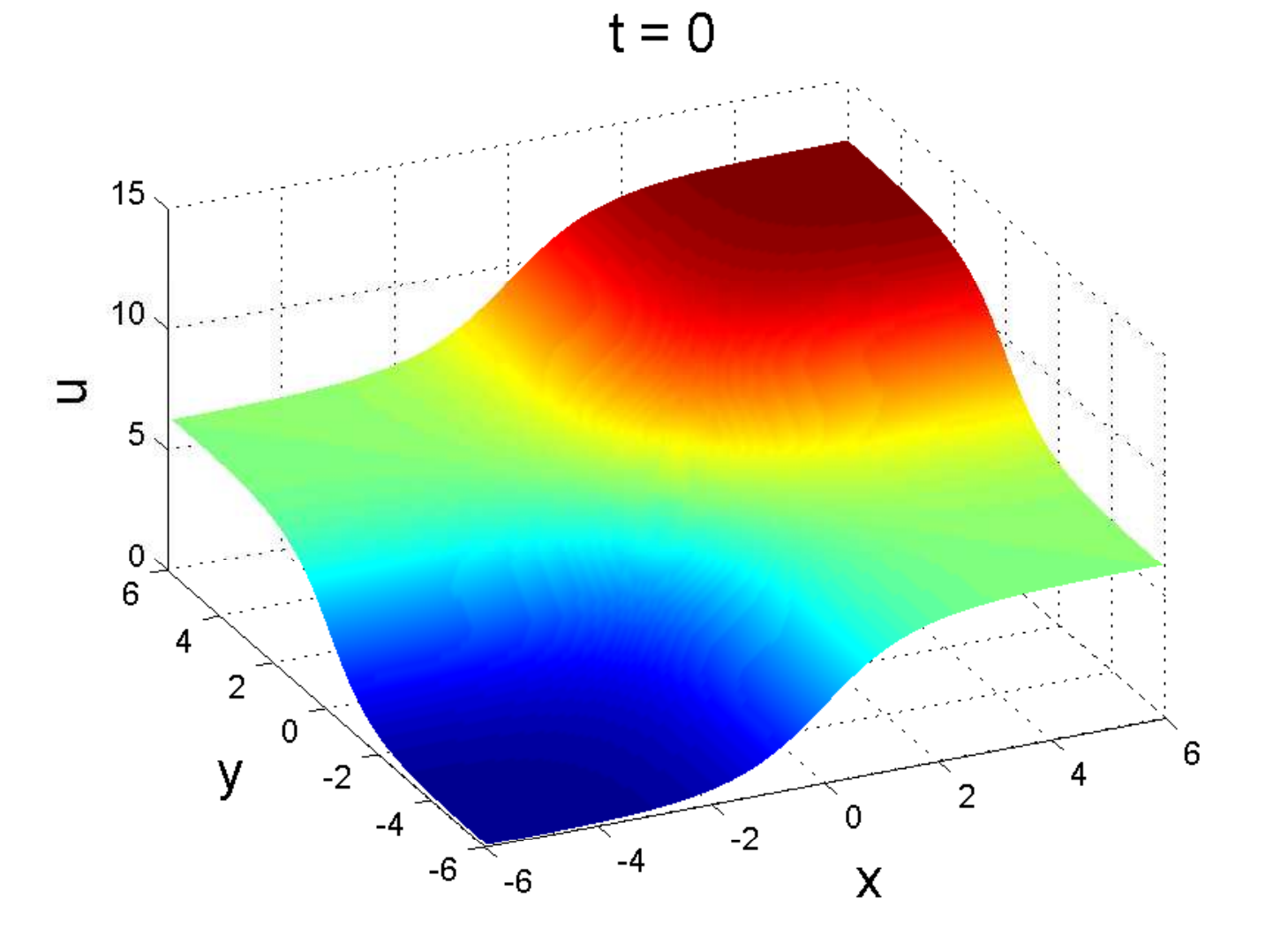}
	\includegraphics[width=0.24\linewidth]{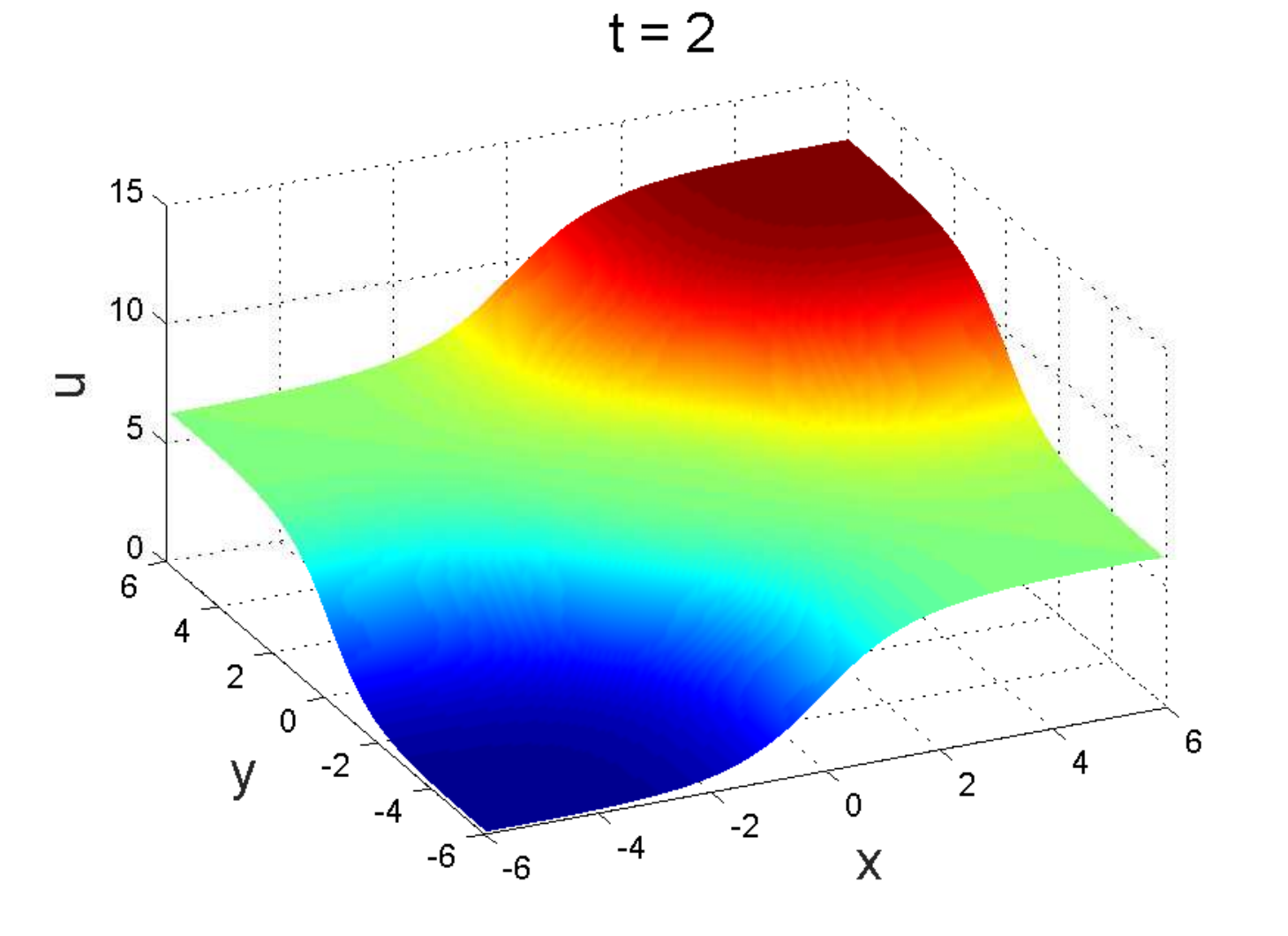}
	\includegraphics[width=0.24\linewidth]{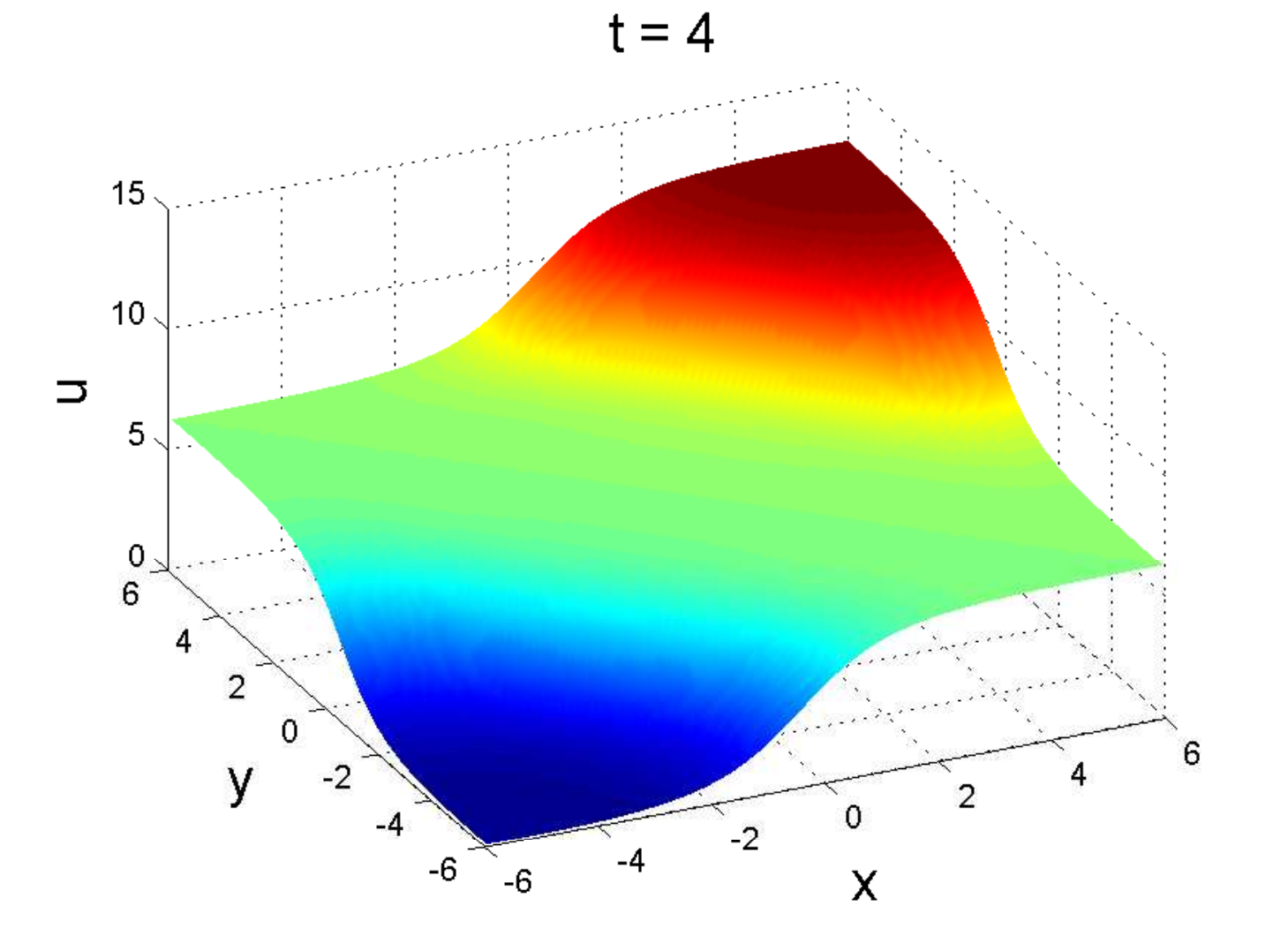}
	\includegraphics[width=0.24\linewidth]{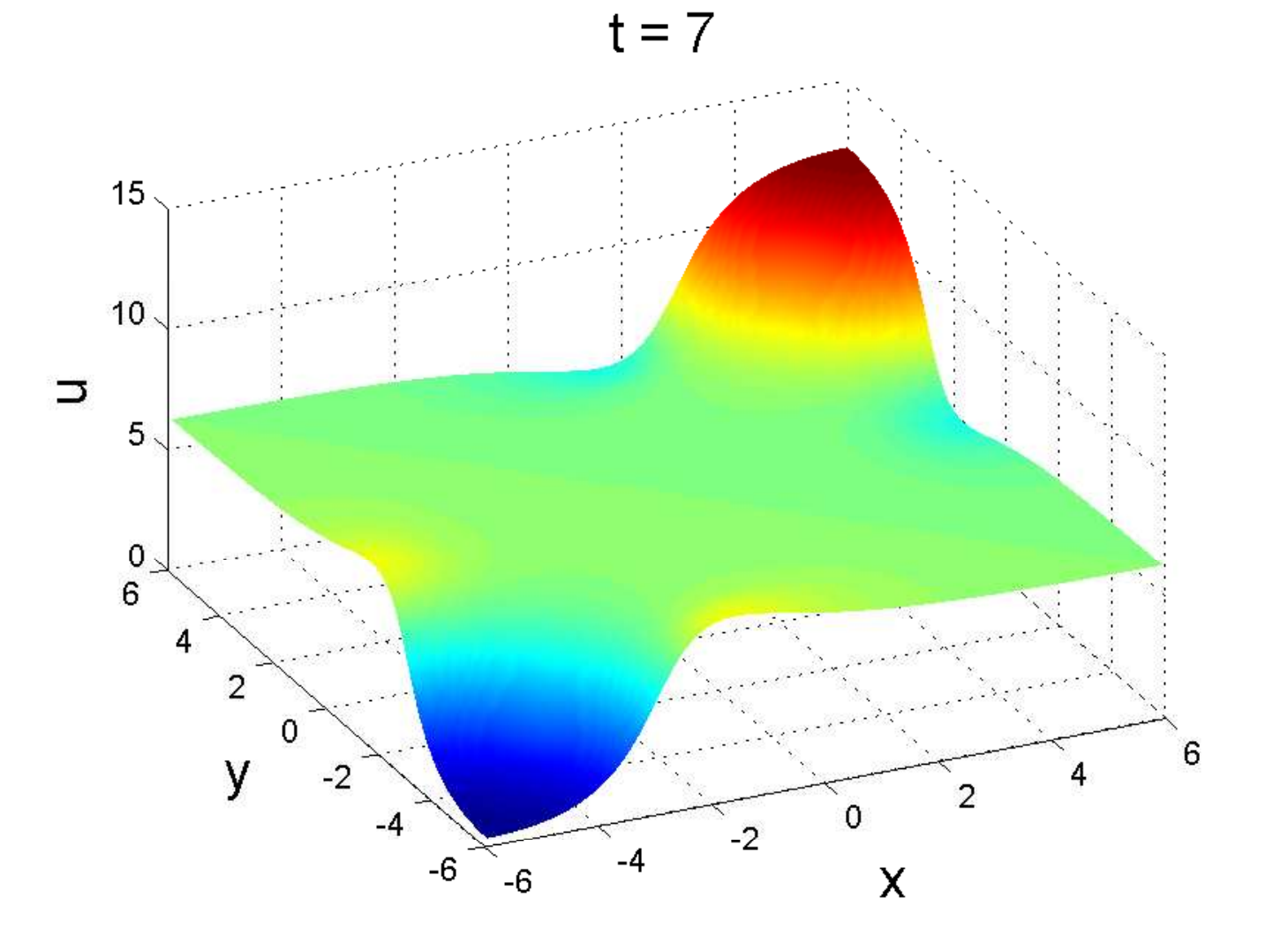}\\[2ex]
	\includegraphics[width=0.24\linewidth]{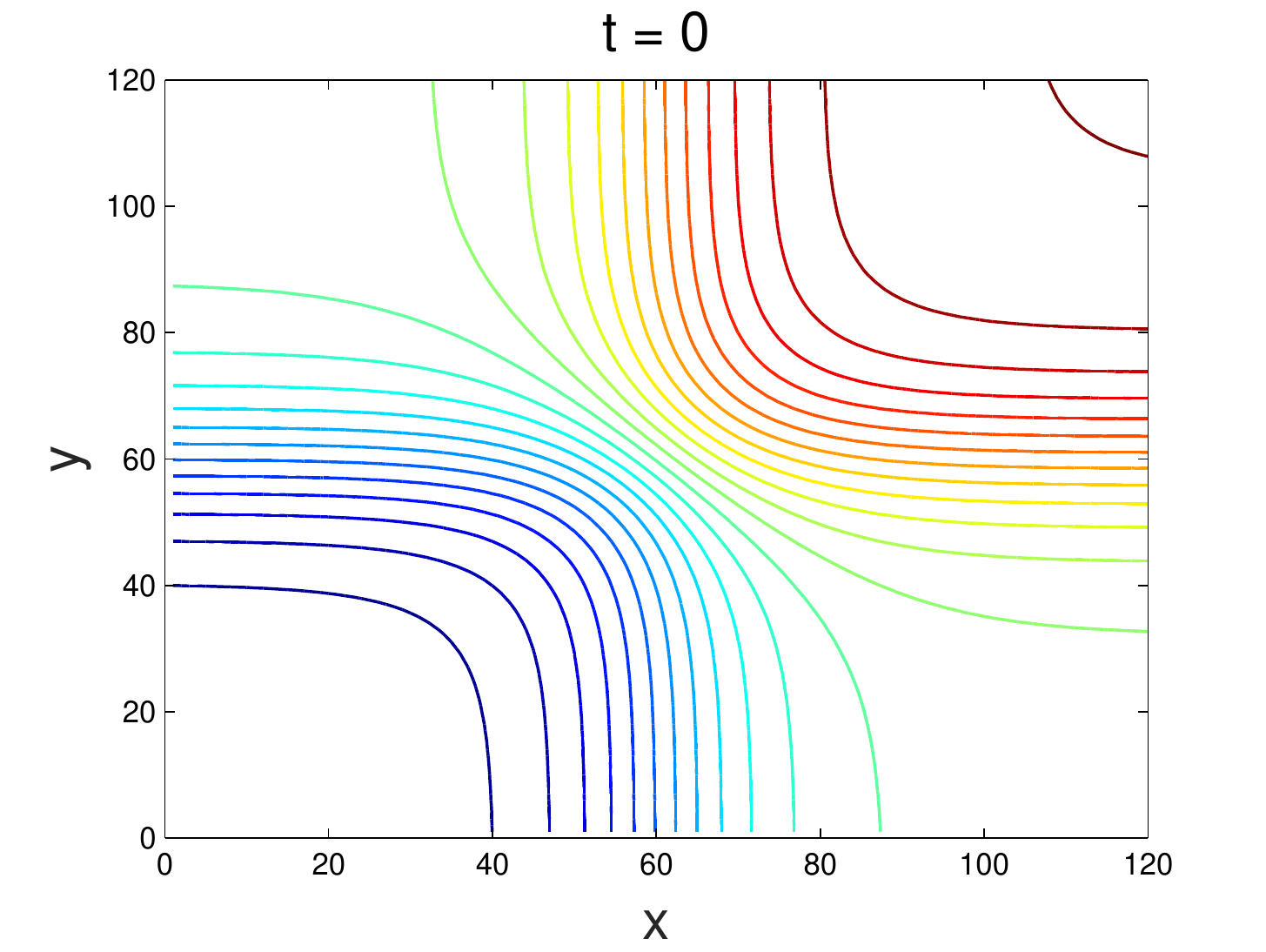}	
	\includegraphics[width=0.24\linewidth]{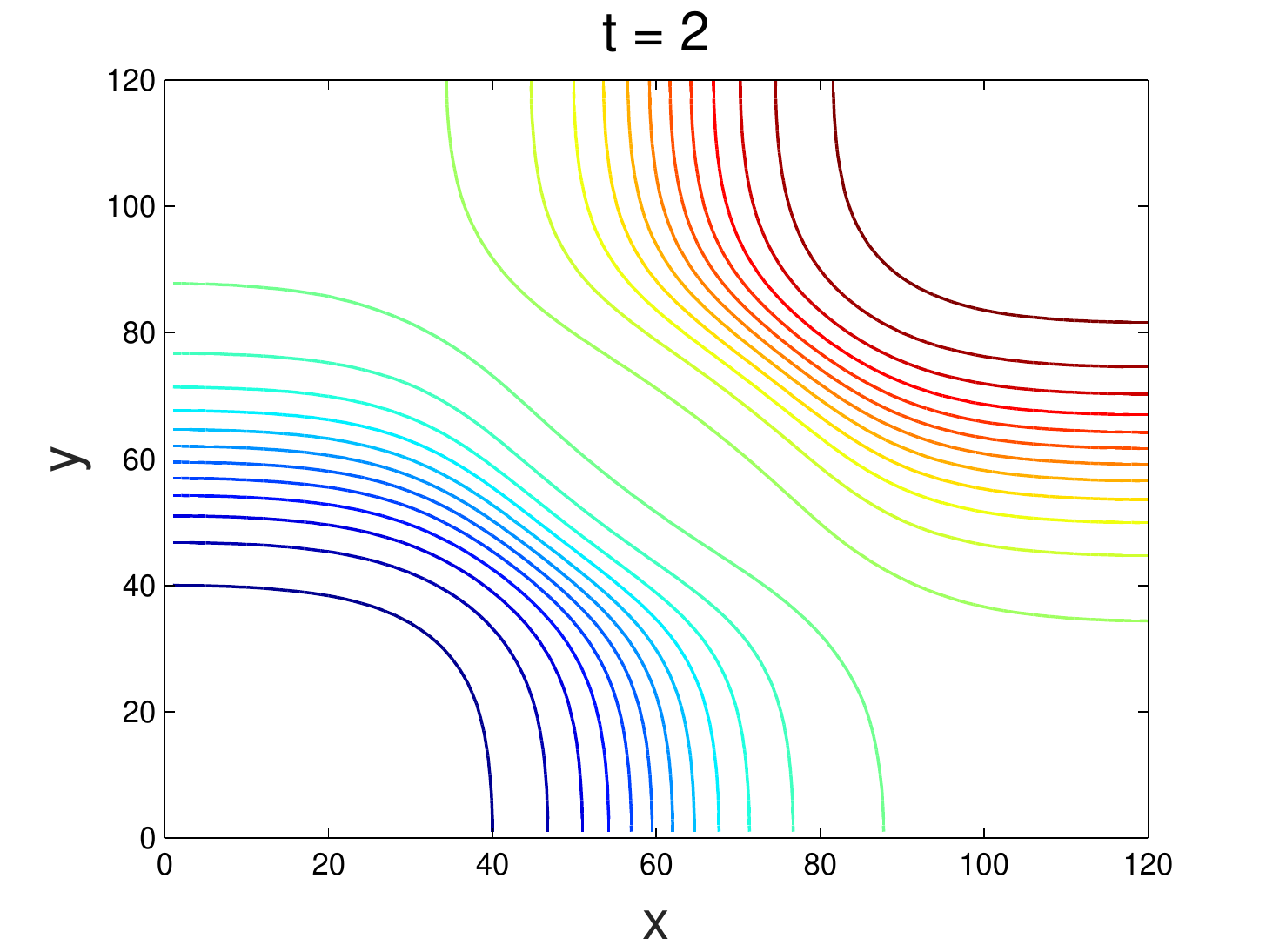}
	\includegraphics[width=0.24\linewidth]{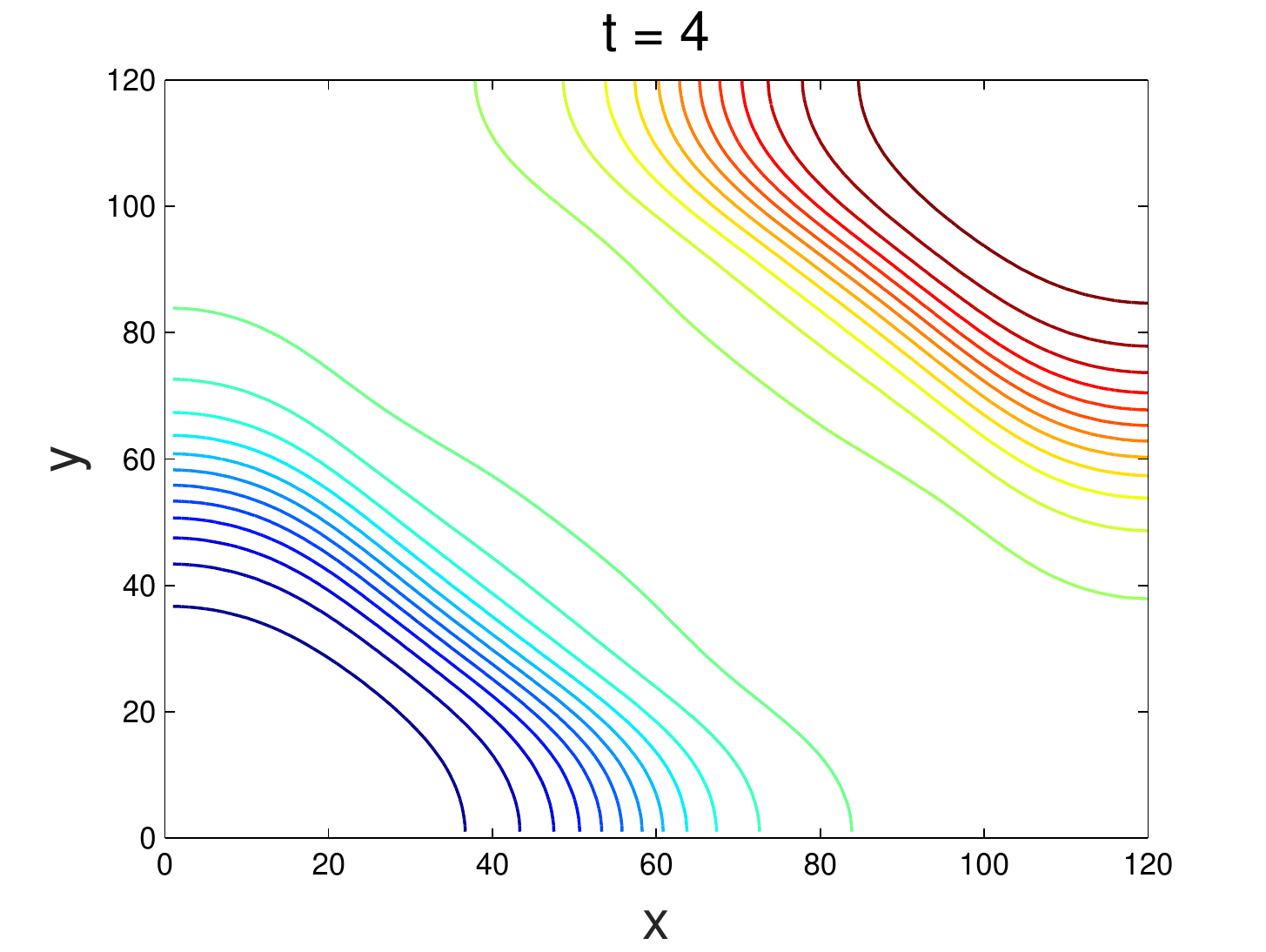}
	\includegraphics[width=0.24\linewidth]{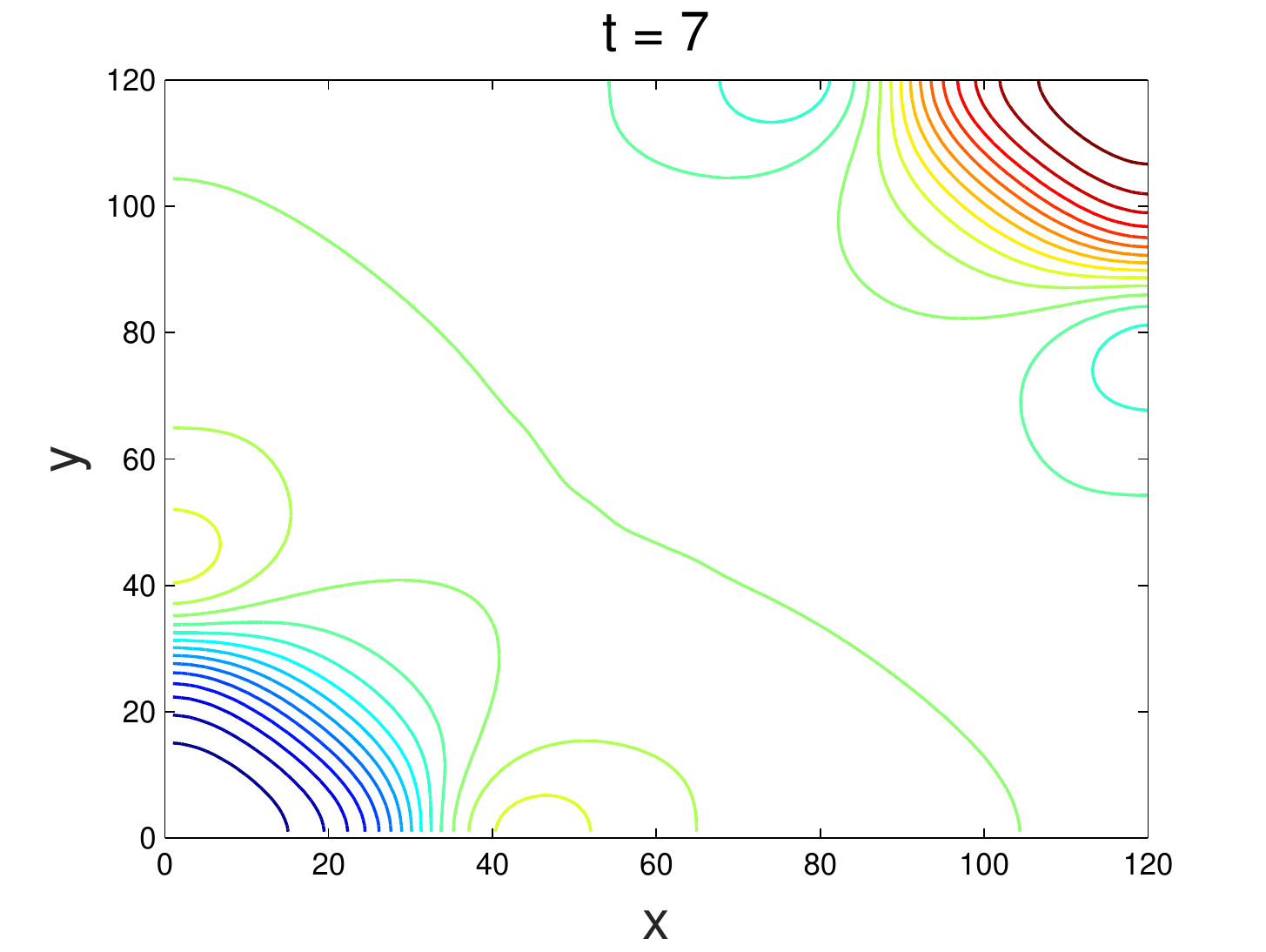}
	\caption{Superposition of two orthogonal line solitons and relevant contours: initial condition and numerical solutions at $t=2,4,7$. Spatial and temporal step sizes are taken as $h=0.5$, $\tau=0.01$.}
	\label{fig:ex2-t}
\end{figure}

Next we consider the superposition of two line solitons, which relates to the case of $\phi(x,y)=1$, with initial conditions
\[
\begin{aligned}
&f_1(x,y) = 4\tan^{-1}\exp(x)+4\tan^{-1}\exp(y),\\
&f_2(x,y) = 0,~-6\leq x,y\leq 6.
\end{aligned}
\] 
Figure~\ref{fig:ex2-t} presents the initial condition as well as numerical solutions at time $t=2,4,7$, and demonstrates that the breakup of two orthogonal line solitons which are parallel to the $y=-x$ line and moving away from each other on $y=x$ line. It can be deduced from their relevant contours that until $t=4$ this separation occurs without any deformation, while at $t=7$ a deformation has already appeared. All solutions are in complete agreement with the corresponding results in \cite{cl81}. 

\begin{figure}[H]
	\centering
	\includegraphics[width=0.45\linewidth]{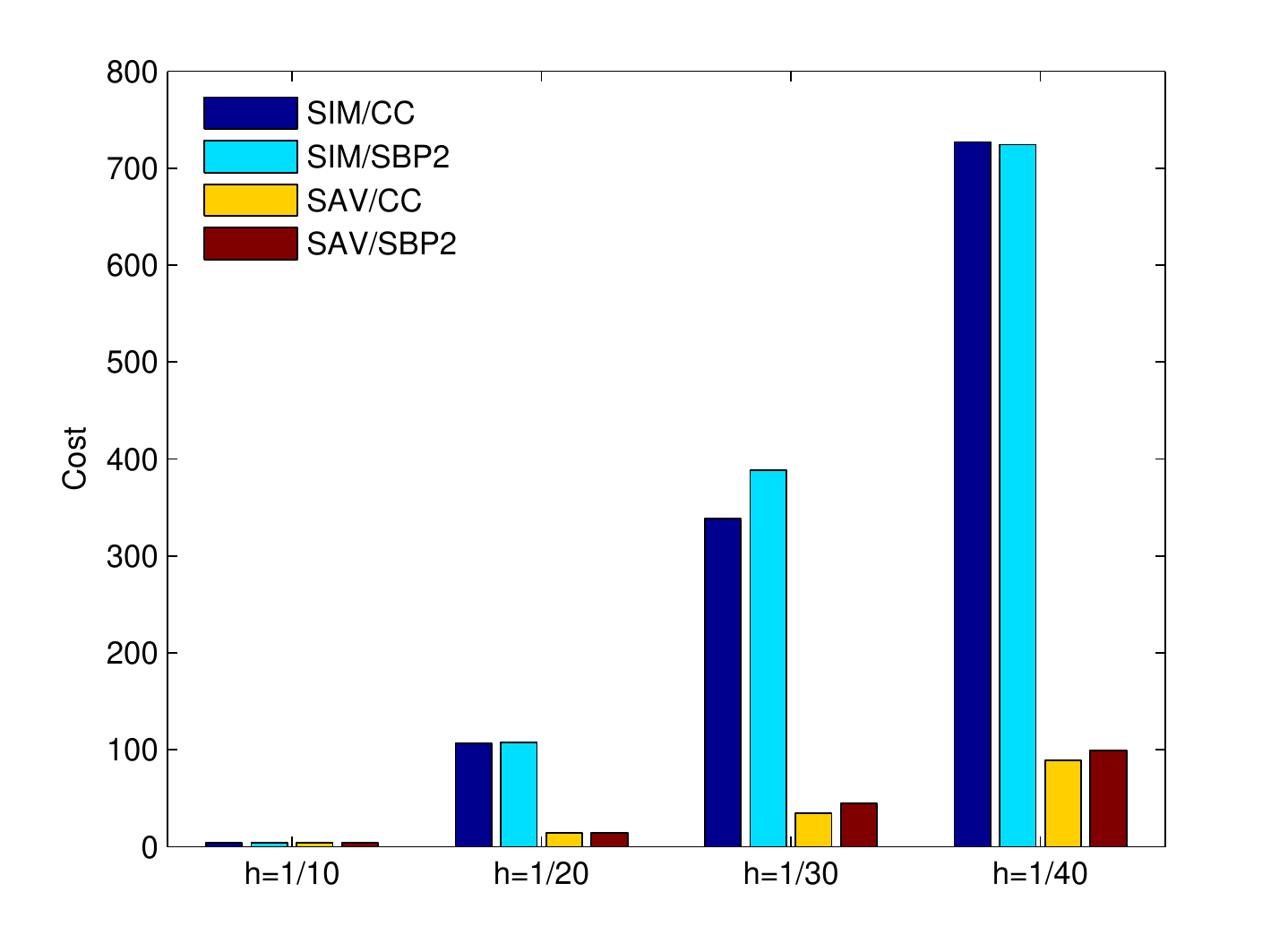}
	\caption{Computational cost of the proposed schemes for the superposition of two line solitons with different mesh size till $t=7$.}
	\label{fig:ex2-cost}
\end{figure}

As demonstrated in Section 4, the energy-preserving schemes based on the SAV approach are not only linearly implicit, but also induce algebraic systems with constant coefficient matrix which can be efficiently solved. While the two symplectic integrators {\bf SIM/CC} and {\bf SIM/SBP2} are fully implicit and need nonlinear iterations to obtain the numerical solution, which can approximately represent the general implicit method, regarding the computational cost. In Figure~\ref{fig:ex2-cost}, we present comparisons on the computational costs among the four schemes by refining the mesh grid gradually. The biconjugate gradient stabilized method (named {\em bicgstab} in MATLAB functions) is adopt as the solver of linear systems given by schemes {\bf SAV/CC} and {\bf SAV/SBP2}. It is noticed that when we take a relatively small grid size, i.e. $h=0.5$, the costs of two fully implicit schemes {\bf SIM/CC} and {\bf SIM/SBP2} are actually considerable with that of the two linearly implicit schemes. However, as the refinement of mesh grid, the advantage of scheme {\bf SAV/CC} and {\bf SAV/SBP2} emerges. From the bar plot in Figure~\ref{fig:ex2-cost}, the costs of these two schemes are about 10 times cheaper than that of schemes {\bf SIM/CC} and {\bf SIM/SBP2}. Therefore, it is preferable to construct linearly implicit schemes through the SAV approach for large scale simulations, keeping the system energy being preserved as well.

\begin{figure}[H]
	\centering
	\includegraphics[width=0.45\linewidth]{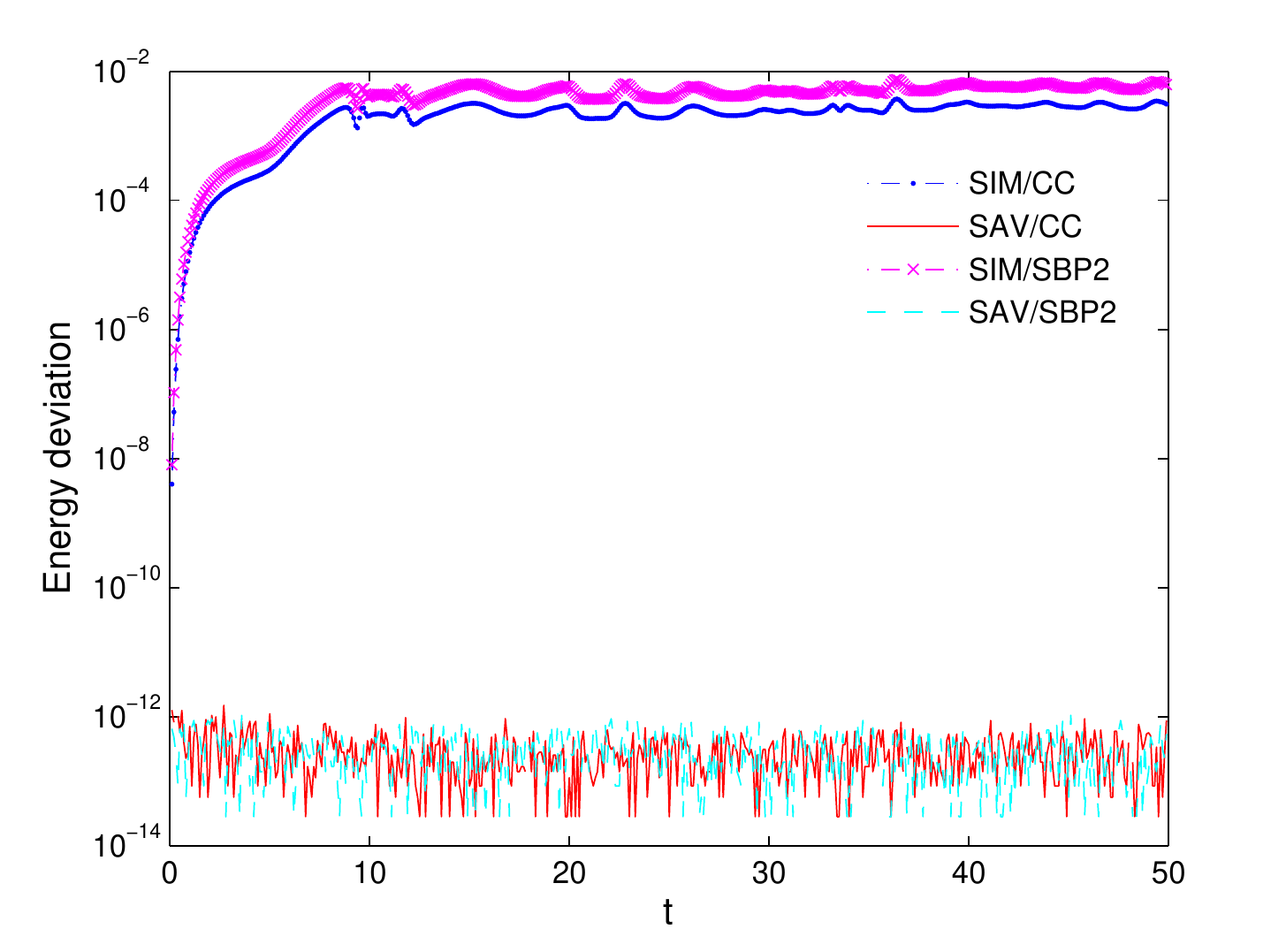}
	\caption{Long time energy deviation of the proposed schemes for the superposition of two line solitons.}
	\label{fig:ex2-ene}
\end{figure}

We also run a long time simulation till $t=50$ and plot the energy deviation in Figure~\ref{fig:ex2-ene}, corresponding to all four schemes, which shows that the SAV approach based schemes {\bf SAV/CC} and {\bf SAV/SBP2} preserve the discrete energy to round-off errors. Although the two symplectic integrators {\bf SIM/CC} and {\bf SIM/SBP2} cannot preserve the energy exactly, the deviations always remain oscillated around a small order of magnitude. 

\subsection{Perturbation of a static line soliton}

\begin{figure}[H]
	\centering
	\includegraphics[width=0.32\linewidth]{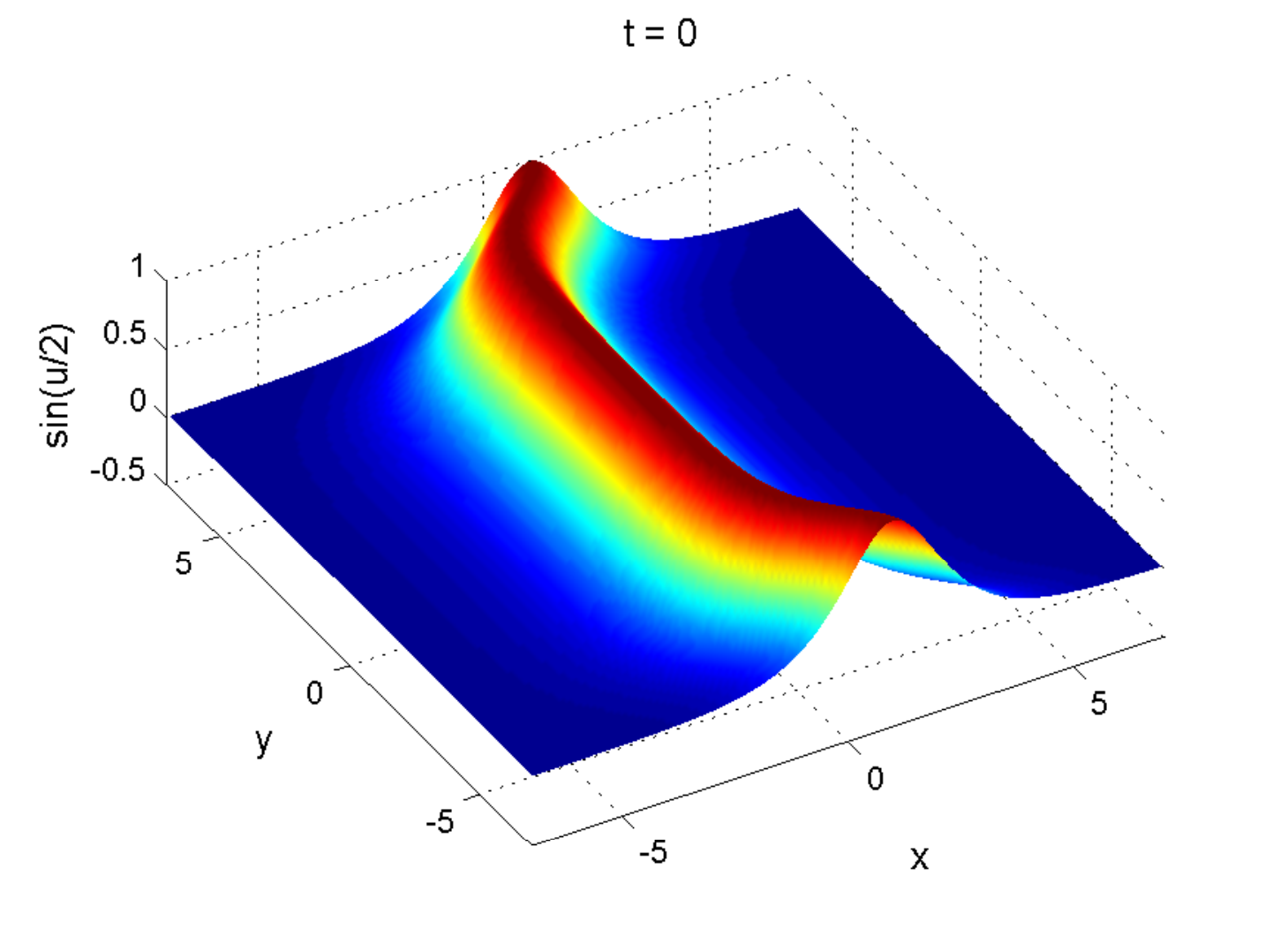}
	\includegraphics[width=0.32\linewidth]{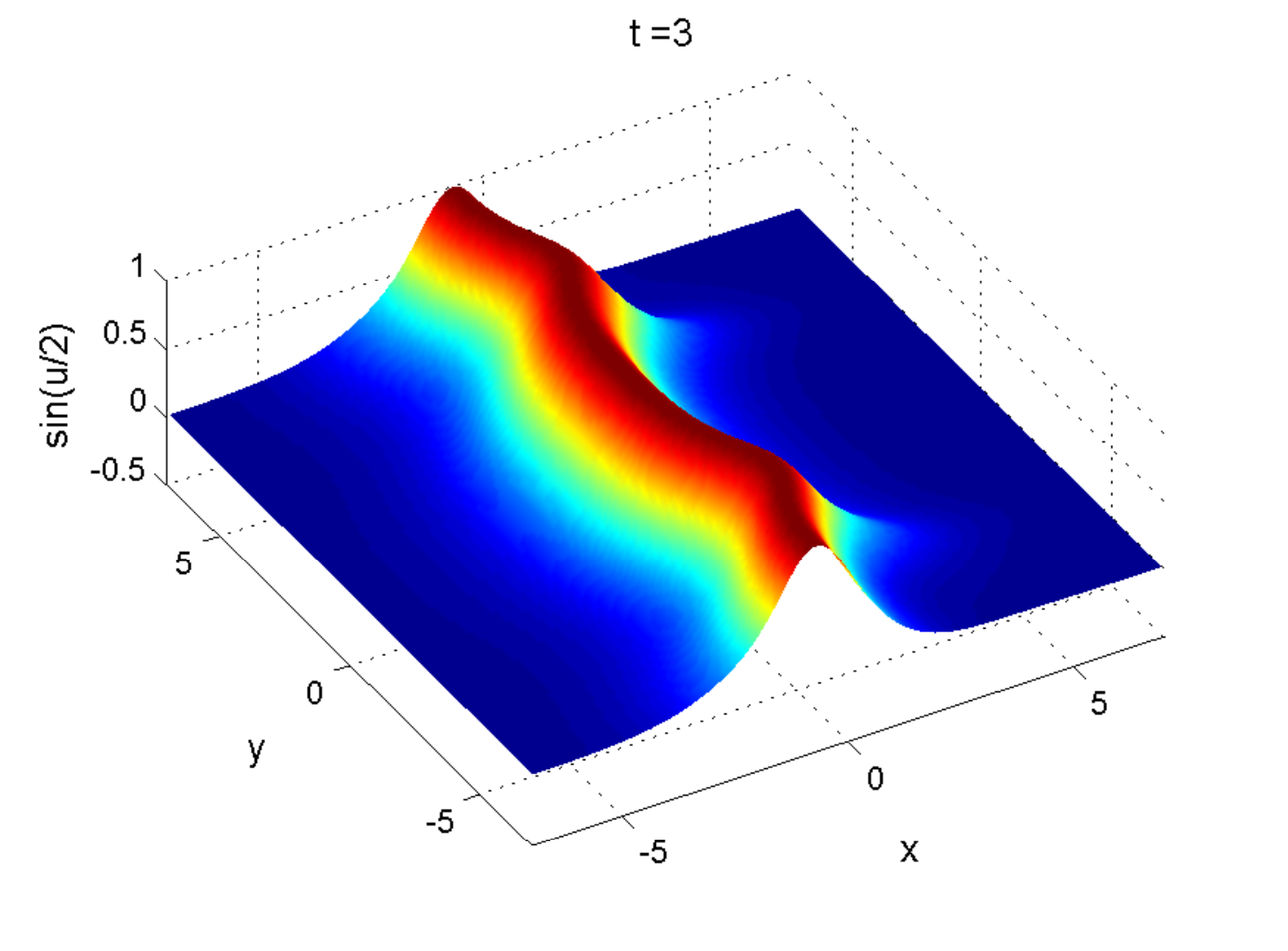}
	\includegraphics[width=0.32\linewidth]{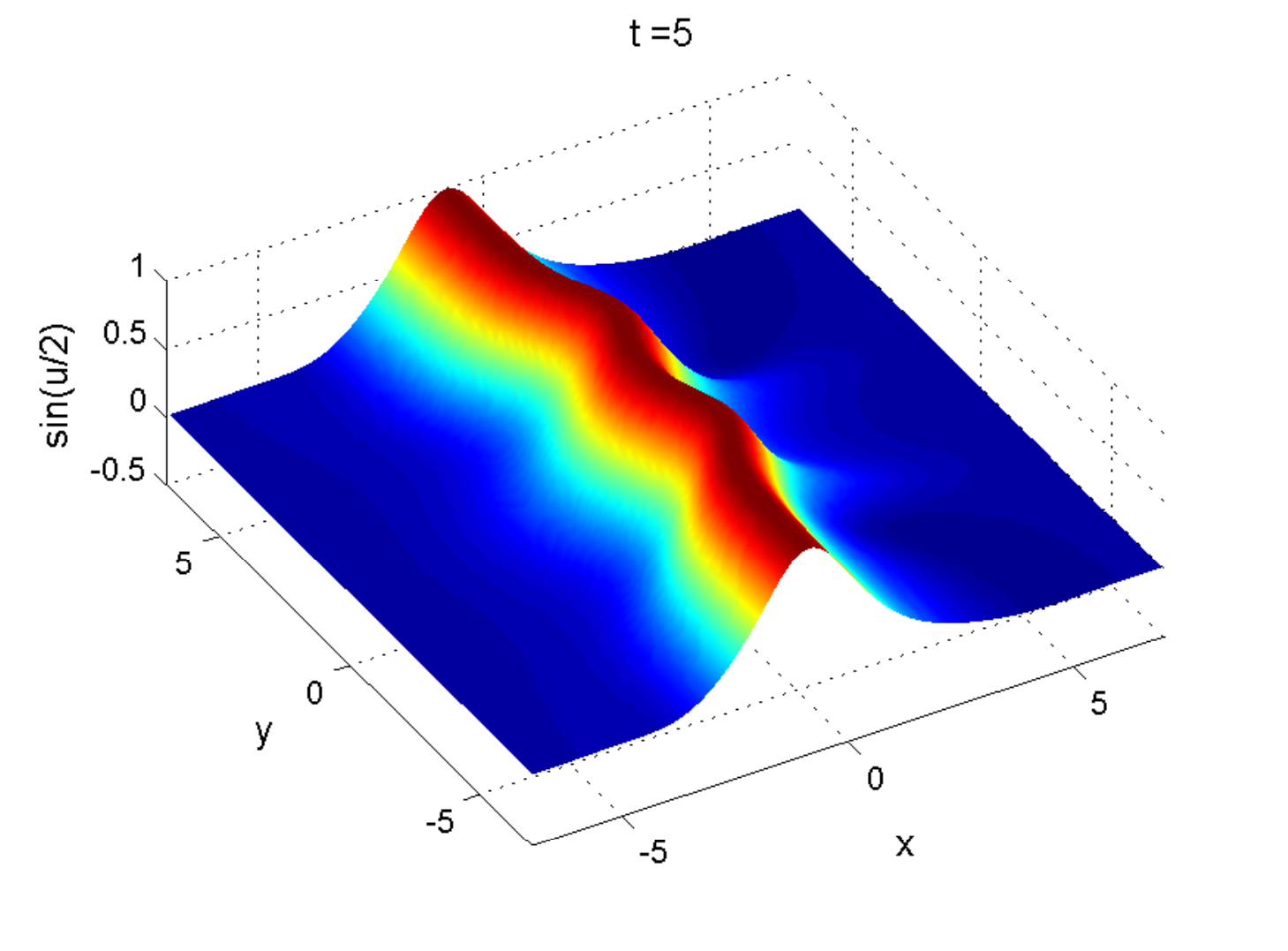}\\
	\includegraphics[width=0.32\linewidth]{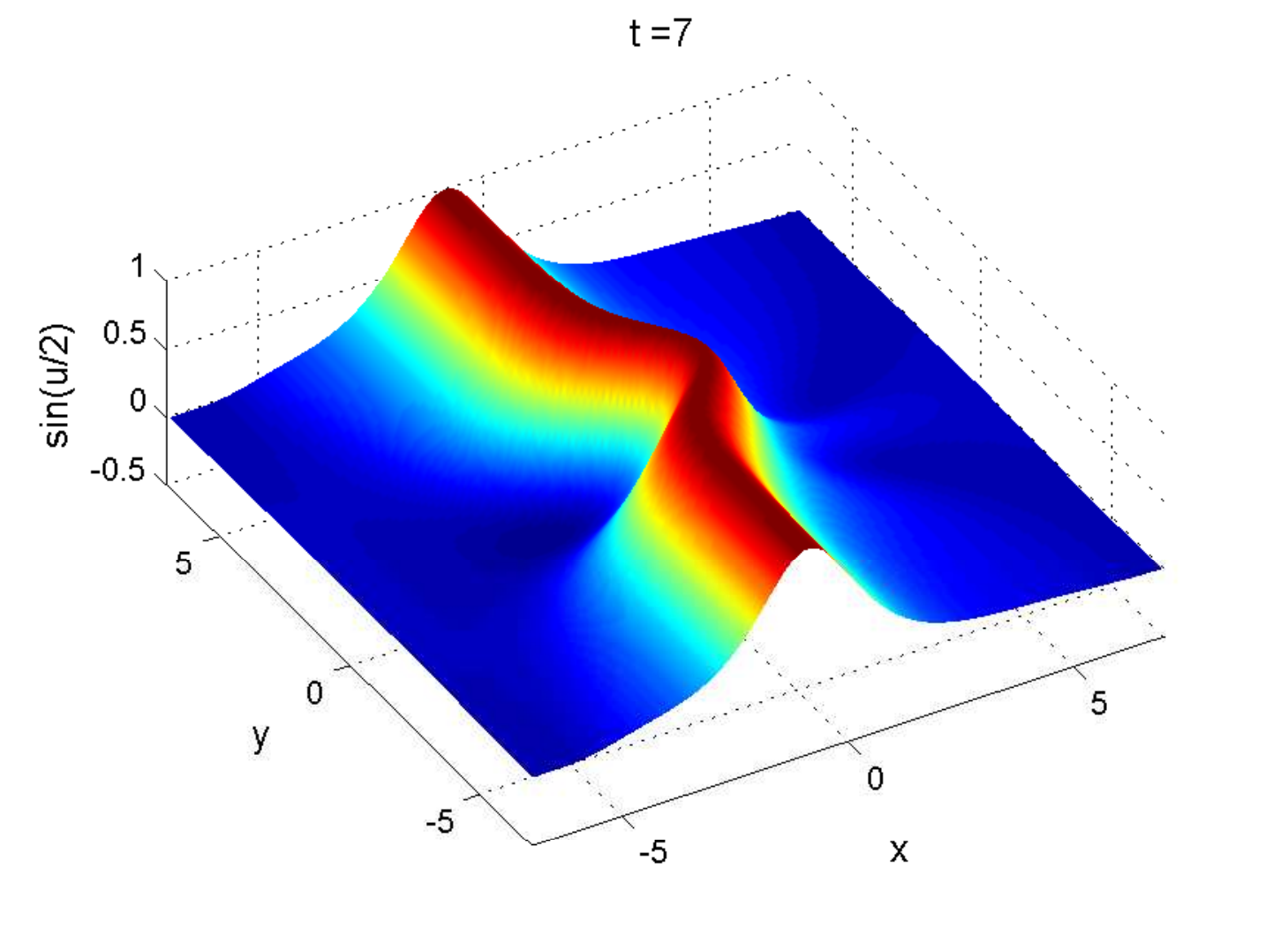}
	\includegraphics[width=0.32\linewidth]{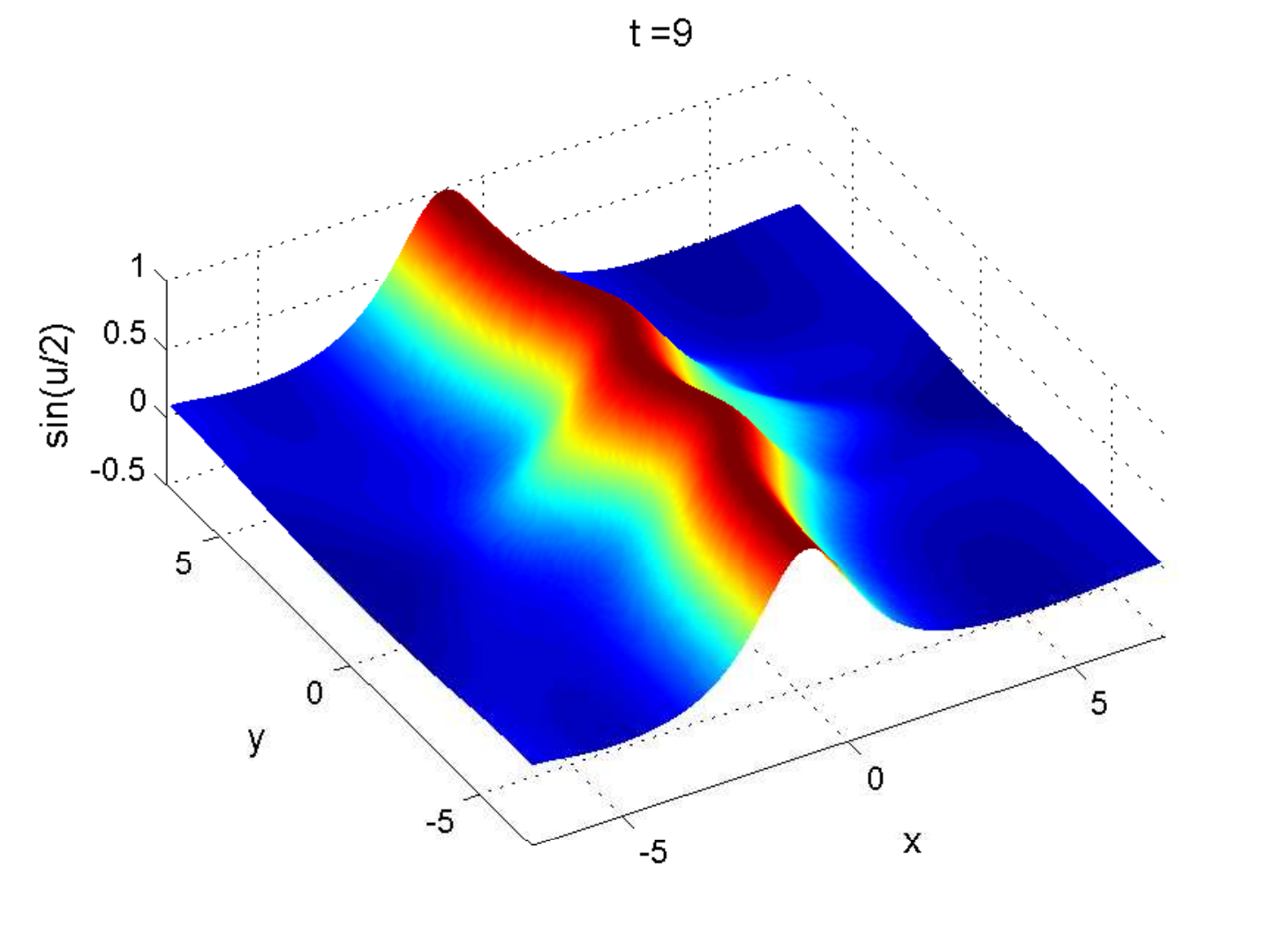}
	\includegraphics[width=0.32\linewidth]{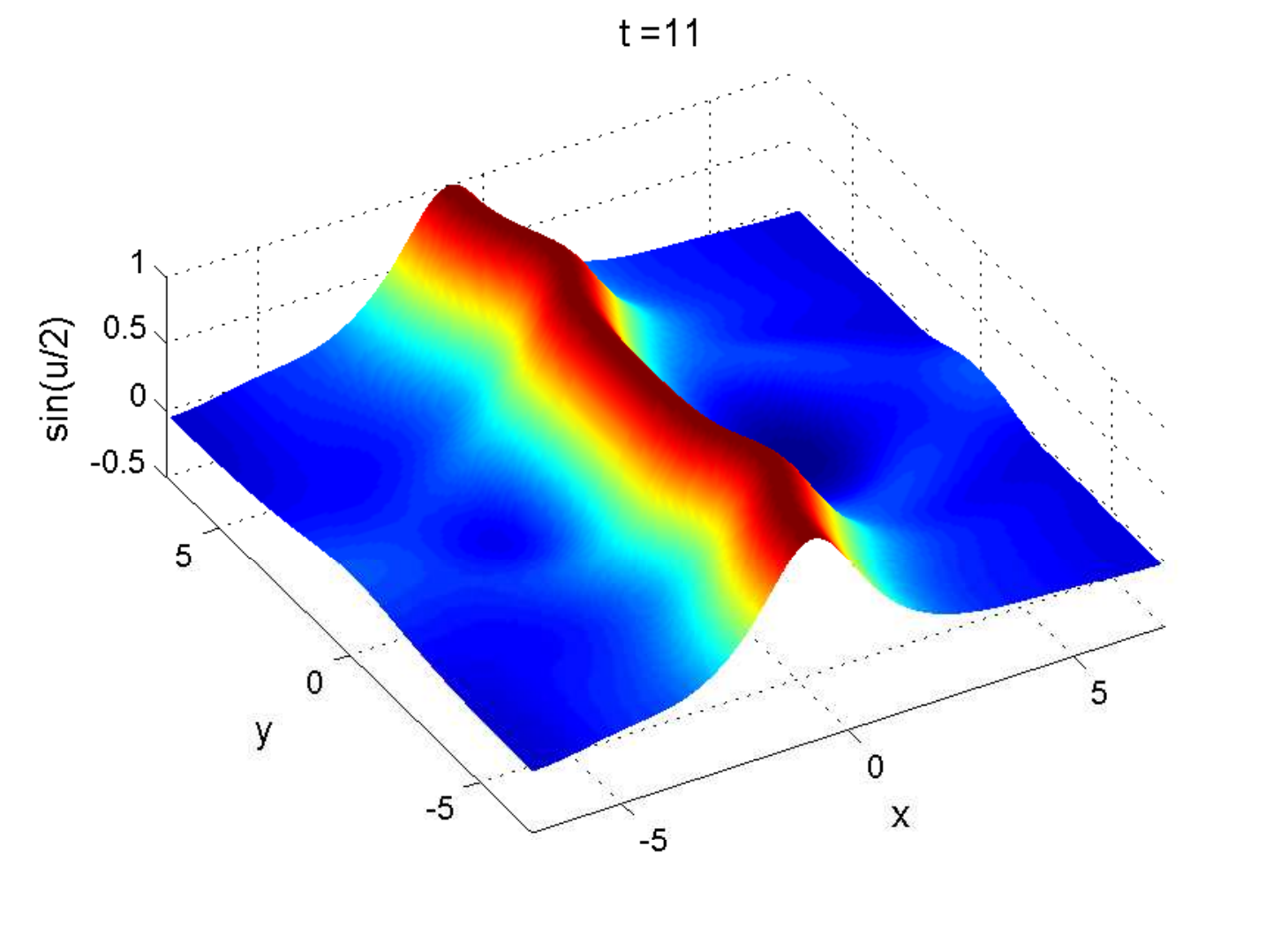}
	\caption{Perturbation of a static line soliton: initial condition and numerical solutions at $t = 3,5,7,9,11$, in terms of  $\sin(u/2)$. Spatial and temporal step sizes are taken as $h=0.5$, $\tau=0.01$.}
	\label{fig:ex3-t}
\end{figure}

\begin{figure}[H]
	\centering
	\includegraphics[width=0.45\linewidth]{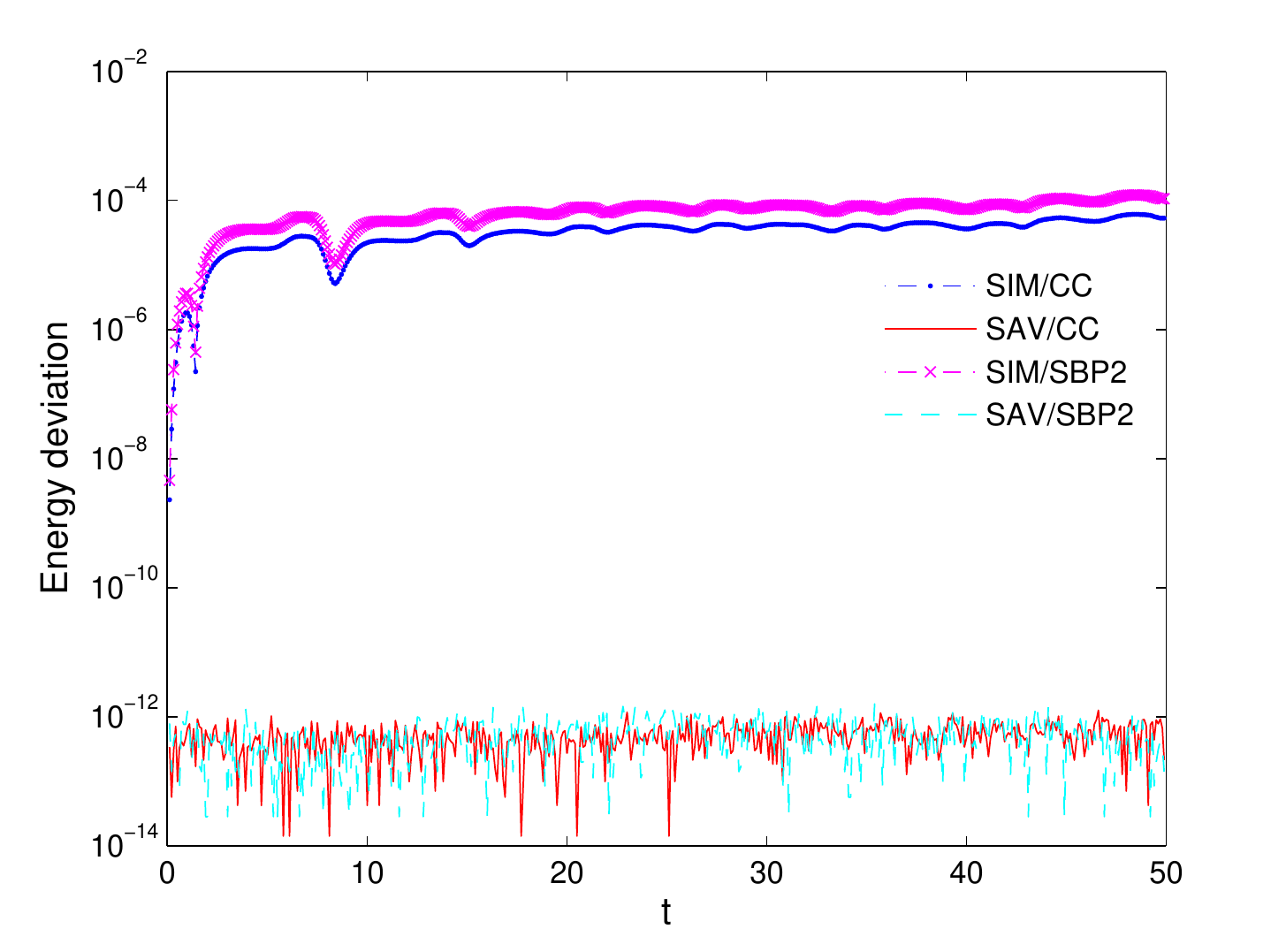}
	\caption{Long time energy deviation of the proposed schemes for the perturbation of a static line soliton.}
	\label{fig:ex3-ene}
\end{figure}

A static line soliton is perturbed to produce two symmetric dents moving towards each other with a constant unit velocity. These dents collide and continue to travel at the same velocity without any shift. The waves are defined by $\phi(x,y)=1$ and the initial conditions
\[
\begin{aligned}
&f_1(x,y) = 4\tan^{-1}\exp[x+1-2\mbox{sech}(y+7)-2\mbox{sech}(y-7)],\\
&f_2(x,y) = 0,~ -7\leq x,y\leq 7.
\end{aligned}
\]
The results in Figure~\ref{fig:ex3-t} show that two symmetric dents moving toward each other, collapsing at $t=7$ and continuing to move away from each other. It can be deduced that after
the collision these dents retain their shape which verifies the conclusions of Christiansen and Lomdahl \cite{cl81}. The energy deviations of all four schemes are also given in Figure~\ref{fig:ex3-ene} which demonstrate a consistent result as that in the superposition of two line solitons.

\subsection{Line soliton in an inhomogeneous medium}

\begin{figure}[H]
	\centering
	\includegraphics[width=0.24\linewidth]{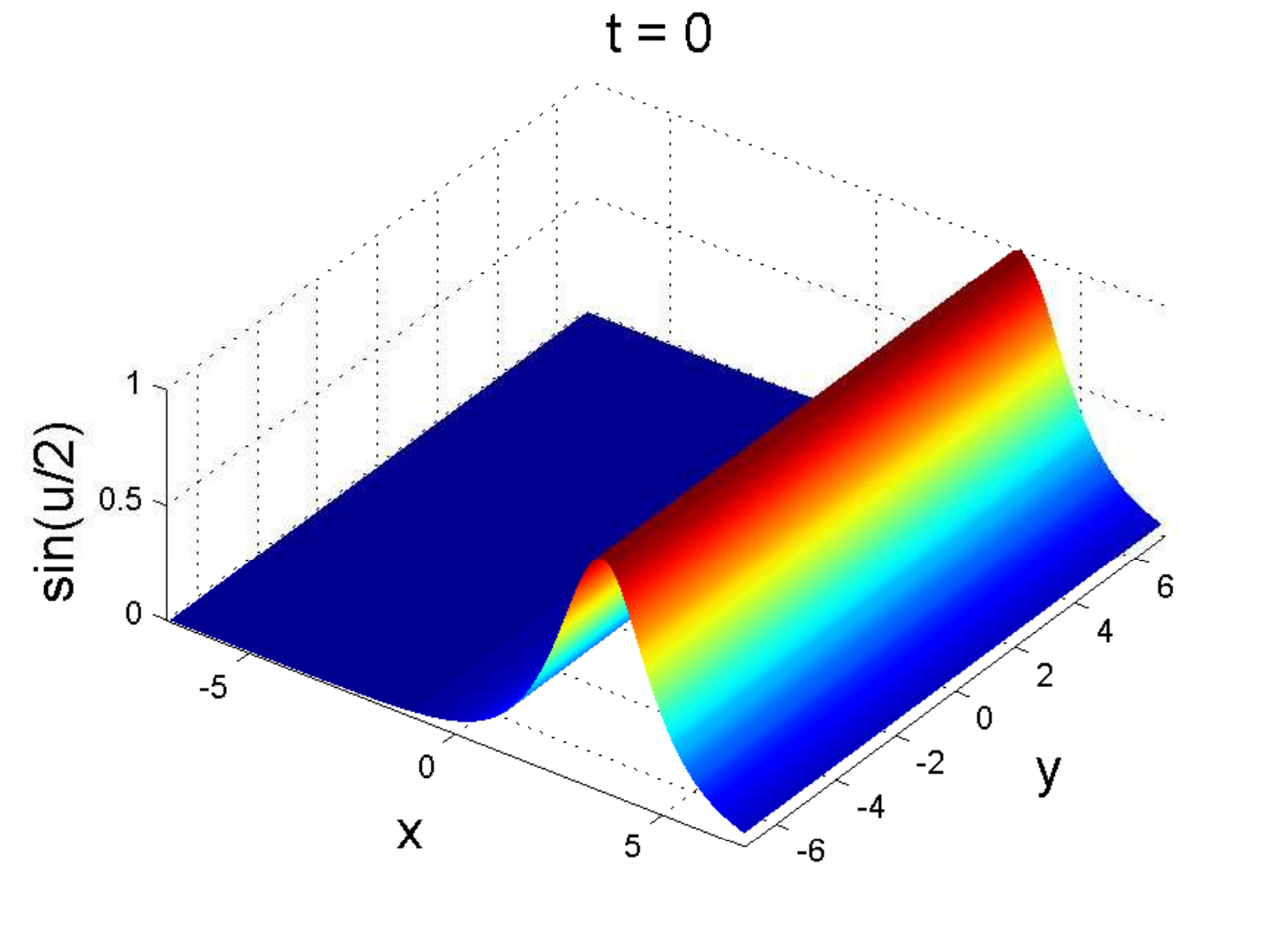}
	\includegraphics[width=0.24\linewidth]{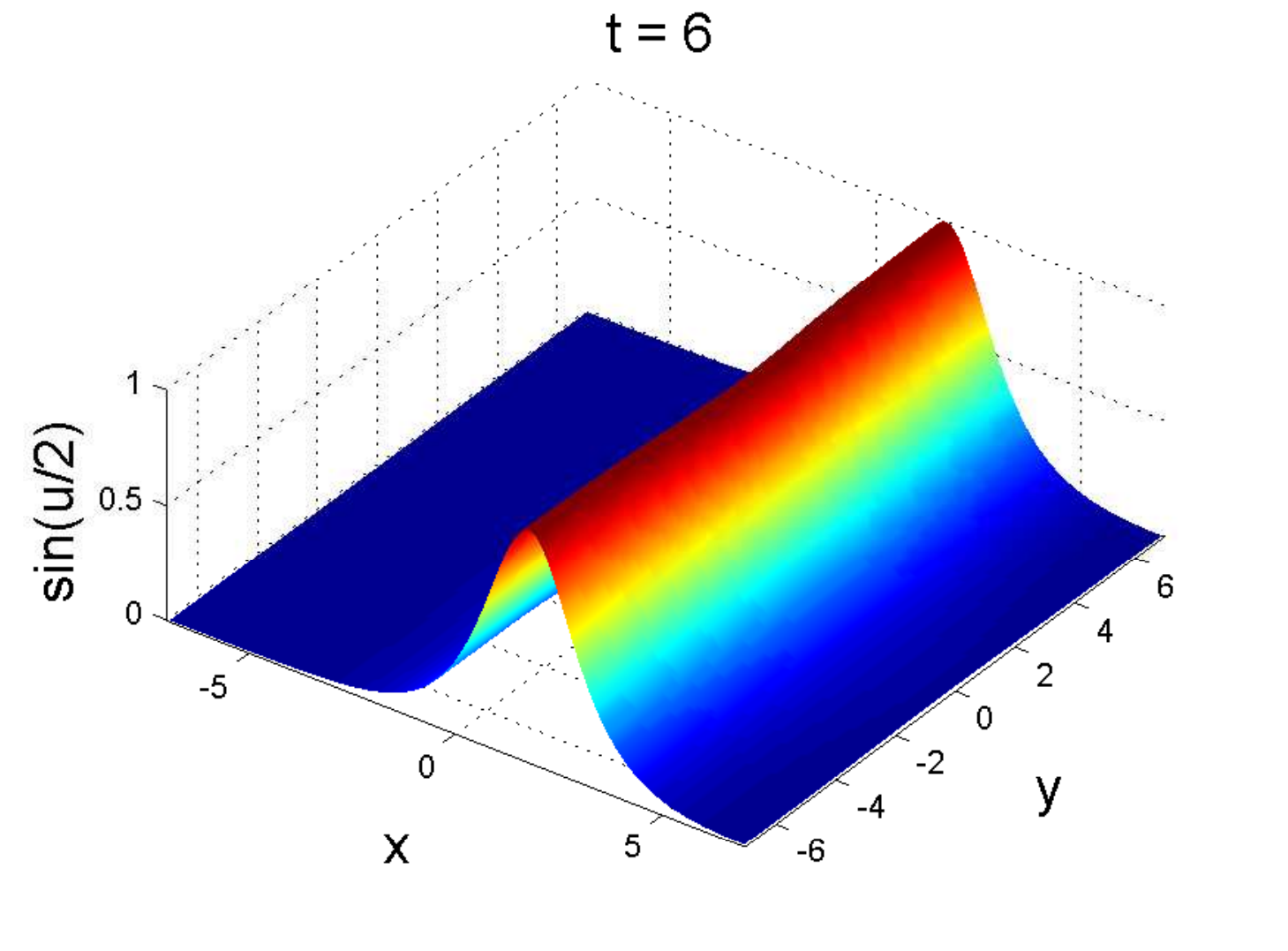}
	\includegraphics[width=0.24\linewidth]{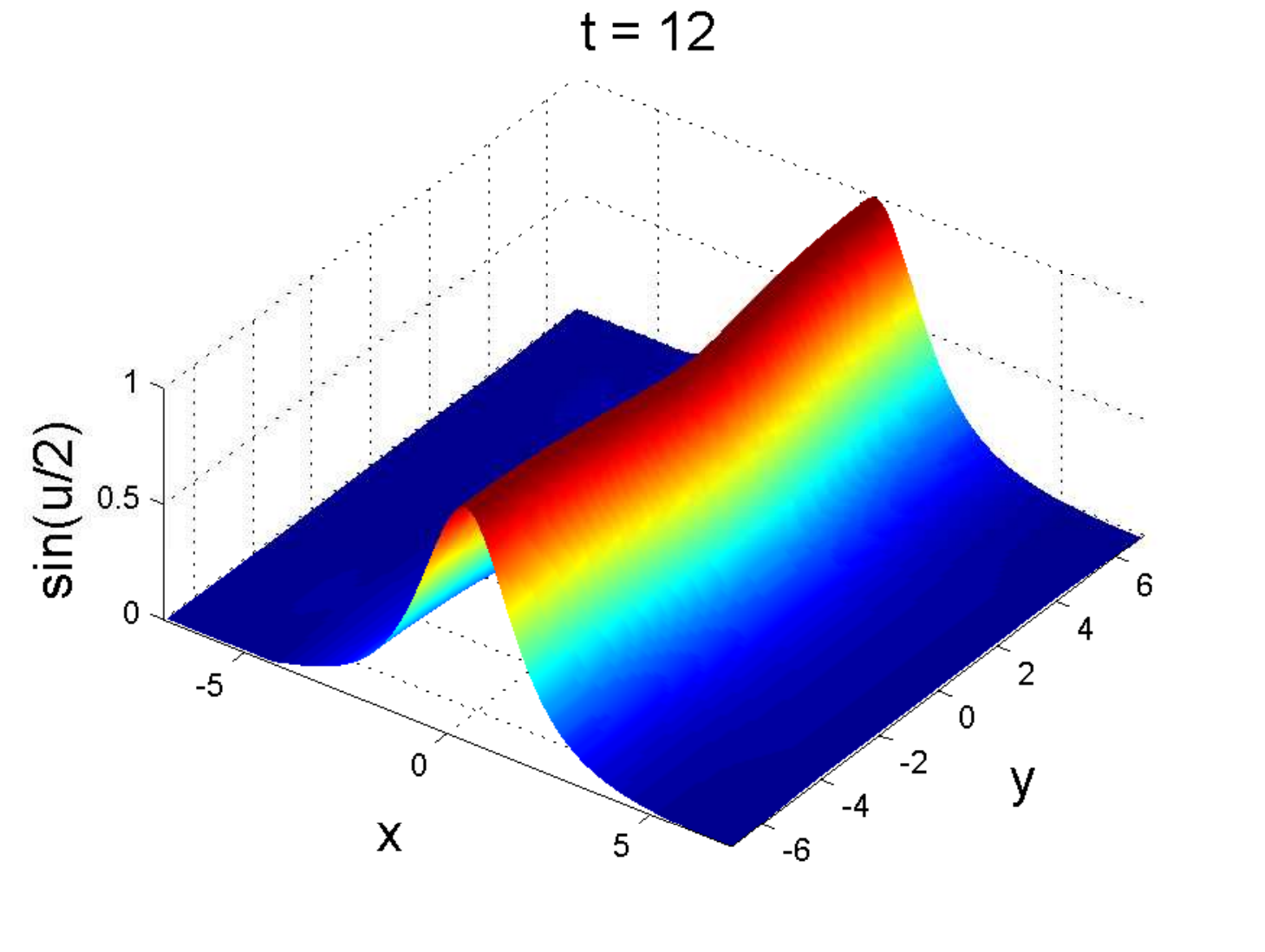}
	\includegraphics[width=0.24\linewidth]{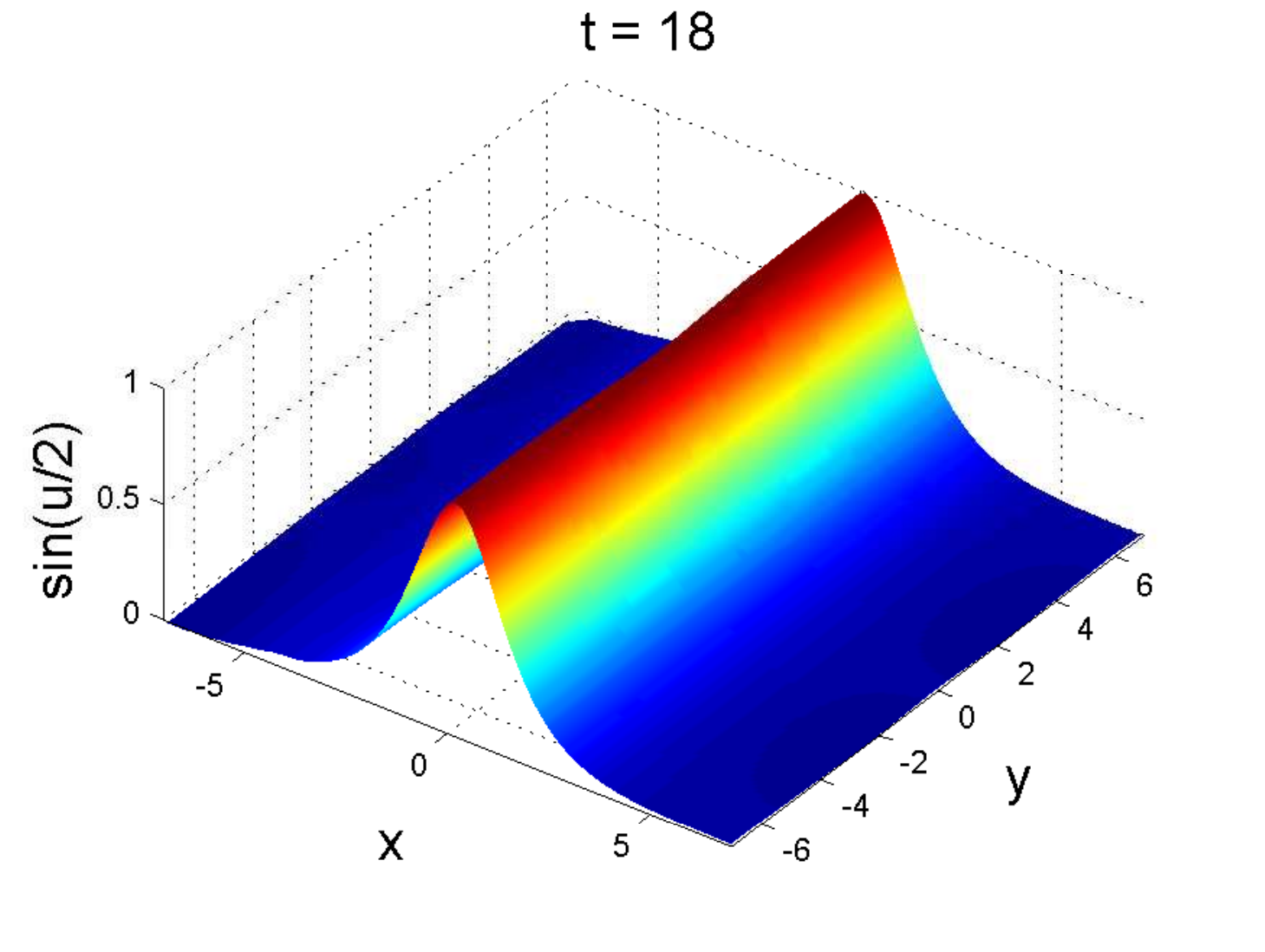}\\[2ex]
	\includegraphics[width=0.24\linewidth]{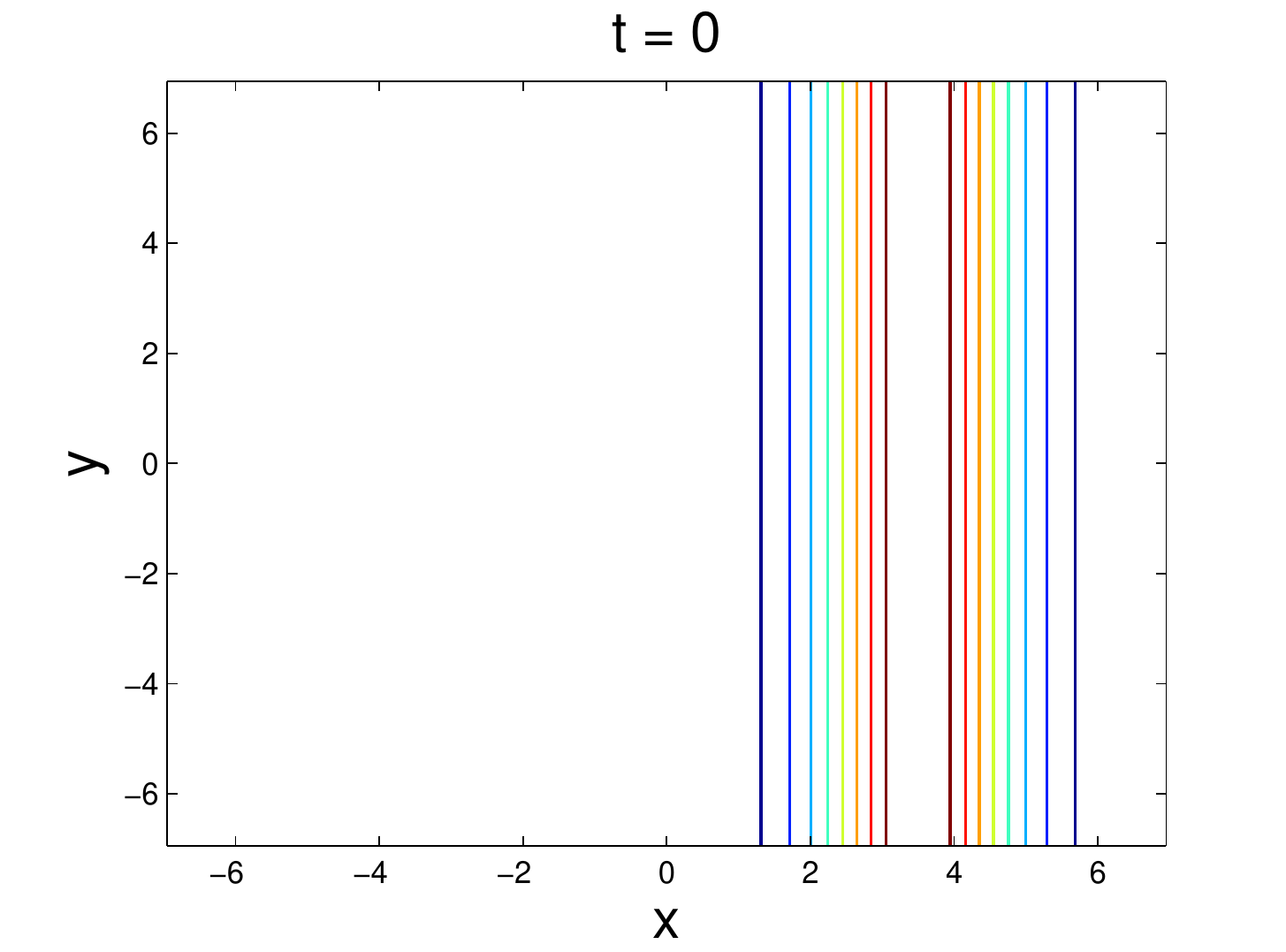}
	\includegraphics[width=0.24\linewidth]{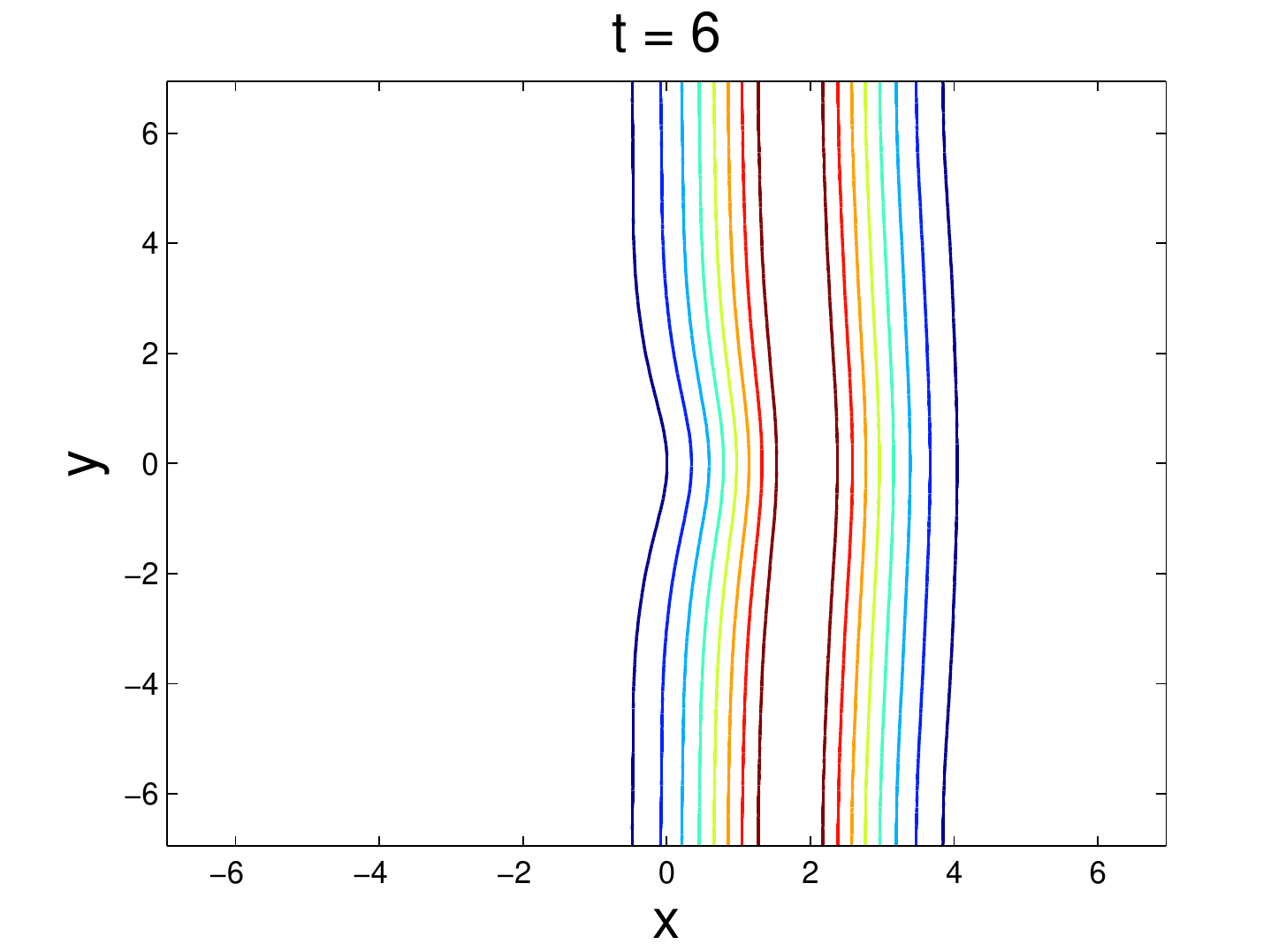}
	\includegraphics[width=0.24\linewidth]{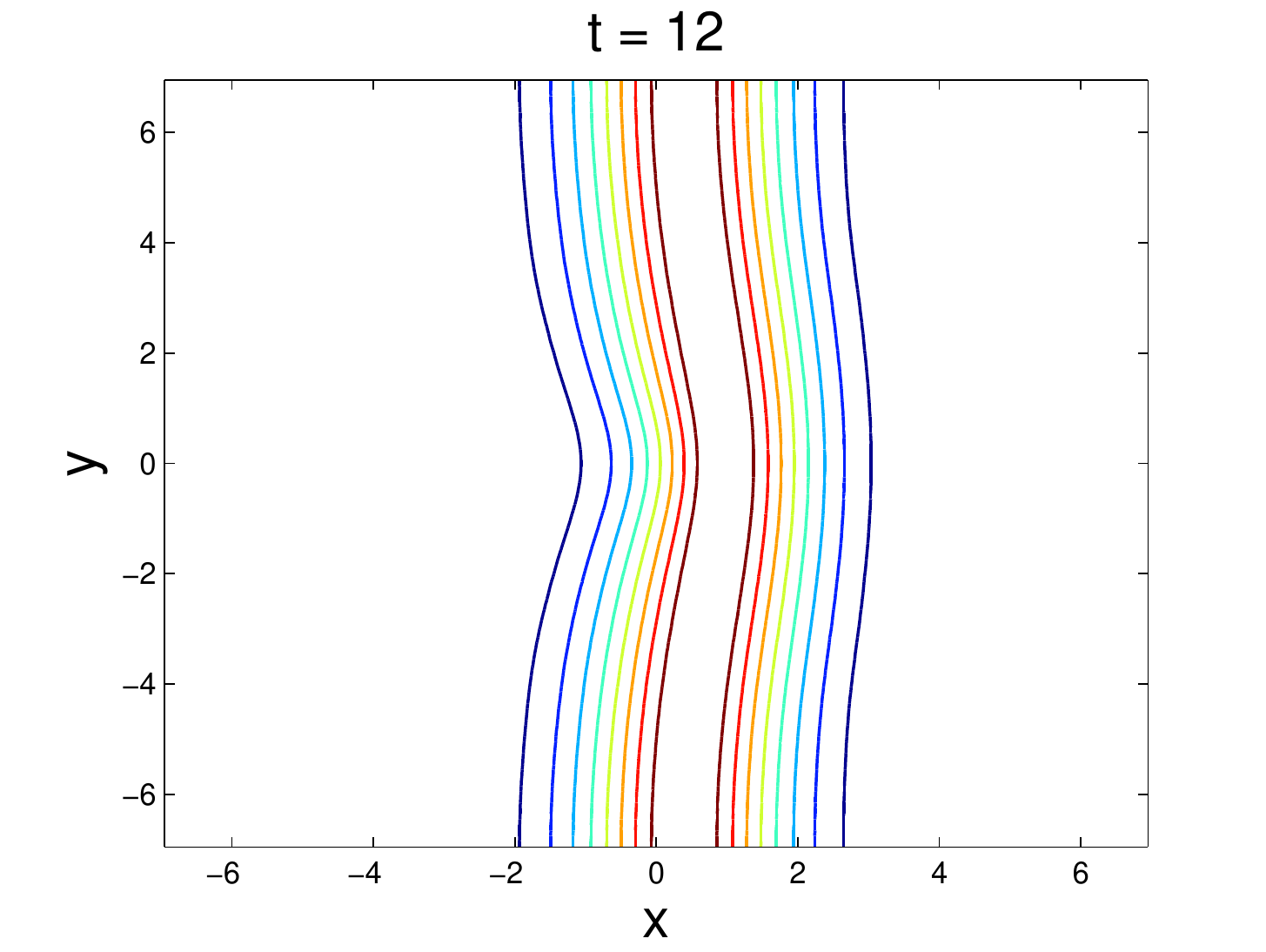}
	\includegraphics[width=0.24\linewidth]{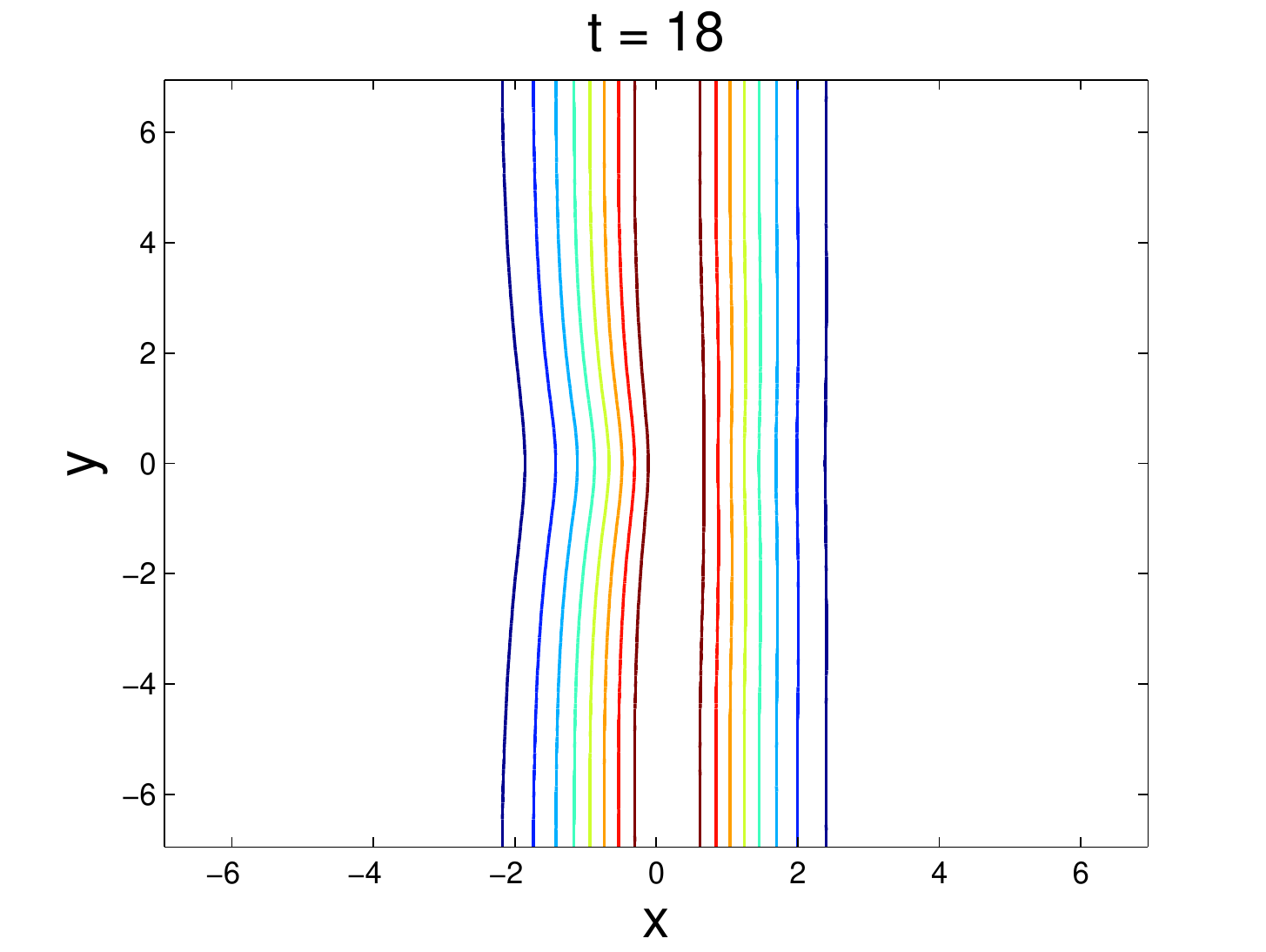}
	\caption{Line soliton in an inhomogeneous medium and relevant contours: initial condition and numerical solutions at $t =6, 12, 18$, in terms of  $\sin(u/2)$. Spatial and temporal step sizes are taken as $h=0.5$, $\tau=0.01$.}
	\label{fig:ex33-t}
\end{figure}

\begin{figure}[H]
	\centering
	\includegraphics[width=0.45\linewidth]{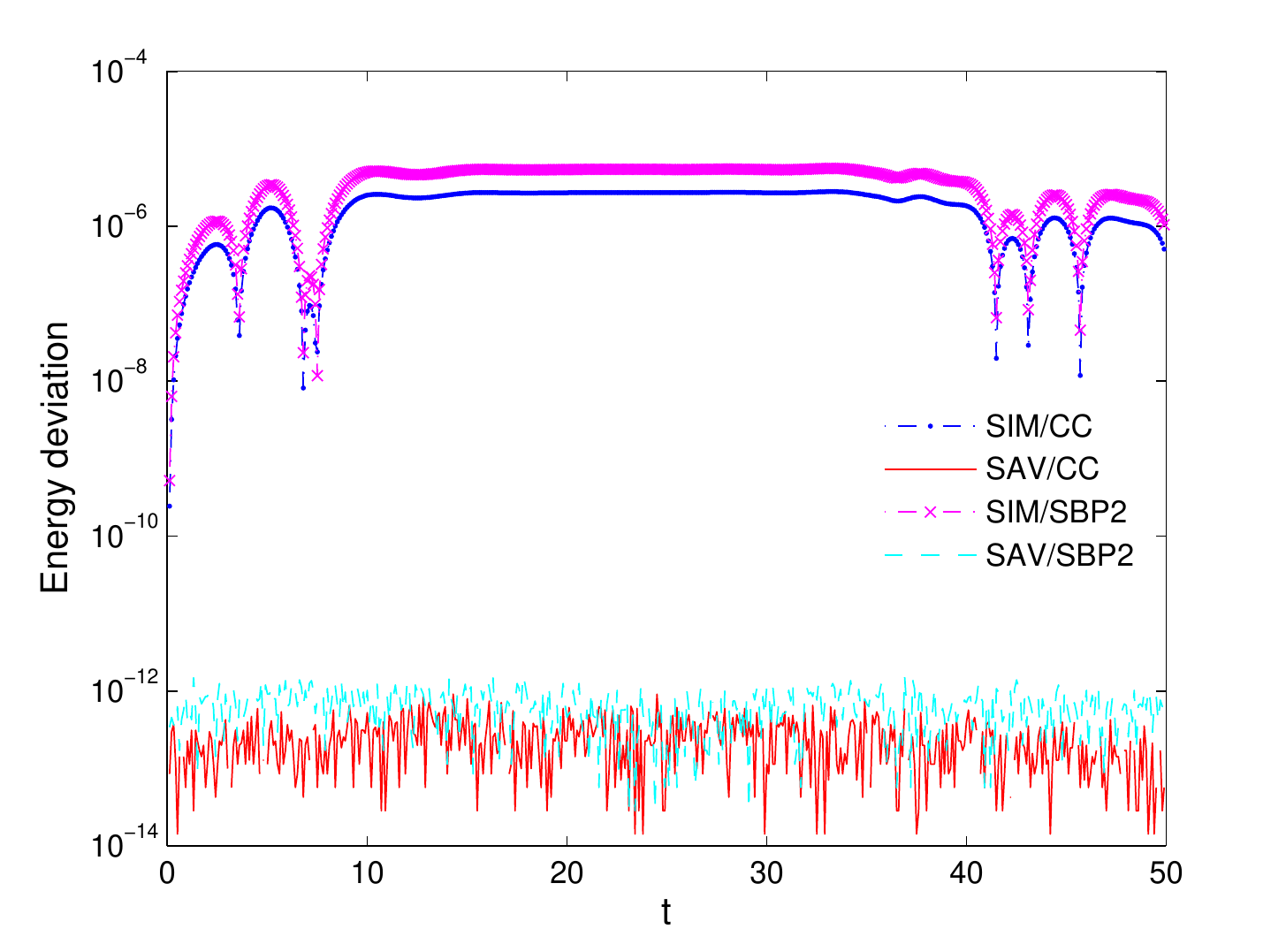}
	\caption{Long time energy deviation of the proposed schemes for the line soliton in an inhomogeneous medium.}
	\label{fig:ex33-ene}
\end{figure}

A model for an inhomogeneous on large-area Josephson junction is given by the Josephson current density 
\[
\phi(x,y)=1+\mbox{sech}^2\sqrt{x^2+y^2},
\]
with initial conditions
\[
\begin{aligned}
&f_1(x,y) = 4\tan^{-1}\exp\Big[\frac{x-3.5}{0.954}\Big],\\
&f_2(x,y) = 0.629\;\mbox{sech}\Big[\frac{x-3.5}{0.954}\Big],~ -7\leq x,y\leq 7.
\end{aligned}
\]
The numerical solutions are plotted at $t=0,6,12,18$ as shown in Figure~\ref{fig:ex33-t}. The line soliton is moving in the direction $x$ through the inhomogeneity. As $t$ tends to 6 a deformation in
its straightness appears. For $t\in[6, 12]$ and as $t$ tends to 12, because of the inhomogeneity, this movement seems to be prevented, while for $t\in[12, 18]$ it is stopped. Finally, when $t =18$ the soliton almost recovers its straightness. Analogous phenomenon was observed in the results of \cite{cl81,ahh91,b07}. Despite the inhomogeneity, our schemes still preserve the discrete energy very well. From Figure~\ref{fig:ex33-ene}, we can verify that the energy deviations have no affect by the inhomogeneity and maintain the similar evolution tendency.

\subsection{Circular ring soliton}

\begin{figure}[H]
	\centering
	\includegraphics[width=0.32\linewidth]{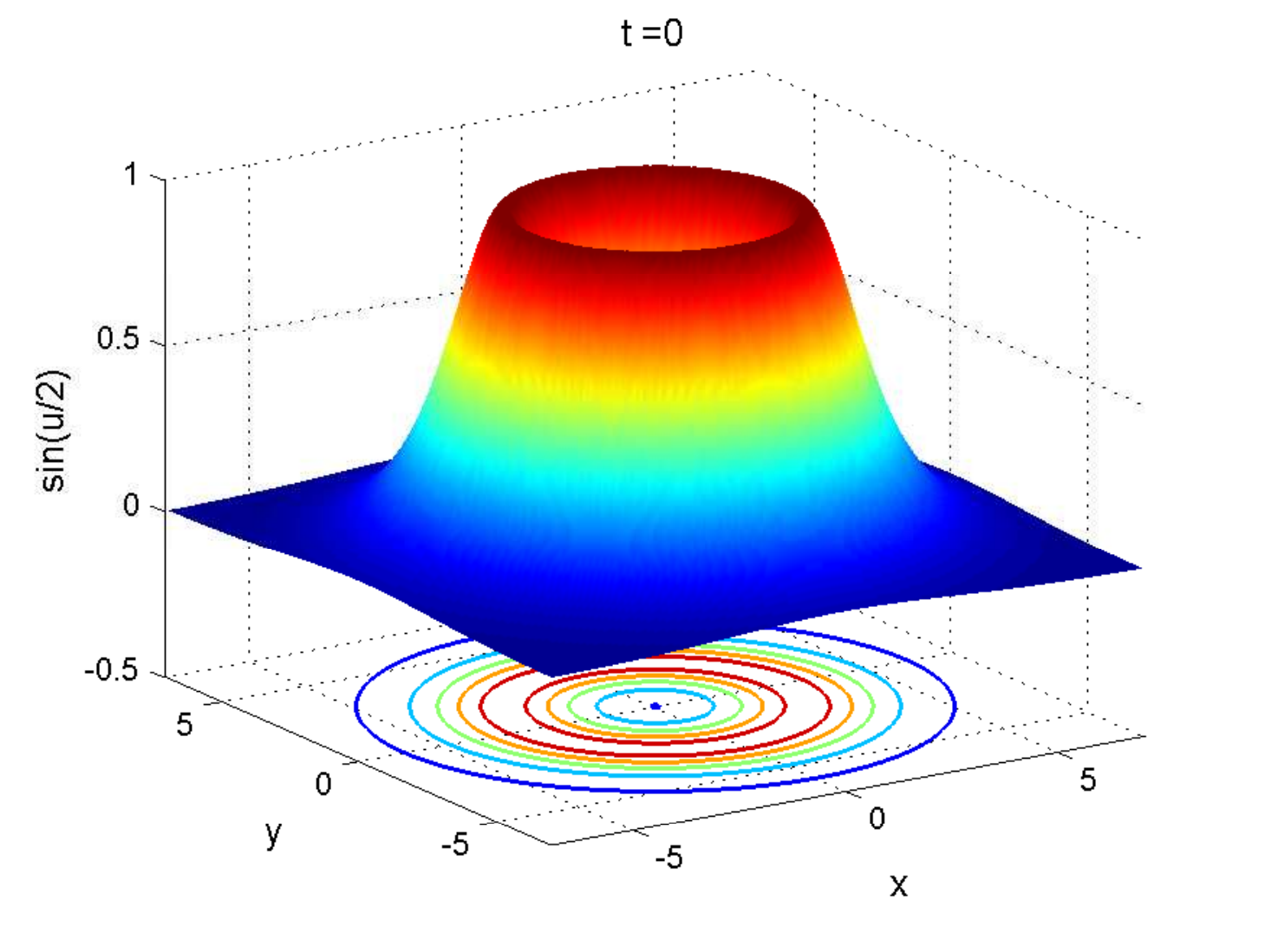}
	\includegraphics[width=0.32\linewidth]{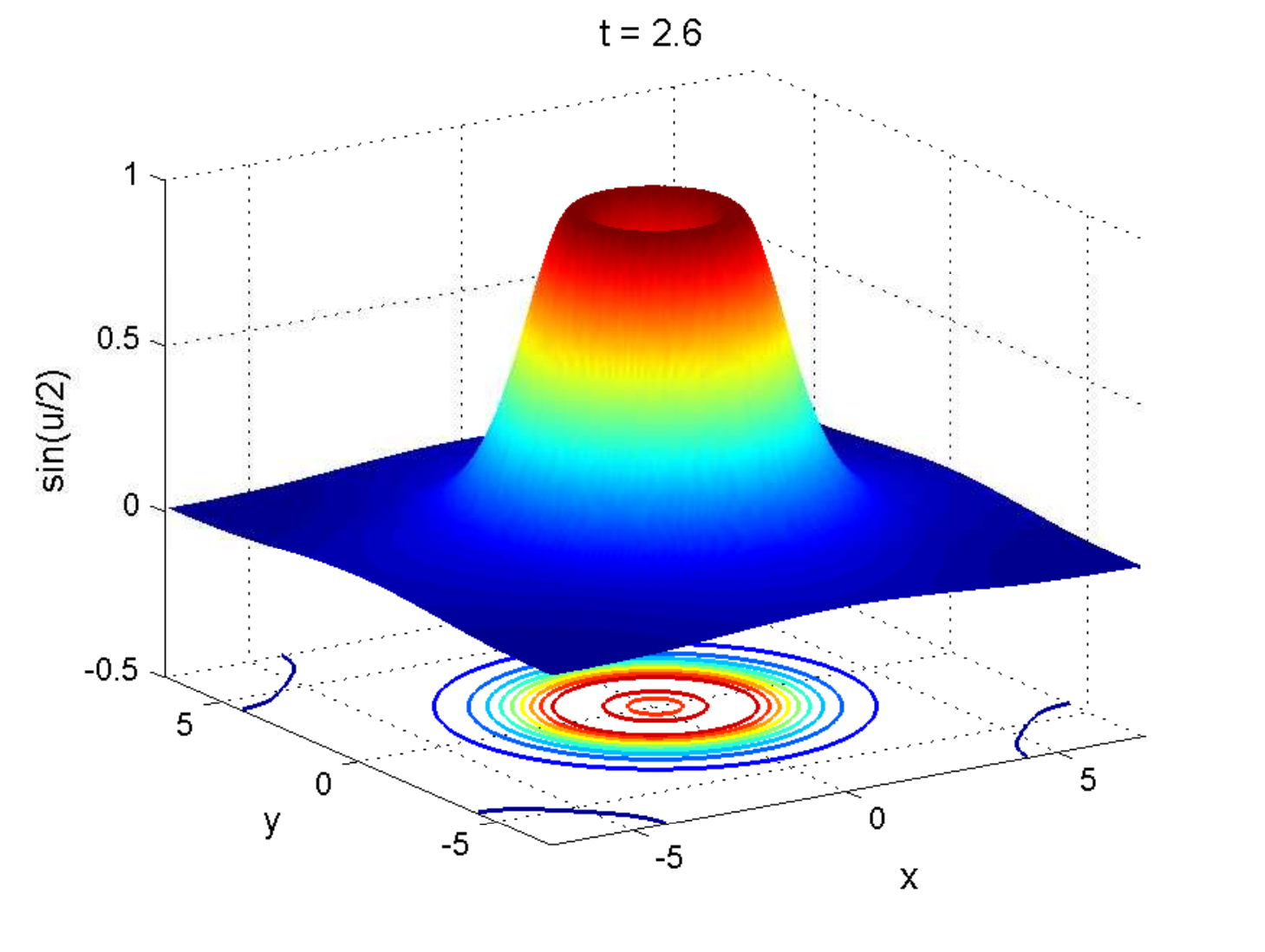}
	\includegraphics[width=0.32\linewidth]{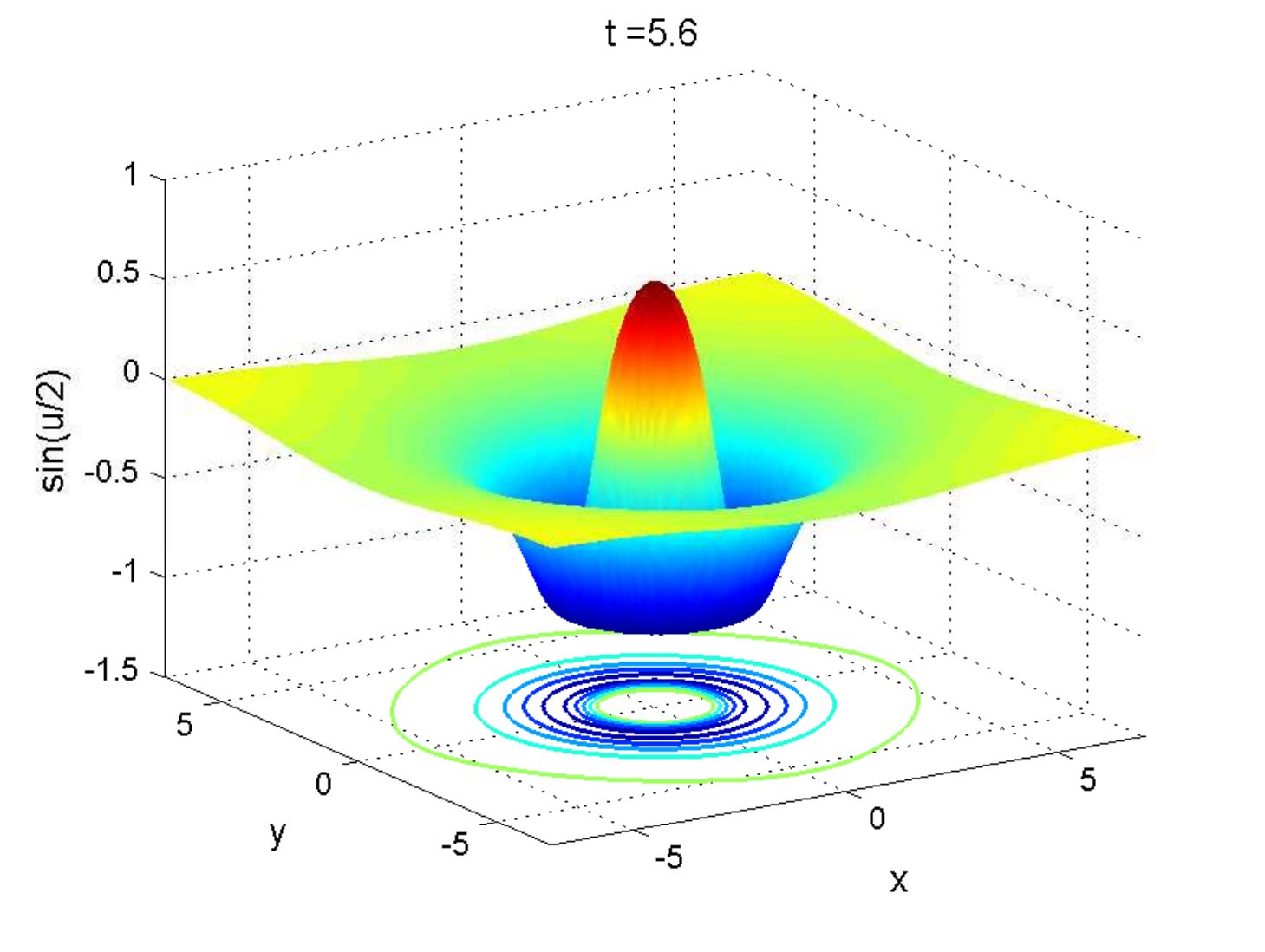}\\
	\includegraphics[width=0.32\linewidth]{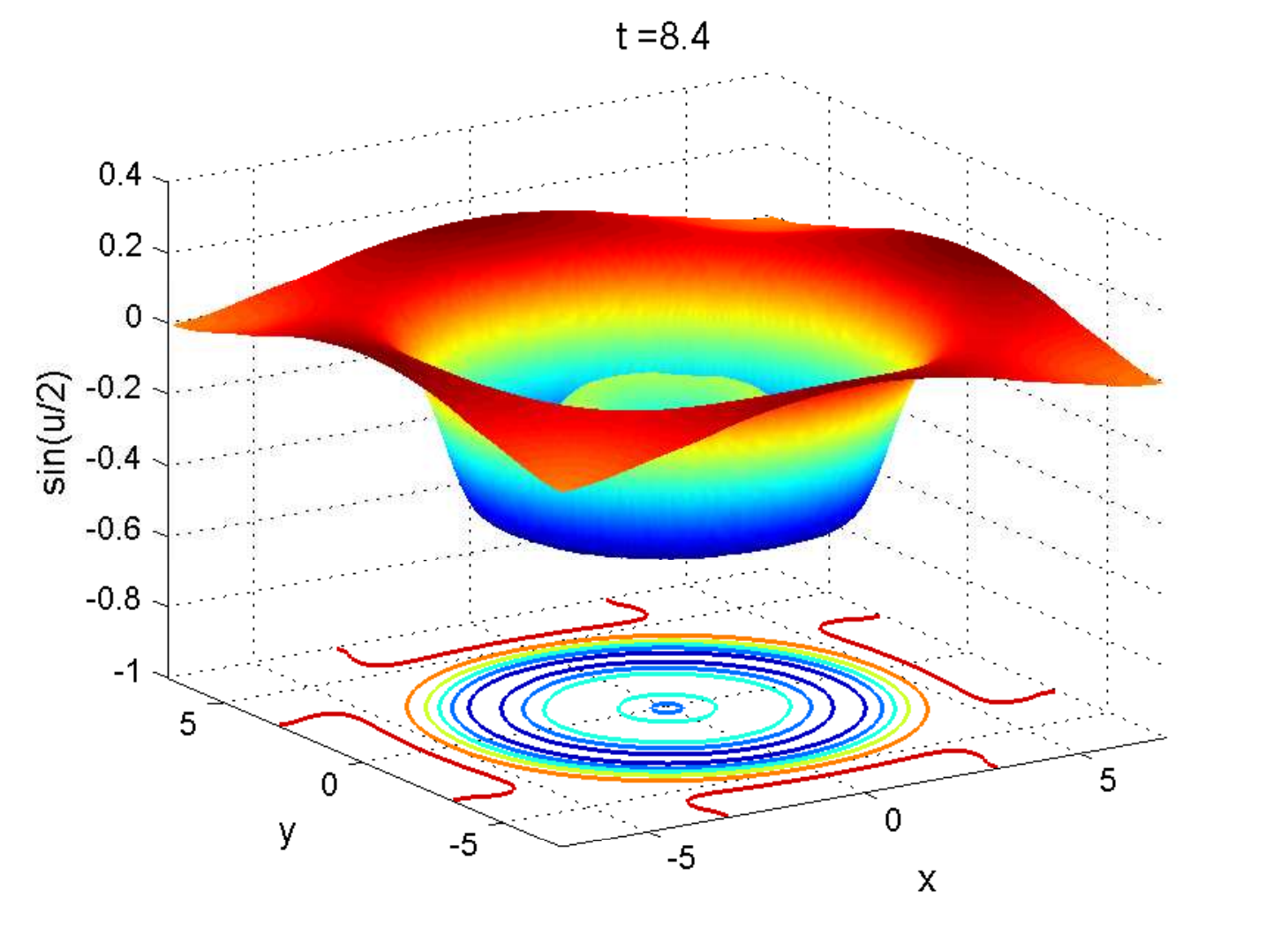}
	\includegraphics[width=0.32\linewidth]{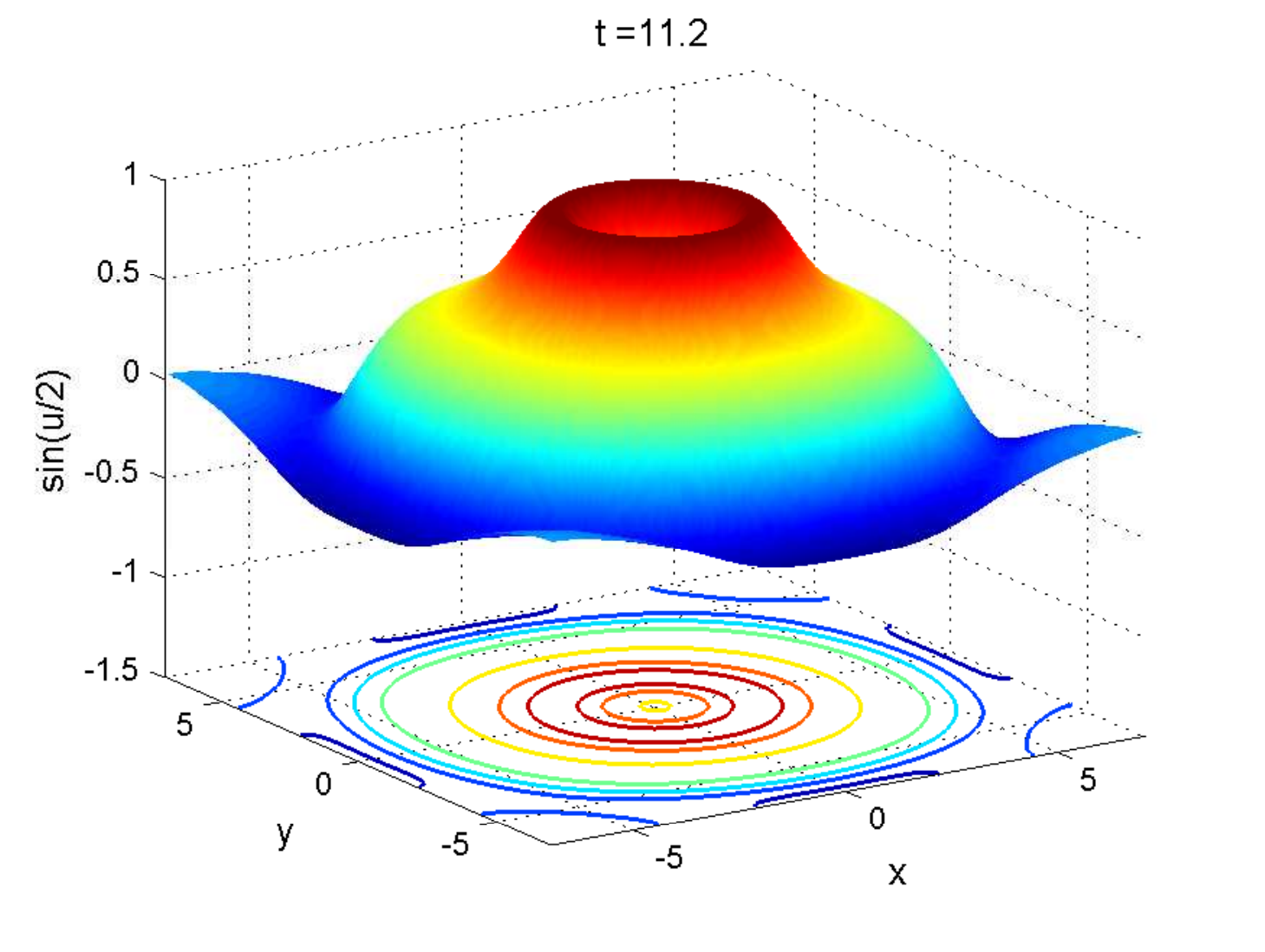}
	\includegraphics[width=0.32\linewidth]{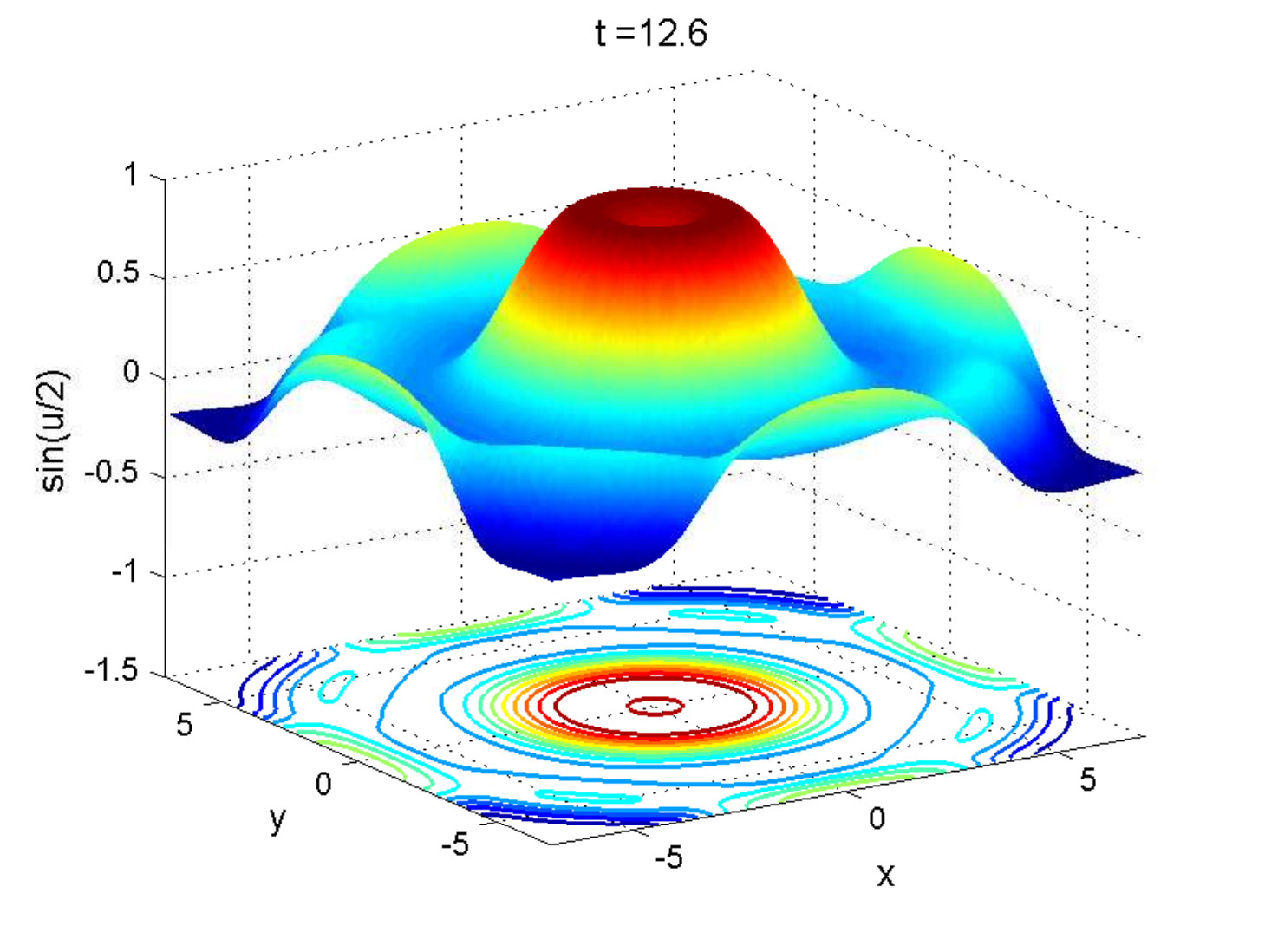}
	\caption{Circular ring soliton: surface and contours plots of initial condition and numerical solutions at $t = 2.8, 5.6,8.4,11.2, 12.6$, in terms of  $\sin(u/2)$. Spatial and temporal step sizes are taken as $h=0.5$, $\tau=0.01$.}
	\label{fig:ex4-t}
\end{figure}

Following \cite{cl81,kdt17} for this experiment it is considered $\phi(x,y)=1$ with initial conditions
\[
\begin{aligned}
&f_1(x,y) = 4\tan^{-1}(\exp(3-\sqrt{x^2+y^2})),\\
&f_2(x,y) = 0,~ -7\leq x,y\leq 7.
\end{aligned}
\]
The initial condition and numerical results are presented in Figures~\ref{fig:ex4-t}  at different times, in terms of $\sin(u/2)$. The soliton from its initial position, shrinks until $t=2.8$ appearing as a single-ring soliton. From $t=5.6$, which could be considered as the beginning of the expansion phase, an radiation appears, which is followed by oscillations at the boundaries. This expansion is continued until $t=11.2$, where the soliton is almost reformed. Finally, since $t=12.6$ it appears to be again in its shrinking phase. During all the above transformations no displacement of the center of
the soliton occurred, which can be clearly observed from the corresponding contour plots. Notice that the expansion strongly interacts with the homogeneous boundary condition, however, the energy deviations in Figure~\ref{fig:ex4-t} still behaves as expected in long time simulation, which again verifies the effectiveness of our proposed schemes to deal with the Neumann boundary conditions.

\begin{figure}[H]
	\centering
	\includegraphics[width=0.45\linewidth]{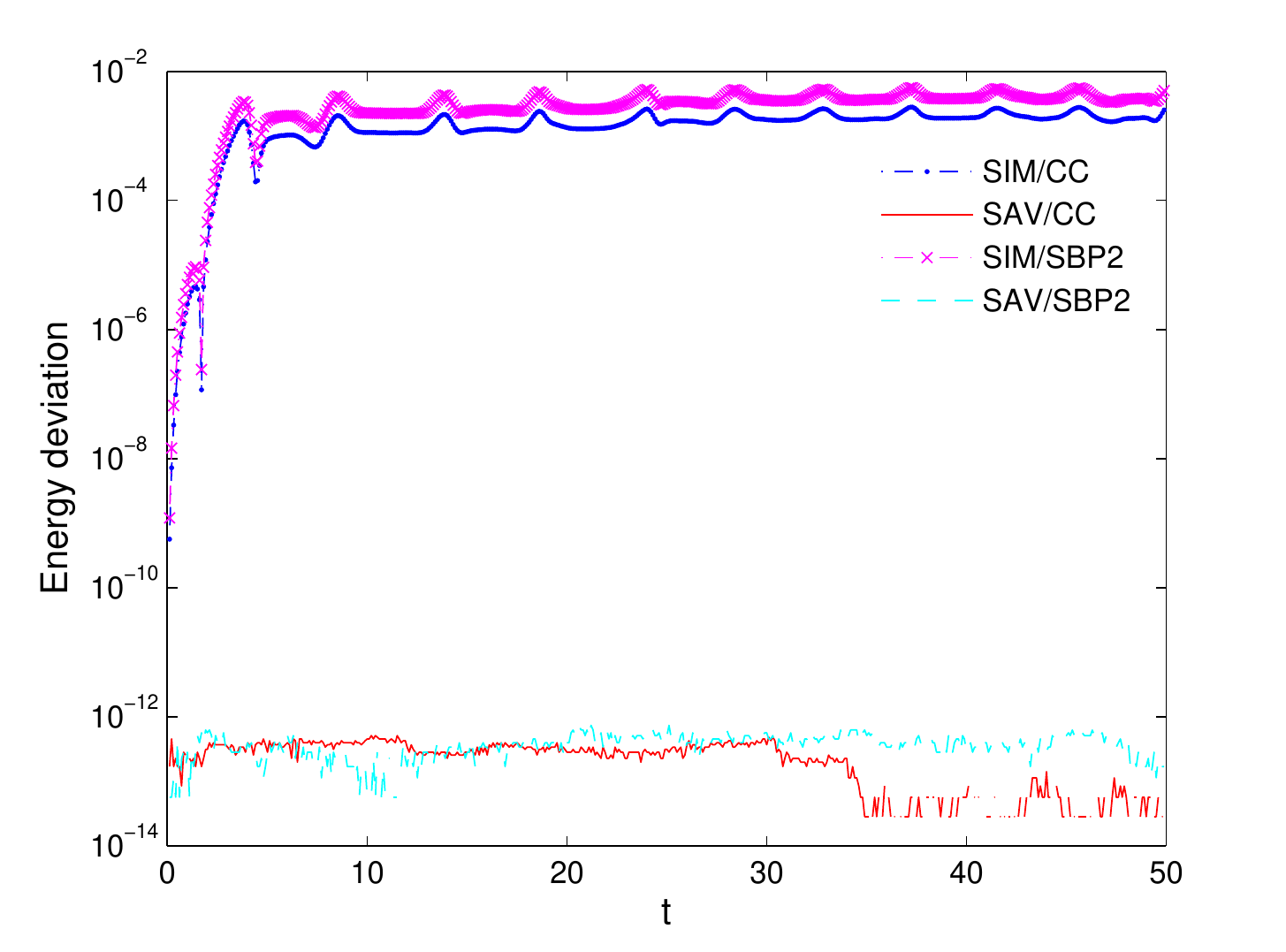}
	\caption{Long time energy deviation of the proposed schemes for the circular ring soliton.}
	\label{fig:ex4-ene}
\end{figure}

\subsection{Collision of two circular solitons}

\begin{figure}[H]
	\centering
	\includegraphics[width=0.32\linewidth]{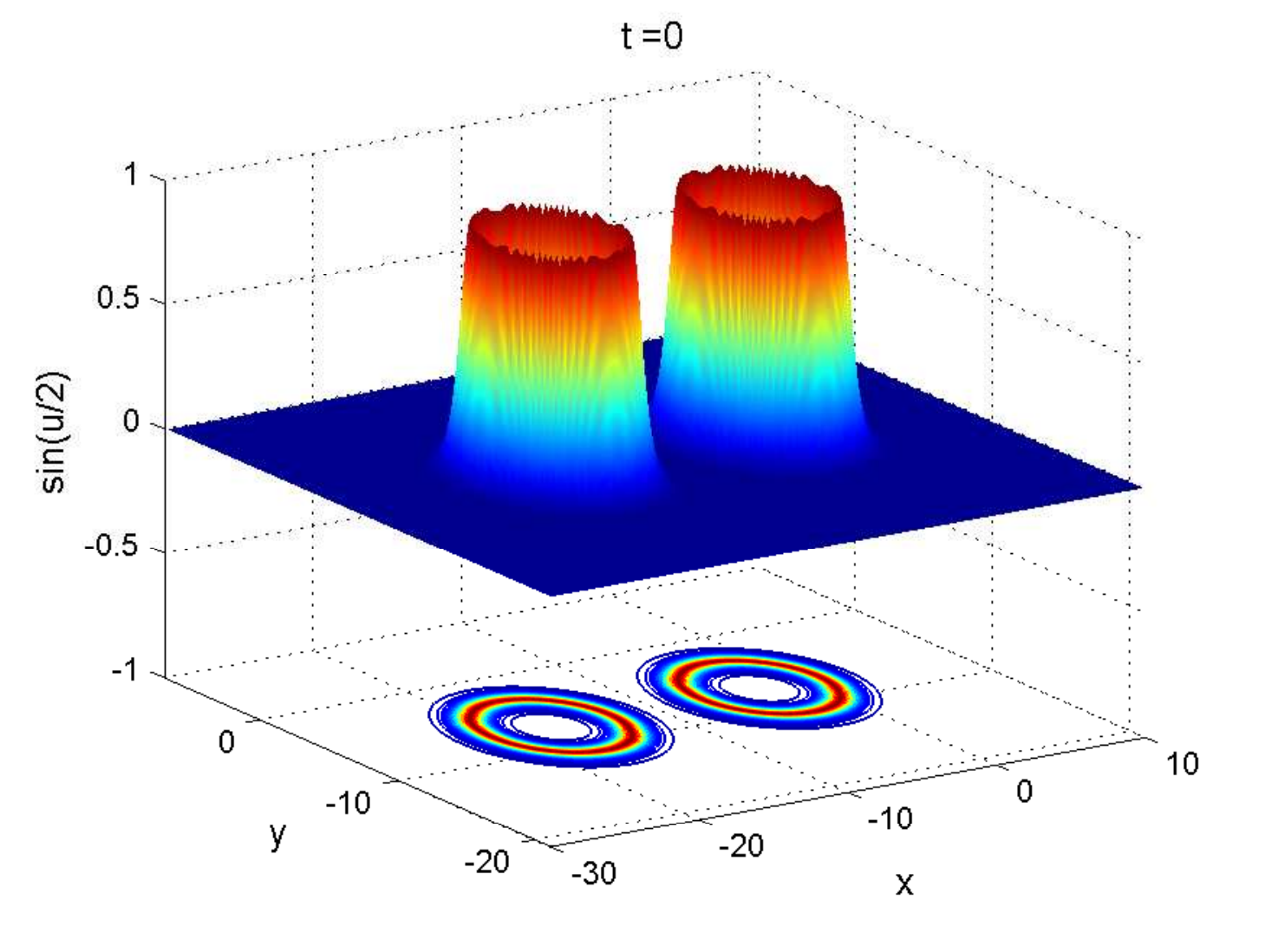}
	\includegraphics[width=0.32\linewidth]{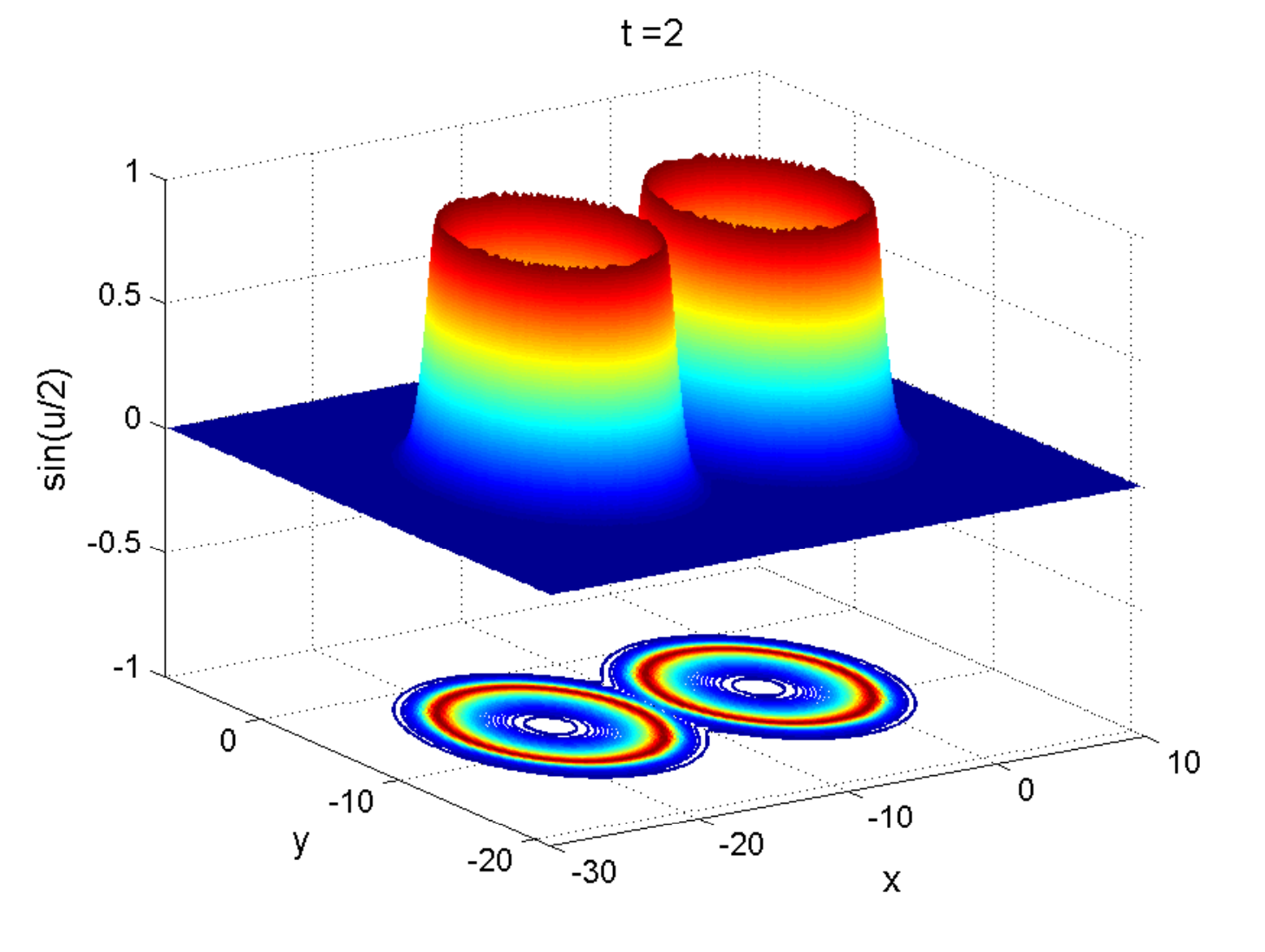}\\
	\includegraphics[width=0.32\linewidth]{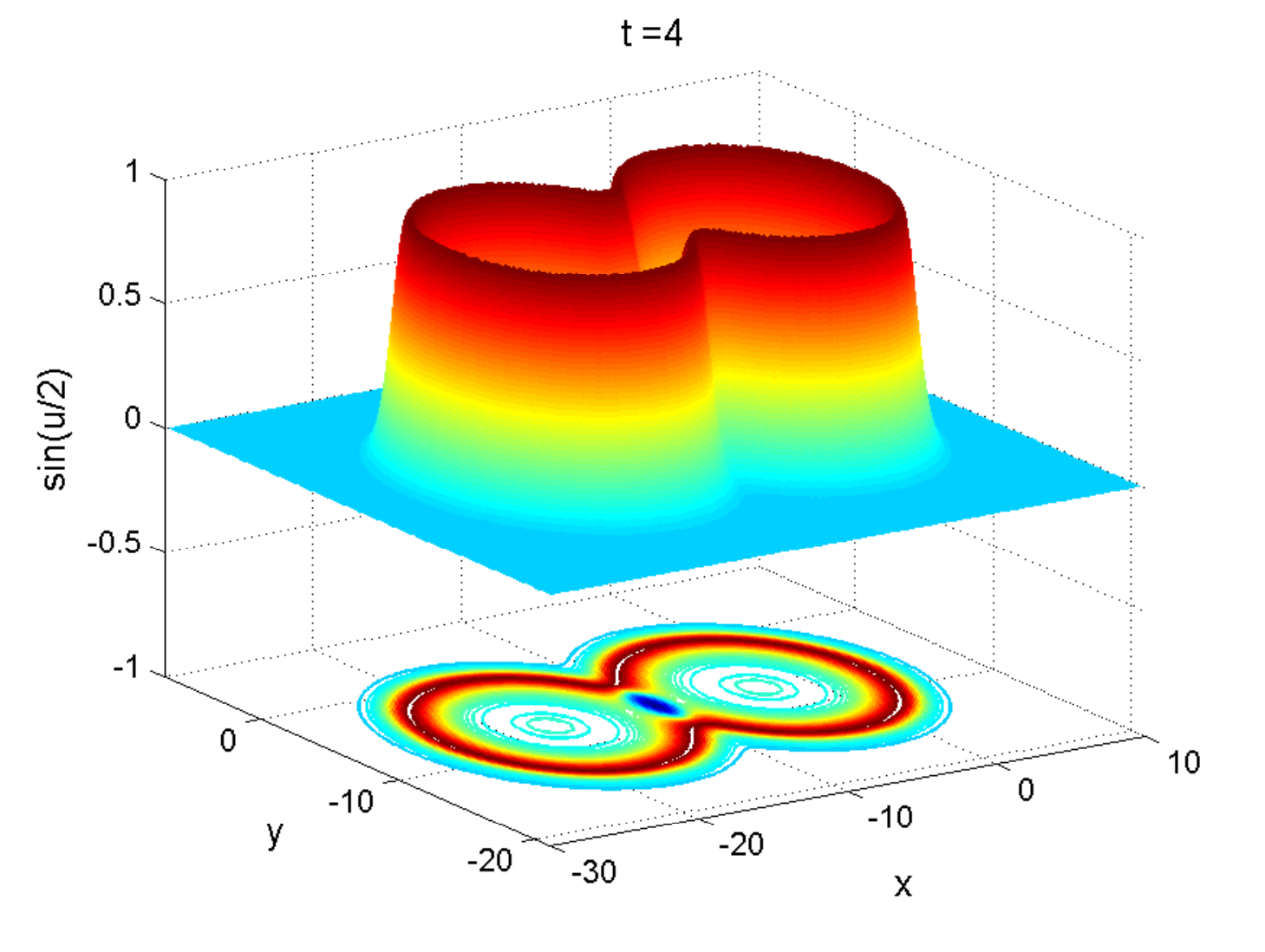}
	\includegraphics[width=0.32\linewidth]{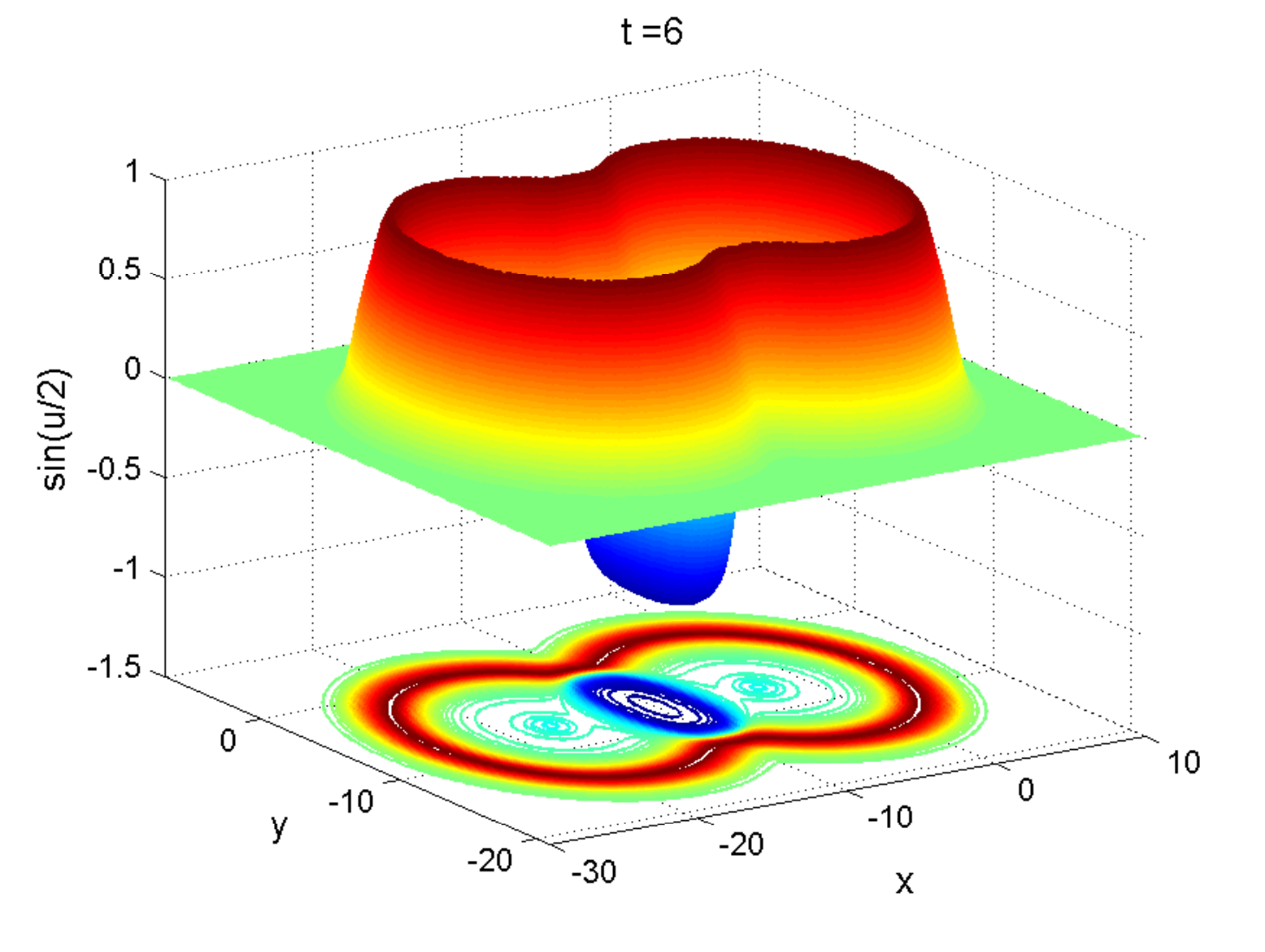}
	\includegraphics[width=0.32\linewidth]{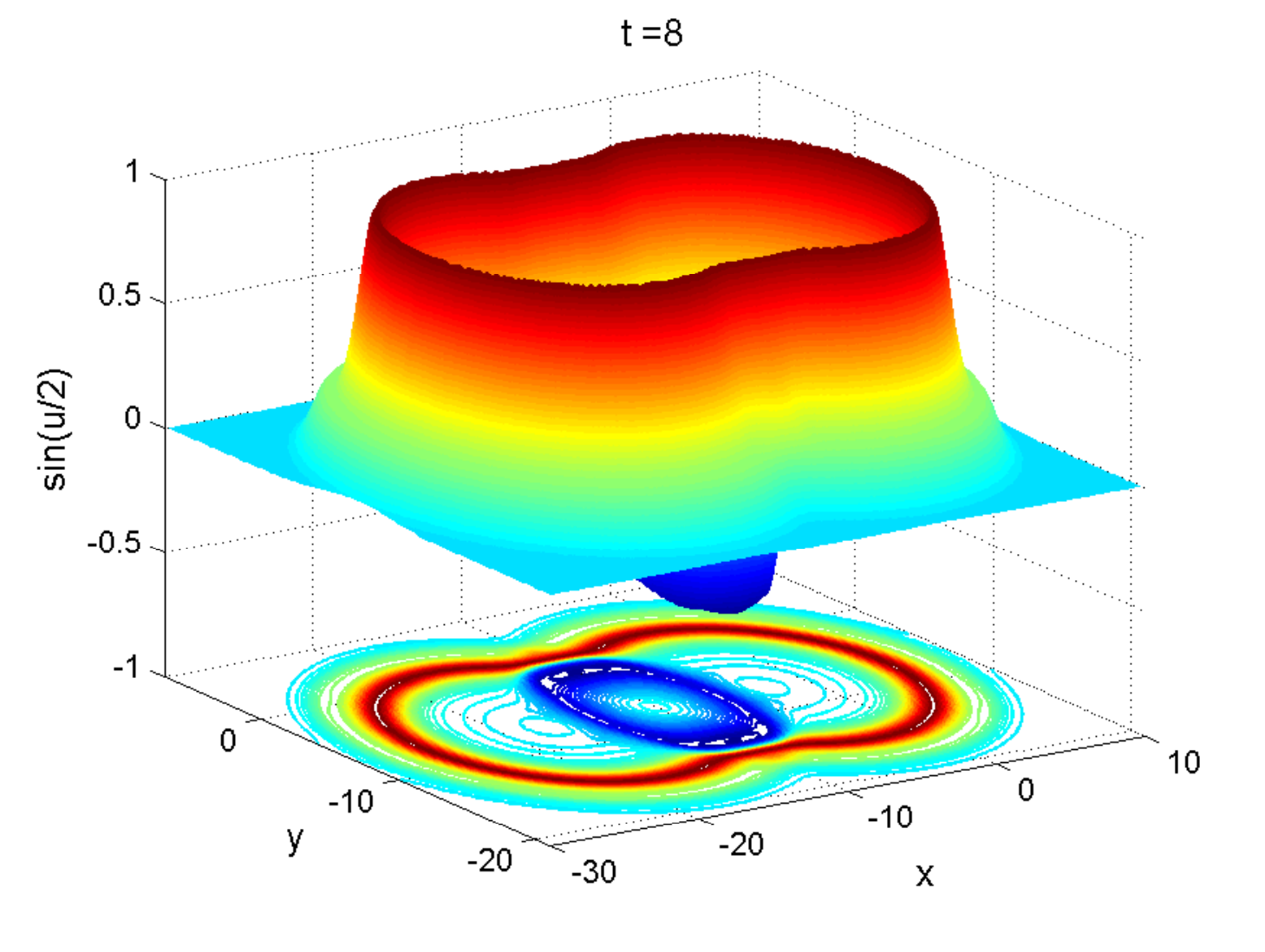}
	\caption{Collision of two expanding circular ring solitons: surface and contours plots of initial condition and numerical solutions at $t = 2, 4, 6, 8$, in terms of  $\sin(u/2)$. Spatial and temporal step sizes are taken as $h=0.5$, $\tau=0.01$.}
	\label{fig:ex5-t}
\end{figure}

\begin{figure}[H]
	\centering
	\includegraphics[width=0.45\linewidth]{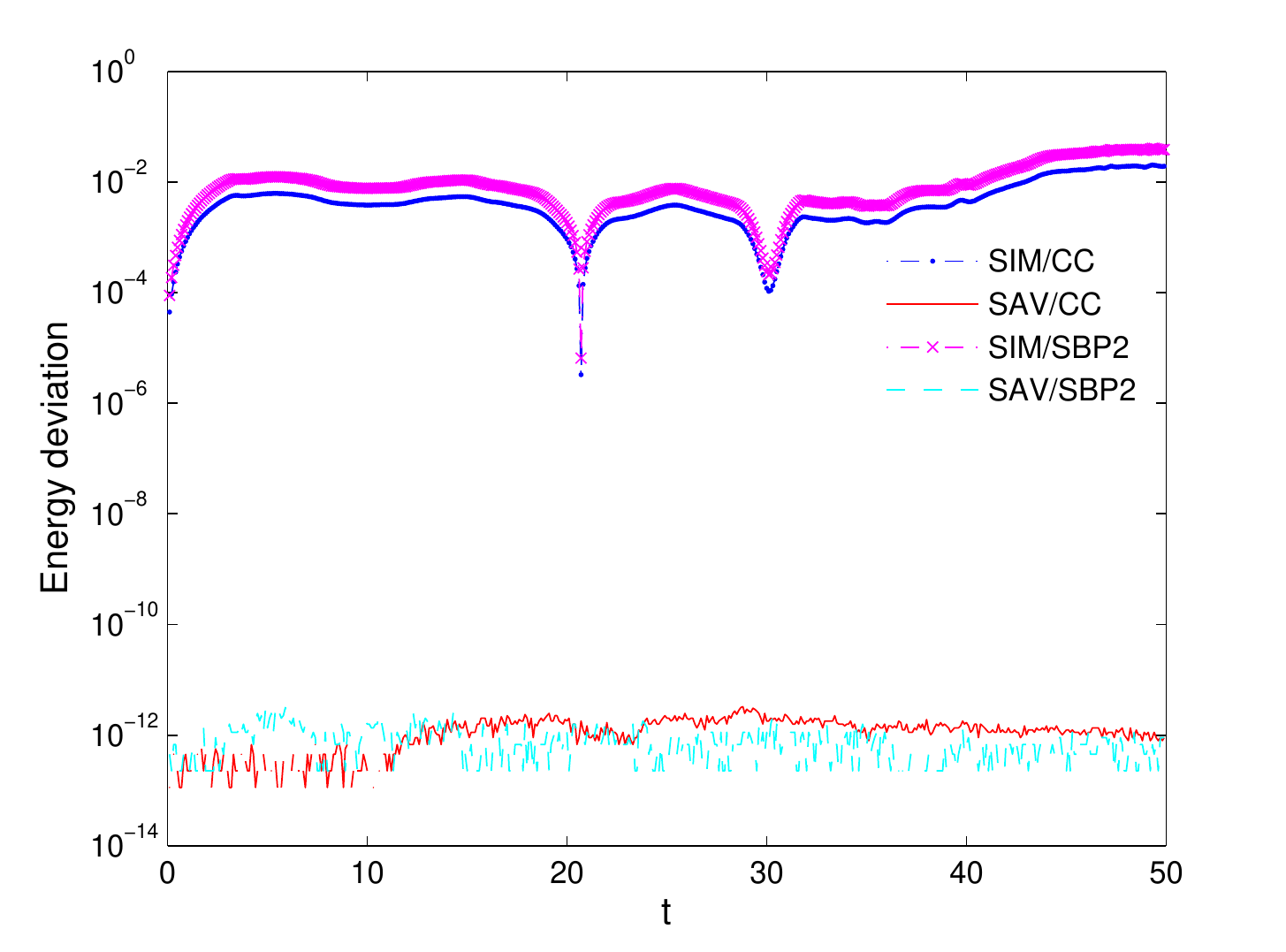}
	\caption{Long time energy deviation of the proposed schemes for the collision of two expanding circular ring solitons.}
	\label{fig:ex5-ene}
\end{figure}

The collision of two expanding circular ring solitons is considered with $\phi(x,y)=1$ and initial conditions
\[
\begin{aligned}
&f_1(x,y) = 4\;\tan^{-1}\exp\big[(4-\sqrt{(x+3)^2+(y+7)^2})/0.436\big],\\
&f_2(x,y) = 4.13\;\mbox{sech}\big[(4-\sqrt{(x+3)^2+(y+7)^2})/0.436\big],~ -30\leq x\leq 10,~-21\leq y\leq 7.
\end{aligned}
\]
The solution shown includes the extension across $x=-10$ and $y=-7$ by symmetry properties of the problem. Figure~\ref{fig:ex5-t} demonstrates the collision between two expanding circular ring solitons in which two smaller ring solitons bounding an annular region emerge into a large ring soliton. Contour maps are given to show more clearly the movement of solitons. The simulated solution is again precisely consistent to existing results \cite{ahh91,skv05}. The long time performance in Figure~\ref{fig:ex5-ene} also indicates a good conservation of energy which ensures the nonlinear stability computationally.

\subsection{Collision of four circular solitons}

\begin{figure}[H]
	\centering
	\includegraphics[width=0.32\linewidth]{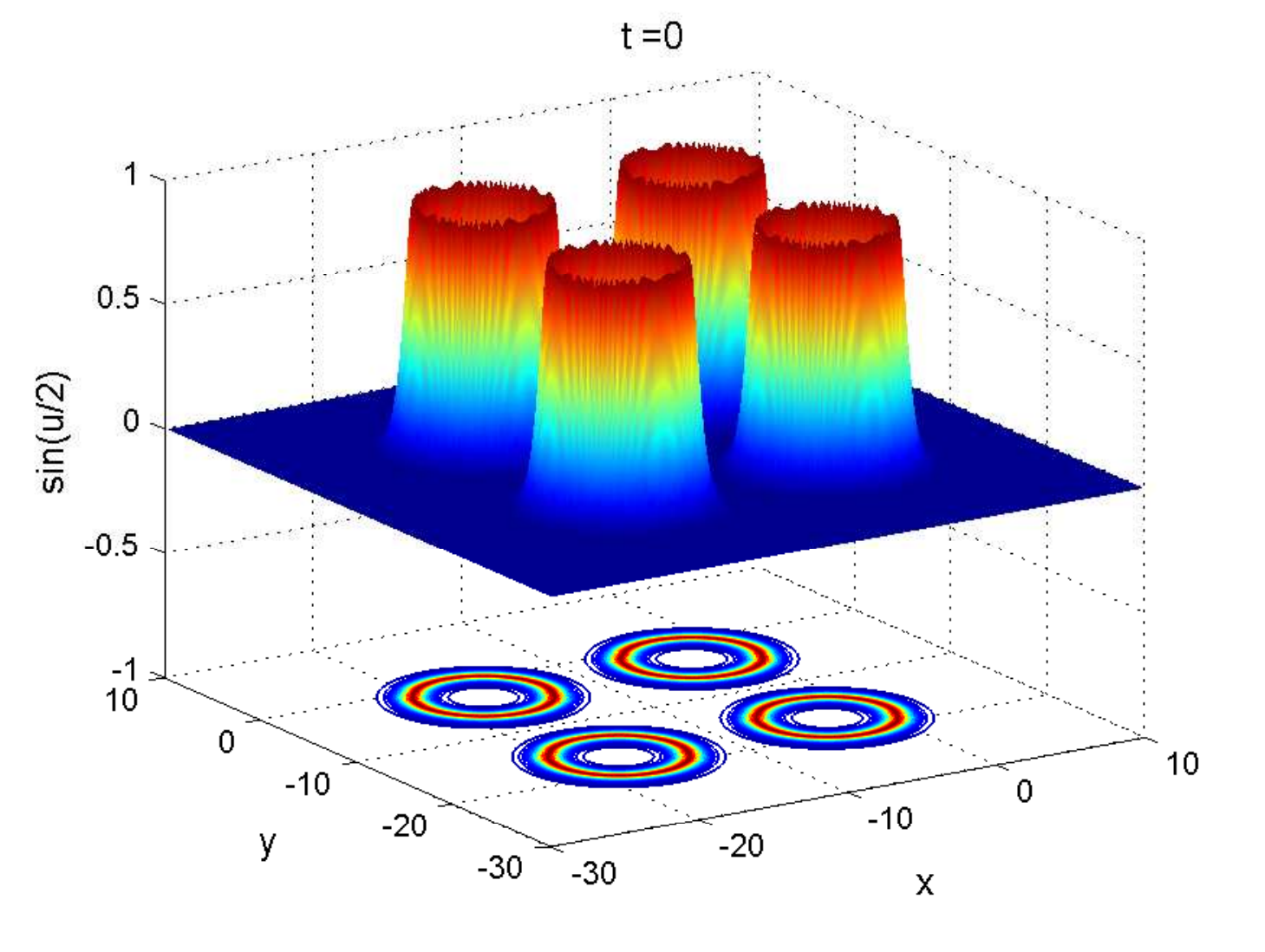}
	\includegraphics[width=0.32\linewidth]{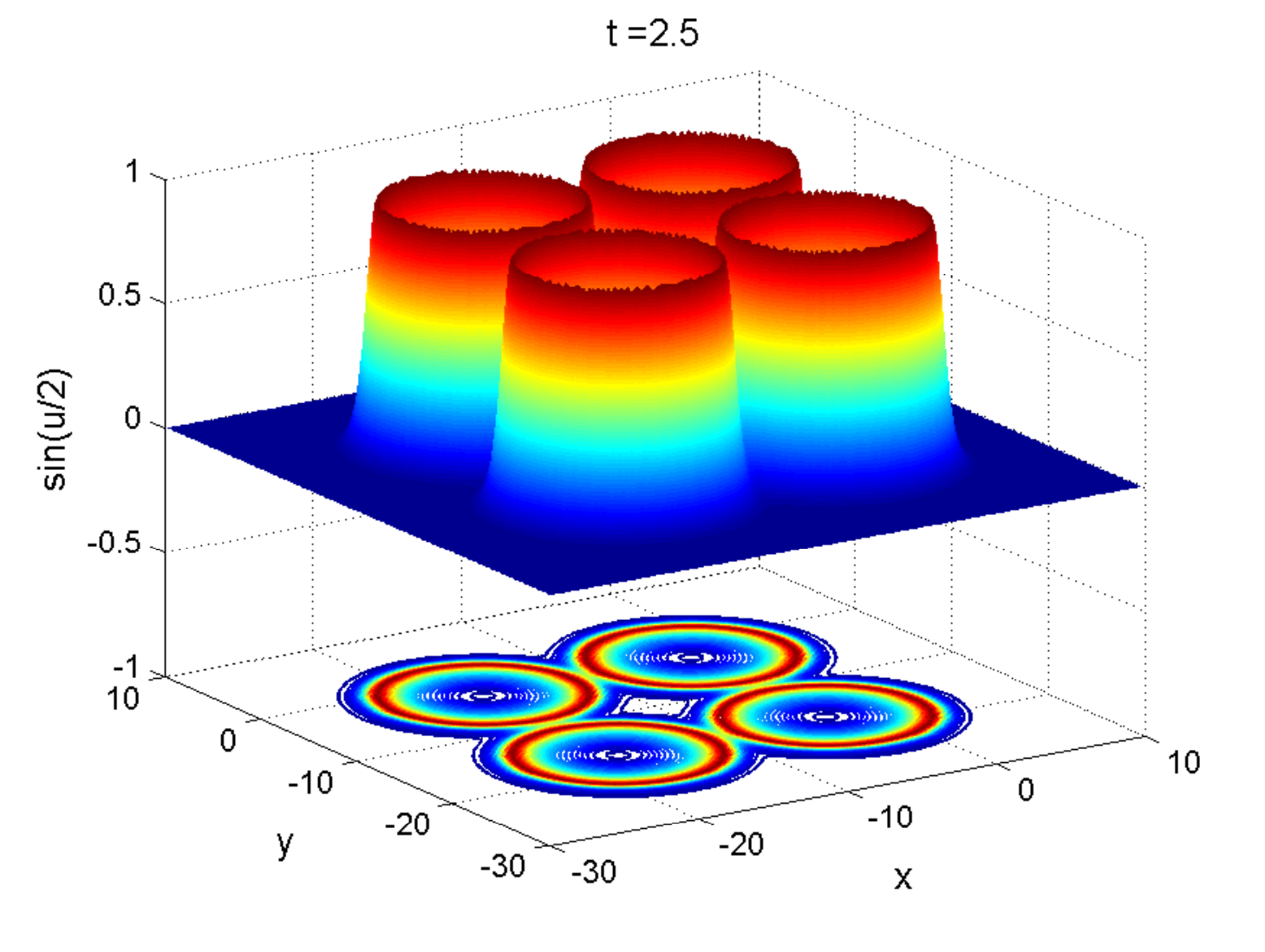}\\
	\includegraphics[width=0.32\linewidth]{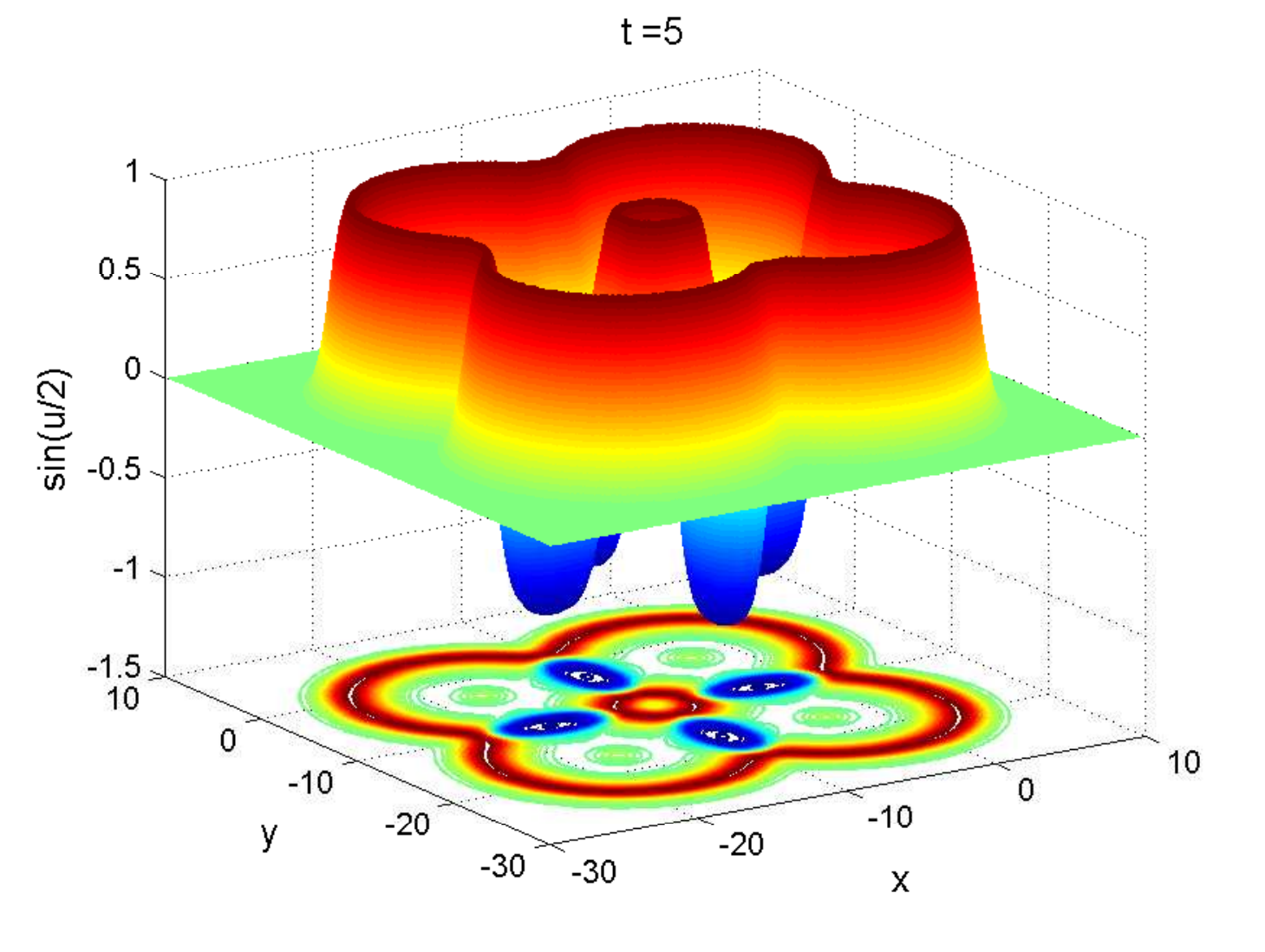}
	\includegraphics[width=0.32\linewidth]{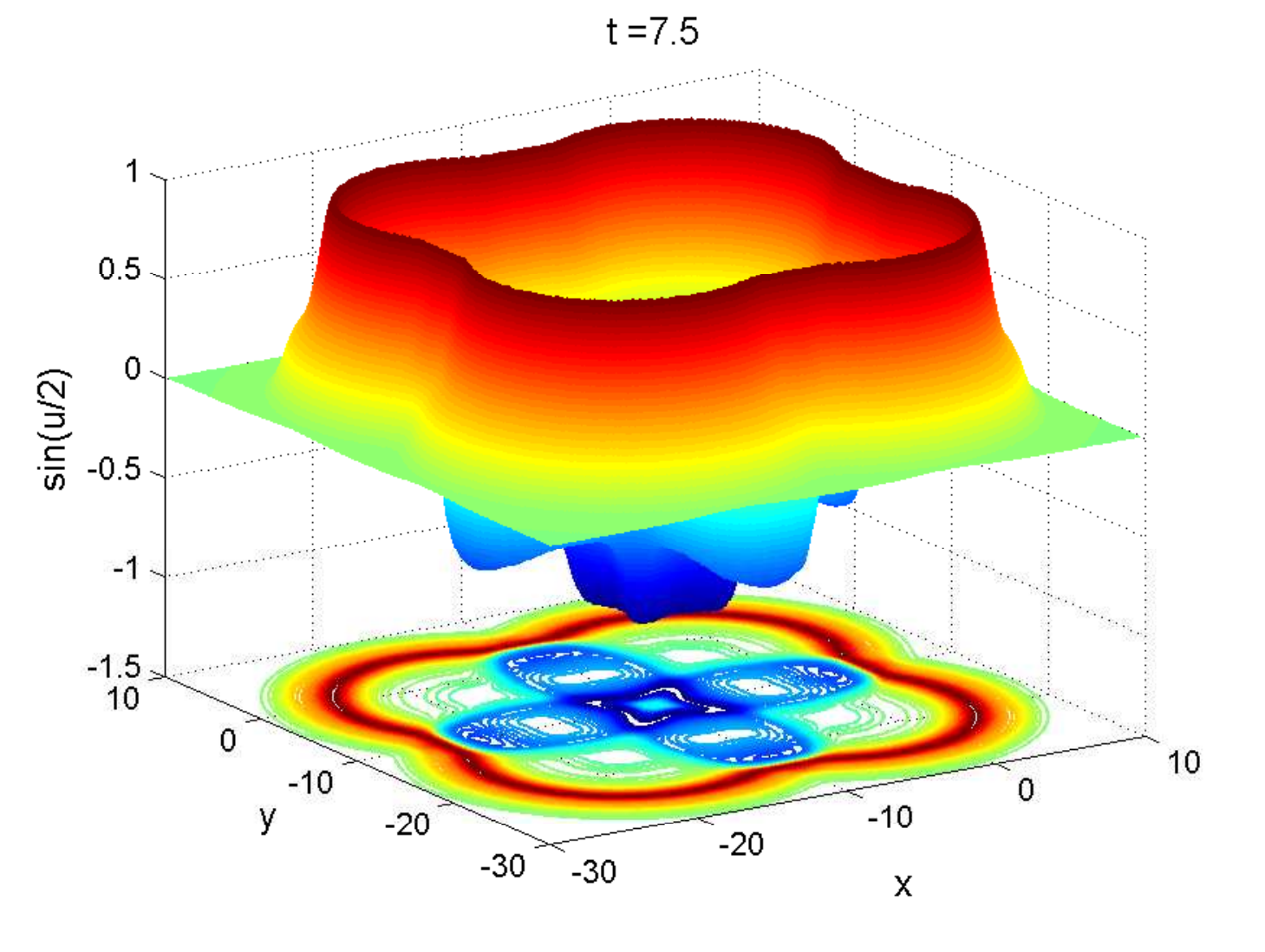}
	\includegraphics[width=0.32\linewidth]{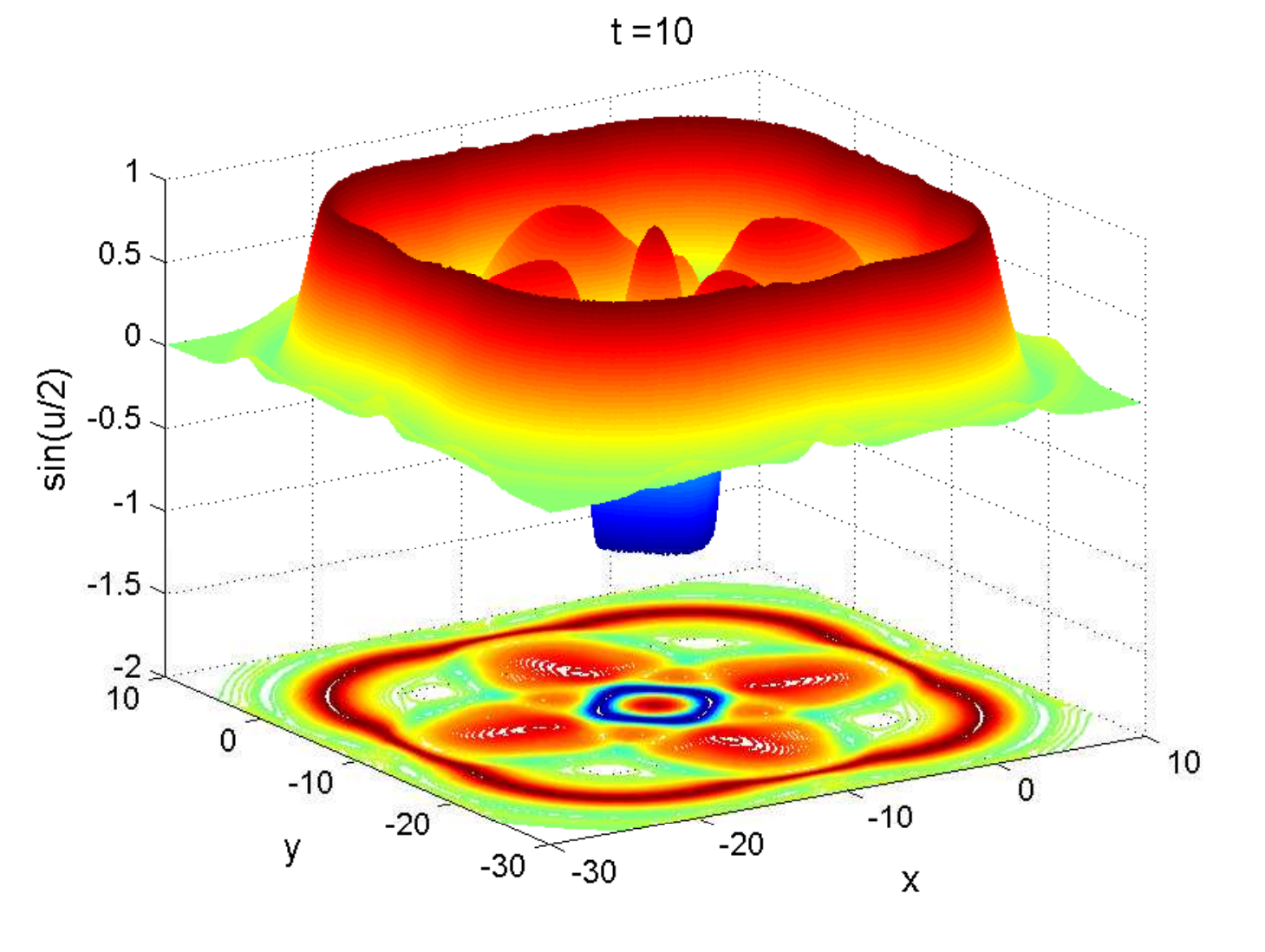}
	\caption{Collision of four expanding circular ring solitons: surface and contours plots of initial condition and numerical solutions at $t = 2.5, 5, 7.5, 10$, in terms of  $\sin(u/2)$. Spatial and temporal step sizes are taken as $h=0.5$, $\tau=0.01$.}
	\label{fig:ex6-t}
\end{figure}

\begin{figure}[H]
	\centering
	\includegraphics[width=0.45\linewidth]{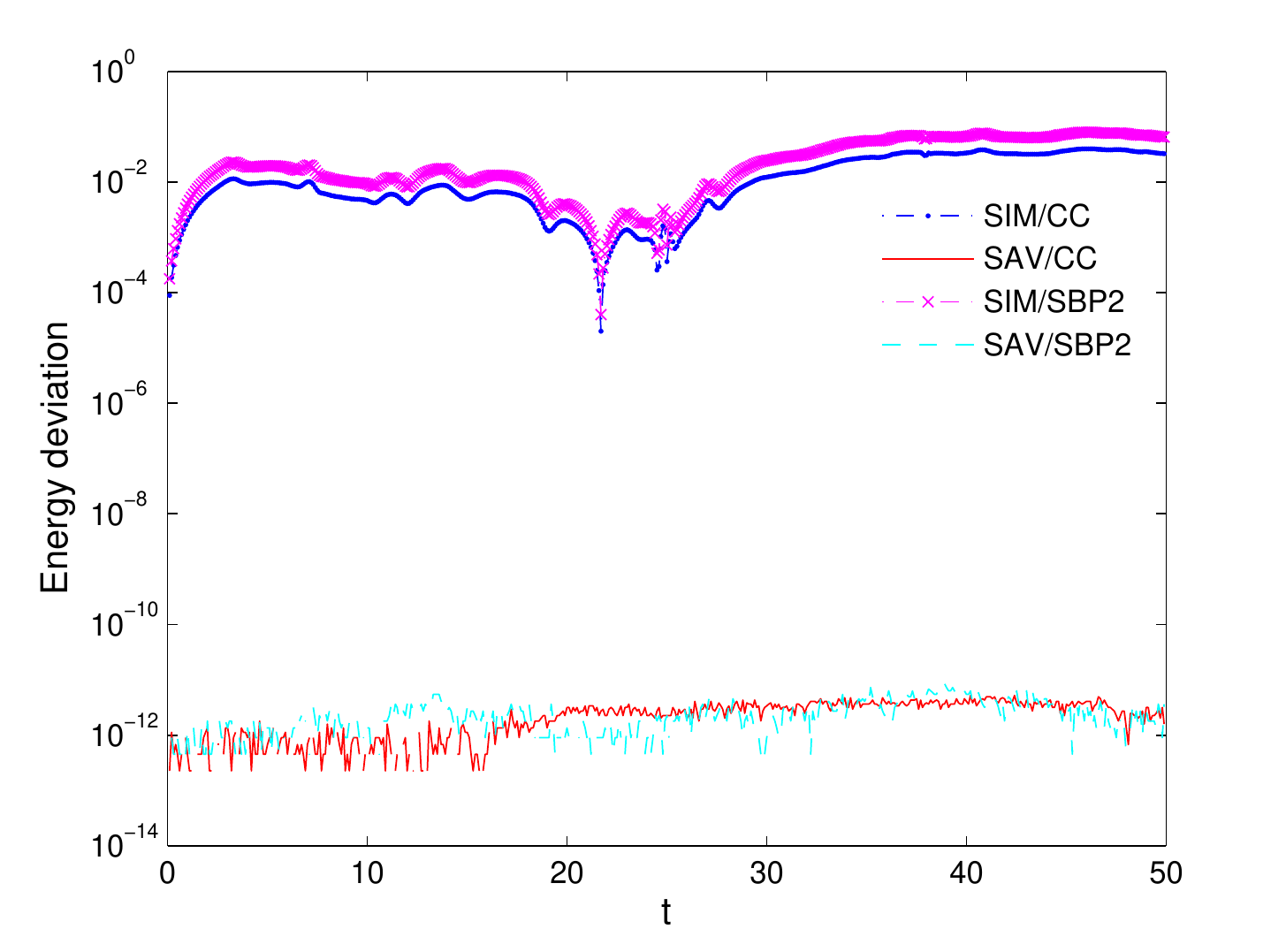}
	\caption{Long time energy deviation of the proposed schemes for the collision of four expanding circular ring solitons.}
	\label{fig:ex6-ene}
\end{figure}

Finally, for collisions of four expanding circular ring solitons, we take $\phi(x,y)=1$ and initial conditions
\[
\begin{aligned}
&f_1 = 4\;\mbox{atan}^{-1}\exp[(4-\sqrt{(x+3)^2+(y+3)^2})/0.436],\\
&f_2 = 4.13\;\mbox{sech}[(4-\sqrt{(x+3)^2+(y+3)^2})/0.436],~ -30\leq x,y\leq 10.
\end{aligned}
\]
The simulation is based on an extension across $x=−10$ and $y=−10$ due to the symmetry. Figure~\ref{fig:ex6-t} demonstrates precisely the collision between four expanding circular ring solitons in which the smaller ring solitons bounding an annular region emerge into a large ring soliton. Again, contour maps are given to illustrate more clearly the movement of the solitons. Those results show an extremely complex interaction with rapidly varying values of $u$ in the center and are in good agreement with the corresponding surfaces given in \cite{ahh91,skv05}. Furthermore, the energy deviations during the strong interaction maintain bounded to round-off errors for schemes {\bf SAV/CC} and {\bf SAV/SBP2}, and a small order of magnitude for schemes {\bf SIM/CC} and {\bf SIM/SBP2}, as demonstrated in Figure~\ref{fig:ex6-ene}.

\section{Conclusion}

In this paper, taking the two-dimensional sine-Gordon equation for example, we investigate the structure-preserving algorithms for problems with homogeneous Neumann boundary conditions. Two kinds of semi-discrete strategies are proposed, one is based on the cell-centered grid and standard central difference schemes for both spatial derivatives and boundary conditions, the other utilize the well-developed summation by parts operators. Both the resulting systems possess discrete energy conservation laws and Hamiltonian structures, but the equipped norms in the discrete energy function are distinguished. Also the relevant Hamiltonian structures exhibit canonical and non-canonical forms, respectively. Thereafter, symplectic integrators are obtained by applying the midpoint rule on time integration and linearly implicit energy-preserving schemes are constructed through the scalar auxiliary variable approach. Numerical experiments verify the convergence order as well as the effectiveness in dealing with homogeneous Neumann boundary condition, reflected by the long-time performance in the conservation of discrete energy.

\vspace{0.3cm}
\hspace{-0.5cm}{\bf Acknowledgements}\\
This work is supported by the National Key Research and Development Project of China (2016YFC0600310),  the National Natural Science Foundation of China (11771213, 41504078), the Natural Science Foundation of Jiangsu Province (BK20171480), the Jiangsu Collaborative Innovation Center for Climate Change and the Priority Academic Program Development of Jiangsu Higher Education Institutions.

\appendix
\renewcommand\thesection{\appendixname~\Alph{section}}
\renewcommand\theequation{\Alph{section}.\arabic{equation}}

\section{SBP operators with general Neumann boundary conditions}

For the implementation of general boundary conditions, the SBP finite difference operators are combined with the Simultaneous Approximation Term (SAT) method \cite{cga93}. Consider the two-dimensional sine-Gordon equation
\begin{equation}\label{eq-A-1}
u_{tt}=\Delta u-\phi(x,y)\sin u,
\end{equation}
on the domain $\Omega=[a,b]\times[c,d]$ with initial conditions \eqref{eq-1-2} and non-homogeneous Neumann boundary conditions
\begin{equation}\label{eq-A-3}
\begin{aligned}
&\frac{\partial  u}{\partial x}=g_1(x,y,t),\quad\mbox{for}~x=a ~\mbox{and}~ x=b, ~c\leq y\leq d,~ t>0,\\
&\frac{\partial  u}{\partial y}=g_2(x,y,t),\quad\mbox{for}~y=c ~\mbox{and}~ y=d, ~a\leq x\leq b,~ t>0.\\
\end{aligned}
\end{equation}
Taking the continuous inner product of \eqref{eq-A-1} with $u_t$ yields
\begin{equation}\label{eq-A-1-1}
\begin{aligned}
(u_{tt},u_t)=&(\Delta u,u_t)-(\phi\sin u,u_t)\\
=&-(\nabla u,\nabla u_t)-(\phi\sin u,u_t)+\int_{c}^{d}(u_xu_t)\big|_{x=a}^{x=b}dy+\int_{a}^{b}(u_yu_t)\big|_{y=c}^{y=d}dx,
\end{aligned}
\end{equation}
which is equivalent to the following energy identity
\begin{equation}\label{eq-A-2}
\frac{1}{2}\frac{d}{dt}\Big(\|u_t\|^2+\|\nabla u\|^2+2(\phi(1-\cos u),1)\Big)=\int_{c}^{d}(u_xu_t)\big|_{x=a}^{x=b}dy+\int_{a}^{b}(u_yu_t)\big|_{y=c}^{y=d}dx.
\end{equation}

The SBP operators, together with the SAT method, is designed to preserve a discrete analog to \eqref{eq-A-2}. The general approximation of \eqref{eq-A-1} using the SBP operator and the SAT method for the boundary is
\begin{equation}\label{eq-6-1}
\begin{aligned}
U_{tt}=&\Lambda_x^{-1}(D_x^2+BS_x)U+U(D_y^2+BS_y)^T\Lambda_y^{-1}-\Phi\cdot\sin U\\
&+\sigma_1\Lambda_x^{-1}(E_{0_x}S_xU-G_1^a)+\sigma_2\Lambda_x^{-1}(E_{N_x}S_xU-G_1^b)\\
&+\sigma_3(US_y^TE_{0_y}-G_2^c)\Lambda_y^{-1}+\sigma_4(US_y^TE_{N_y}-G_2^d)\Lambda_y^{-1},
\end{aligned}
\end{equation}
where $\Lambda_x^{-1}(D_x^2+BS_x)$ and $(D_y^2+BS_y)^T\Lambda_y^{-1}$ represent the SBP approximations of the second-order derivatives with $\Lambda_\alpha=\Lambda_\alpha^T>0$, $D_\alpha^2=(D_\alpha^2)^T\geq0$, $\alpha=x$ or $y$. $S_x$ and $S_y$ include approximations of first derivatives at boundaries. The last four penalty-like terms correspond to the SAT method by imposing boundary conditions weakly, where the relevant matrices are defined by
\[
\begin{aligned}
&B_\alpha=\mbox{diag}([-1,\underbrace{0,\cdots,0}_{N_\alpha-1},1]),\\
&E_{0_\alpha}=\mbox{diag}([1,\underbrace{0,\cdots,0}_{N_\alpha}]),\quad E_{N_\alpha}=\mbox{diag}([\underbrace{0,\cdots,0}_{N_\alpha},1]),\quad\alpha=x~ \mbox{or}~ y.
\end{aligned}
\]
The four parameters in these terms $\sigma_i$, $i=1,2,3,4$ are determined to achieve a discrete energy estimate. The inhomogeneous boundary conditions are reflected by the corresponding matrices with $N_x\times N_y$ dimensions
{\small
	\[
	G_1^a=\left(\begin{array}{ccccc}
	g_1(a,y_0)& g_1(a,y_1) & \cdots & g_1(a,y_{N_y}) \\ 
	0 & 0 & \cdots & 0 \\ 
	\vdots& \vdots &  &  \vdots \\ 
	0& 0 & \cdots & 0  \\ 
	\end{array} \right),\quad G_1^b=\left(\begin{array}{ccccc}
	0 & 0 & \cdots & 0 \\ 
	\vdots& \vdots &  &  \vdots \\ 
	0& 0 & \cdots & 0  \\ 
	g_1(b,y_0)& g_1(b,y_1) & \cdots & g_1(b,y_{N_y}) \\ 
	\end{array} \right),
	\]
	\[
	G_2^c=\left(\begin{array}{ccccc}
	g_2(x_0,c)& 0 & \cdots & 0 \\ 
	g_2(x_1,c)& 0 & \cdots & 0 \\ 
	\vdots& \vdots &  &  \vdots \\ 
	g_2(x_{N_x},c)& 0 & \cdots & 0  \\ 
	\end{array} \right),\quad G_2^d=\left(\begin{array}{ccccc}
	0 & \cdots & 0 & g_2(x_0,d)\\ 
	0 & \cdots & 0 & g_2(x_1.d)\\ 
	\vdots &  &  \vdots & \vdots\\ 
	0 & \cdots & 0 & g_2(x_{N_x},d) \\ 
	\end{array} \right).
	\]
}

In Section 3, we have introduce a simplest SBP operators with corresponding $\Lambda_\alpha$, $D_\alpha^2$ given by \eqref{eq-3-3} and \eqref{lam} which result schemes of second-order accuracy. For high-order SBP operators approximating the second derivative, Mattsson and Nordstr\"om  \cite{mn04} proposed a systematic means of constructing SBP operators and gave exact expressions for the finite difference coefficients. Therefore, it is straightforward to utilize the existing high-order SBP operators to construct structure-preserving algorithms with Neumann boundary conditions.

To illustrate the stability of the semi-discretization \eqref{eq-6-1},  we take the weighted inner product \eqref{eq-3-0} of \eqref{eq-6-1} with $U_t$, and obtain
\begin{equation}\label{eq-6-1-1}
\begin{aligned}
(U_{tt},U_t)_\Lambda=&\big((D_x^2+B_xS_x)U\Lambda_y,U_t\big)_h+\big(\Lambda_xU(D_y^2+B_yS_y)^T,U_t\big)_h-(\Phi\cdot\sin U,U_t)_\Lambda\\
&+\sigma_1\big((E_{0_x}D_xU-G_1^a)\Lambda_y,U_t\big)_h+\sigma_2\big((E_{N_x}D_xU-G_1^b)\Lambda_y,U_t\big)_h\\
&+\sigma_3\big(\Lambda_x(UD_y^TE_{0_y}-G_2^c),U_t\big)_h+\sigma_4\big(\Lambda_x(UD_y^TE_{N_y}-G_2^d),U_t\big)_h\\
=&\big(D_x^2U\Lambda_y+(B_x+\sigma_1E_{0_x}+\sigma_2E_{N_x})D_xU\Lambda_y-(\sigma_1G_1^a+\sigma_2G_1^b)\Lambda_y,U_t\big)_h\\
&+\big(\Lambda_xUD_y^2+\Lambda_xUD_y^T(B_y+\sigma_3E_{0_y}+\sigma_4E_{N_y})-\Lambda_x(\sigma_3G_2^c+\sigma_4G_2^d),U_t\big)_h.
\end{aligned}
\end{equation}
A choice of $\sigma_1=\sigma_3=1$ and $\sigma_2=\sigma_4=-1$ can eliminate the terms containing $D_x$ and $D_y$ and yields a discrete analog to the original initial boundary problem
\begin{equation}\label{eq-6-1-2}
\begin{aligned}
(U_{tt},U_t)_\Lambda
=&(D_x^2U\Lambda_y+\Lambda_xUD_y^2,U_t)_h-(\Phi\cdot\sin U,U_t)_\Lambda-\big((G_1^a-G_1^b)\Lambda_y+\Lambda_x(G_2^c-G_2^d),U_t\big)_h\\
=&\big(\Lambda_x^{-1}D_x^2U+UD_y^2\Lambda_y^{-1},U_t\big)_\Lambda-(\Phi\cdot\sin U,U_t)_\Lambda-\big(\Lambda_x^{-1}(G_1^a-G_1^b)+(G_2^c-G_2^d)\Lambda_y^{-1},U_t\big)_\Lambda,
\end{aligned}
\end{equation}
which can be further rearranged as
\begin{equation}
\frac{1}{2}\frac{d}{dt}\Big(\|U_t\|_\Lambda^2+\|\nabla_h U\|_\Lambda^2+2(\Phi\cdot(1-\cos U),1)\Big)=\big(\Lambda_x^{-1}(-G_1^a+G_1^b)+(-G_2^c+G_2^d)\Lambda_y^{-1},U_t\big)_\Lambda,
\end{equation}
where a similar identity to \eqref{eq-3-4} can be obtained. Notice that when the boundary conditions are homogeneous, we can recover a same semi-discrete energy conservation law as presented in Section 3.

With the choices of parameters $\sigma_i$, $i=1,2,3,4$, the scheme \eqref{eq-6-1} with arbitrary high-order SBP operators can be simplified as
\begin{equation}\label{eq-6-4}
U_{tt}=\Lambda_x^{-1}D_x^2U+UD_y^2\Lambda_y^{-1}-\Phi\cdot\sin U+\Lambda_x^{-1}(-G_1^a+G_2^b)+(-G_2^c+G_2^d)\Lambda_y^{-1}.
\end{equation}
In the following, we further present the associated matrices of a fourth-order accuracy SBP operator with
\[
\Lambda_\alpha=\left(\begin{array}{ccccccp{0.5cm}}
\frac{17}{48}&  &  &  &  & & \\ 
& \frac{59}{48} &  &  &  & &\\ 
&  & \frac{43}{48} &  &  & &\\ 
&  &  & \frac{49}{48} &  & &\\ 
&  &  &  & 1 &  &\\ 
&  &  &  &  & \ddots &
\end{array} \right),
\]
and
\[
D_\alpha^2=\frac{1}{h_\alpha^2}\left(\begin{array}{cccccccc}
\frac{9}{8} & -\frac{59}{48} & \frac{1}{12} & \frac{1}{48} & 0 & 0 & 0 &\\[1ex]
-\frac{59}{48} & \frac{59}{24} & -\frac{59}{48} & 0 & 0  & 0 & 0&\\ [1ex]
\frac{1}{12} & -\frac{59}{48} & \frac{55}{24} & -\frac{59}{48} & \frac{1}{12} & 0 & 0 &\\ [1ex]
\frac{1}{48} & 0 & -\frac{59}{48} & \frac{59}{24} & -\frac{4}{3} & \frac{1}{12} & 0 &\\ [1ex]
0 & 0 & \frac{1}{12} & -\frac{4}{3} & \frac{5}{2} & -\frac{4}{3} & \frac{1}{12} &\\ [1ex]
&  &  & \ddots & \ddots & \ddots & \ddots & \ddots
\end{array} \right),
\]
where both matrices are bisymmetric and have corresponding entries in the lower right-hand corner. Expressions of higher-order SBP operators for the second derivative up to eighth order accuracy can be directly found in \cite{mn04}.



\end{document}